\documentclass[noinfoline,11pt]{imsart}
\usepackage{amstext}
\usepackage{amssymb}

\usepackage[top = 1 in, bottom = 1in, left = 1in, right=1in]{geometry}
\usepackage[numbers]{natbib}
\usepackage{hyperref}
\hypersetup{colorlinks,allcolors=blue}
\usepackage[section]{placeins}

\RequirePackage[OT1]{fontenc}
\RequirePackage{amsthm,amsmath}
\RequirePackage{natbib}

\usepackage{booktabs}
\usepackage{graphicx}
\usepackage{caption}
\usepackage{subcaption}
\usepackage[table,xcdraw]{xcolor}
\usepackage{relsize}

\startlocaldefs
\numberwithin{equation}{section}
\theoremstyle{plain}
\newtheorem{thm}{Theorem}[section]
\endlocaldefs
\newtheorem{prop}{Proposition}[section]
\newtheorem{lem}{Lemma}[section]
\newtheorem{remark}{Remark}[section]

\newcommand{\y}{\mathbf{y}}
\newcommand{\X}{\mathbf{X}}
\newcommand{\bb}{\boldsymbol{\beta}}
\newcommand{\ep}{\boldsymbol{\varepsilon}}
\newcommand{\bl}{\boldsymbol{\lambda}}
\newcommand{\Bl}{\boldsymbol{\Lambda}}
\newcommand{\E}{\boldsymbol{E}}

\newcommand\Myciteauthor[1]{\citeauthor{#1} \cite{#1}}
\newcommand{\lmDelta}{\Delta}
\newcommand{\lmdelta}{\delta}
\newcommand{\lmLambda}{\Lambda}
\newcommand{\lmlambda}{\lambda}
\newcommand{\lmX}{X}

\newcommand{\lmy}{y}
\newcommand{\lmU}{U}
\newcommand{\lmu}{u}
\newcommand{\lmV}{V}

\newcommand{\lmD}{D}
\newcommand{\lmd}{d}
\newcommand{\R}{\mathbb{R}}
\newcommand{\norm}[1]{ \Vert #1 \Vert }
\newcommand{\lme}{e}
\newcommand{\lmA}{A}
\newcommand{\lma}{a}

\newcommand{\lmxbf}{{\bf x}}
\newcommand{\lmxv}{{\bf v}}

\begin{document}
  
\begin{frontmatter}

\title{Geometric ergodicity of Gibbs samplers for the Horseshoe and its regularized variants}
\runtitle{Geometric ergodicity for the Horseshoe}

\author[a]{ Suman K. Bhattacharya},
\author[a]{ Kshitij Khare}
\and 
\author[b]{ Subhadip Pal}
\address[a]{Department of Statistics, University of Florida}
\address[b]{Department of Bioinformatics and Biostatistics, University of Louisville}

\runauthor{S. Bhattacharya, K. Khare, and S. Pal}

\begin{abstract}
The Horseshoe is a widely used and popular continuous shrinkage prior 
for high-dimensional Bayesian linear regression. Recently, regularized versions of 
the Horseshoe prior have also been introduced in the literature. Various Gibbs 
sampling Markov chains have been developed in the literature to generate approximate 
samples from the corresponding intractable posterior densities. Establishing 
geometric ergodicity of these Markov chains provides crucial technical justification 
for the accuracy of asymptotic standard errors for Markov chain based estimates of 
posterior quantities. In this paper, we establish geometric ergodicity for various 
Gibbs samplers corresponding to the Horseshoe prior and its regularized variants in 
the context of linear regression. First, we establish geometric ergodicity 
of a Gibbs sampler for the original Horseshoe posterior under strictly weaker 
conditions than existing analyses in the literature. Second, we consider the 
regularized Horseshoe prior introduced in \cite{piironen2017}, and prove geometric 
ergodicity for a Gibbs sampling Markov chain to sample from the corresponding 
posterior without any truncation constraint on the global and local shrinkage 
parameters. Finally, we consider a variant of this regularized Horseshoe prior 
introduced in \cite{nishimura2019shrinkage}, and again establish geometric 
ergodicity for a Gibbs sampling Markov chain to sample from the corresponding 
posterior. 
\end{abstract}



\end{frontmatter}


\section{Introduction} \label{sec:intro}

\noindent
Consider the linear model $\y = \X \bb + \sigma \ep$, where $\y \in 
\mathbb{R}^n$ is the response vector, $\X$ is the $n \times p$ design matrix, $\bb \in \mathbb{R}^p$ is the vector of regression coefficients, $\ep$ is the error vector with i.i.d. standard normal components, and $\sigma^2$ is the error variance. The goal is to estimate the unknown parameters $(\bb, \sigma^2)$. In modern applications, datasets where the number of predictors $p$ is much larger than the sample size $n$ are commonly encountered. A standard approach for meaningful statistical estimation in these over-parametrized settings is to assume that only a few of the signals are prominent (the others are small/insignificant). This is mathematically formalized by assuming that the underlying regression coefficient vector is sparse. In the Bayesian paradigm, this assumption of sparsity is accommodated either by choosing spike-and-slab priors (mixture of point mass at zero and an absolutely continuous density) or absolutely continuous shrinkage priors which selectively shrink the small/insignificant signals. 

A variety of useful shrinkage priors have been proposed in the literature (see \cite{BPPD:2015, 10.2307/25734098, Polsonandscott2012} and the references therein), and the Horseshoe prior  (\cite{10.2307/25734098}) is a widely used and highly popular choice. The Horseshoe prior for linear regression is specified as follows. 
\begin{eqnarray}
\left.\bb \right. \mid \bl, \sigma^2,\tau^2\sim \mathcal{N}_p(0,\sigma^2\tau^2 \Bl) \nonumber\\
\left.\lambda_i\right.\sim C^+ (0,1) \text{ independently for}\ i=1,2,\cdots,p \nonumber\\
\tau^2\sim \pi_\tau (\cdot) \hspace{0.2in} \sigma^2\sim 
\text{Inverse-Gamma}(a,b) \label{horseshoe_prior}
\end{eqnarray}


\noindent
where $\mathcal{N}_d$ denotes the $d-$variate normal density, $\Bl$ is a diagonal matrix with diagonal entries given by the entries $\left\{\lambda_j^2\right\}_{j=1}^p$, and $\text{Inverse-Gamma}(a,b)$ denotes the Inverse-Gamma density with shape parameter $a$ and rate parameter $b$. The vector $\bl = (\lambda_j^2)_{j=1}^p$ is referred to as the vector of local (component-wise) shrinkage parameters, while $\tau^2$ is referred to as the global shrinkage parameter. 

The resulting posterior distribution for $(\bb, \sigma^2)$ is intractable in the sense that closed form computations or i.i.d. sampling from this distribution are not feasible. Several Gibbs sampling Markov chains have been proposed in the literature to generate approximate samples from the Horseshoe posterior, see for example (\cite{10.1093/biomet/asw042, f521496bc3f649b2a2d7eda0124887a5, johndrow2017bayes, Makalic_2016}). 

The fact that parameter values which are far away from zero are not regularized at all due to the heavy tails is considered to be a key strength of the Horseshoe prior. However, as pointed out in \Myciteauthor{piironen2017}, this can be undesirable when the parameters are only weakly identified. To address this issue, \cite{piironen2017} introduced the regularized Horseshoe prior, given by 
\begin{eqnarray*}
\left.\beta_i \right. \mid \bl, \sigma^2,\tau^2\sim \mathcal{N}_p \left( 0, \left( \frac{1}{c^2} + \frac{1}{\lambda_i^2 \tau^2} \right)^{-1} \sigma^2 \right) \text{ independently for}\ i=1,2,\cdots,p \\
\left.\lambda_i\right.\sim C^+ (0,1) \text{independently for}\ i=1,2,\cdots,p \\
\tau^2\sim \pi_\tau (\cdot) \hspace{0.2in} \sigma^2\sim 
\text{Inverse-Gamma}(a,b)
\end{eqnarray*}

\noindent
Here $c$ is a finite constant which controls additional regularization of all regression parameters (large and small). The original Horseshoe prior can be recovered by letting $c \rightarrow \infty$. \Myciteauthor{piironen2017} use a Hamiltonian Monte Carlo (HMC) based approach to generate approximate samples from the corresponding regularized Horseshoe posterior distribution. Also, any Gibbs sampler 
for the Horseshoe posterior can be suitably adapted in the regularized setting. 

For any practitioner using Markov chain Monte Carlo, it is crucial to understand the accuracy of the resulting MCMC based estimates by obtaining valid standard errors for these estimates. The notion of geometric ergodicity plays an important role in this endeavor, as explained below. Let $(\bb_m, \sigma^2_m)_{m \geq 0}$ denote a Harris ergodic Markov chain with the Horseshoe or regularized Horseshoe posterior density, denoted by $\pi_H (\cdot \mid \y)$, as its stationary density. The Markov chain is said to be geometrically ergodic if 
$$
\left\| K_{\bb_0, \sigma^2_0}^m - \Pi_H \right\|_{\mbox{\tiny{TV}}} \leq C\left(\bb_0, \sigma^2_0\right) \gamma^m 
$$

\noindent
where $K_{\bb_0, \sigma^2_0}^m$ denotes the distribution of the Markov chain started at $(\bb_0, \sigma^2_0)$ after $m$ steps, $\Pi_H$ denotes the stationary distribution, and $\|\cdot\|_{\mbox{\tiny{TV}}}$ denotes the total variation norm. Suppose we wish to evaluate the posterior expectation 
$$
E_{\pi_H (\cdot \mid \y)} g = \int \int g\left(\bb, \sigma^2\right) \pi_H \left(\bb, \sigma^2 
\mid \y\right) d \bb d \sigma^2 
$$

\noindent
for a real-valued measurable function $g$ of interest. Harris ergodicity guarantees that the Markov chain based estimator 
$$
\bar{g}_m := \frac{1}{m+1} \sum_{i=0}^m g\left(\bb_i, \sigma^2_i\right) 
$$

\noindent
is strongly consistent for $E_{\pi (\cdot \mid \y)} g$. An estimate by itself, however, is not quite useful without an associated standard error. All known methods to compute consistent estimates (see for example \cite{flegal2010}, \cite{10.2307/25734098}) of the standard error for $\bar{g}_m$ require the existence of a Markov chain CLT which establishes 
$$
\sqrt{m} \left( \bar{g}_m - E_{\pi_H (\cdot \mid \y)} g \right) \rightarrow \mathcal{N}(0, \sigma_g^2), 
$$

\noindent
for $\sigma_g^2 \in (0, \infty)$. In turn, the standard approach for 
establishing a Markov chain CLT requires proving geometric ergodicity of the underlying Markov chain. To summarize, proving geometric ergodicity helps 
rigorously establish the asymptotic validity of
CLT based standard error estimates used by MCMC practitioners. 

Establishing geometric ergodicity for continuous state space Markov chains encountered in most statistical applications is in general a very challenging task. For a significant majority of Markov chains in statistical applications, the question of whether they are geometrically ergodic or not has not been resolved, although there have been some success stories. In the context of Markov chains arising in Bayesian shrinkage, geometric ergodicity of Gibbs samplers corresponding to various shrinkage priors such as the Bayesian lasso, Normal-Gamma, Dirichlet-Laplace and double Pareto priors has been recently established in (\cite{khare2013, 
RePEc:bla:scjsta:v:44:y:2017:i:2:p:307-323, pal2014}). 

Results for the Horseshoe prior remained elusive until very recently. The marginal Horseshoe prior on entries of $\bb$ (integrating out $\bl$, given $\tau^2$) has an infinite spike near zero and significantly heavier tails than the shrinkage priors mentioned above. This structure, while making it very attractive for sparsity selection, implicitly creates a lot of complications and challenges in the geometric ergodicity analysis using drift and minorization techniques. Recently, the authors 
in \Myciteauthor{johndrow2017bayes} derived a two-block Gibbs sampler for the Horseshoe posterior (the `exact algorithm' in \cite[Section 2.1]{johndrow2017bayes}, henceforth referred to as the JOB Gibbs sampler), and established geometric ergodicity 
(\cite[Theorem 14]{johndrow2017bayes}). However, the truncation assumptions needed for this result are rather restrictive, requiring all the local shrinkage parameters $\lambda_i^2$ to be bounded above by a finite constant, and also requiring the global shrinkage parameter $\tau^2$ to be bounded above and below by finite positive constants. In parallel work (\Myciteauthor{BBJJ:2020}, uploaded on arxiv a few days prior to our submission) geometric ergodicity for the JOB Gibbs sampler has now been established without requiring truncation of the local shrinkage parameters. However, the requirement of the global shrinkage parameter $\tau^2$ to be bounded above {\it and} below remains. 

{\bf Contribution \#1}: The first contribution of this paper is the proof of geometric ergodicity for a Horseshoe Gibbs sampler (see Theorem \ref{Theorem1}) with no truncation assumptions on the local shrinkage parameters, and {\it with the global shrinkage parameter only required to be truncated below by a finite positive constant and to have a finite $\delta^{th}$ prior moment for some $\delta > 0.00081$}. Hence, the 
conditions required for our geometric ergodicity result are strictly weaker than those in \cite{johndrow2017bayes} and \cite{BBJJ:2020}. Infact, as discussed in Remark \ref{Horseshoe:negative:moment}, the assumption of truncation below by a positive constant can be further relaxed to existence of the negative $(p+\delta)/2^{th}$ prior moment for some 
$\delta > 0.00162$. 

The Gibbs sampler analyzed in Theorem \ref{Theorem1} is a slight modification of the JOB Gibbs sampler with latent variables introduced to simplify conditional sampling of the local shrinkage parameters in the Markov chain (see Section \ref{original:horseshoe} for more details). There are also important differences in the technical arguments compared to \cite{johndrow2017bayes, BBJJ:2020}. We focus on the $\bl$-block of the Gibbs sampler and establish a drift condition (Lemma \ref{Theorem1}) using a drift function which is `unbounded off compact sets', and that directly leads to geometric ergodicity. On the other hand, the approaches in \cite{johndrow2017bayes, BBJJ:2020} use other drift functions (using all the parameters or a different parameter block than $\bl$) which are not unbounded off compact sets, and hence need an additional minorization argument. 

Next we move to the regularized Horseshoe setting of \Myciteauthor{piironen2017}. As mentioned previously, \cite{piironen2017} use a Hamiltonian Monte Carlo (HMC) based approach to generate approximate samples from the corresponding regularized Horseshoe posterior distribution, but do not investigate geometric ergodicity of the proposed Markov chain. It is not clear whether the intricate sufficient conditions needed for geometric ergodicity of HMC chains in \Myciteauthor{livingstone2019}  apply to the HMC chain in \cite{piironen2017}. Given the variety of efficient Gibbs samplers available for the original Horseshoe posterior, it is natural to consider an appropriately adapted version of any of these samplers for the regularized Horseshoe posterior. 

{\bf Contribution \#2}: As the second main contribution of this paper, we establish geometric ergodicity for one such Gibbs sampler for the regularized Horseshoe posterior (see Theorem \ref{thm:regularized:Horseshoe:Piironen:Vehtari}) {\it with no truncation assumptions on the global and local shrinkage parameters at all}. The seemingly minor change in the prior structure (compared to the original Horseshoe), leads to crucial changes in our convergence analysis. For example, we need a different drift function for this analysis (Lemma \ref{Theorem1reghorseshoe}) compared to the the original Horseshoe analysis. This drift function is not `unbounded off compact sets', and hence we need 
an additional minorization condition (Lemma \ref{Theorem2reghorseshoe}) to establish 
geometric ergodicity. 

Recently, \Myciteauthor{nishimura2019shrinkage} construct a further variant of the regularized Horseshoe prior of \cite{piironen2017} by changing the algebraic form of the conditional prior density of $\bb$ for computational simplicity. Their prior specification is as follows. 
\begin{eqnarray*}
\pi\left(\beta_j,\lambda_j^2\mid \tau^2,\sigma^2\right) &\propto&  \frac{1}{\sqrt{\tau^2\lambda_j^2}}\exp{\left[-\frac{\beta_j^2}{2\sigma^2}\left(\frac{1}{c^2}+\frac{1}{\tau^2\lambda_j^2}\right)\right]}\pi_\ell\left(\lambda_j\right)\\
& & \text{independently for}\ j=1,2,\cdots,p\\
\tau^2\sim \pi_\tau (\cdot) & & \hspace{0.2in} \sigma^2\sim \text{Inverse-Gamma}(a,b)
\end{eqnarray*}

\noindent
The algebraic modification, in particular removal of the $(c^{-2} + (\lambda_i \tau)^{-2})^{1/2}$ in the conditional prior for $\beta_i$ simplifies posterior 
computation (see Section \ref{regularized:Horseshoe:Nishimura:Suchard} for more details). \Myciteauthor{nishimura2019shrinkage} prove geometric ergodicity for the related  but structurally different setting of Polya-Gamma logistic regression assuming that the global shrinkage parameter $\tau^2$ is bounded above and below by finite positive constants. However, as discussed in Remark \ref{remregular}, several details of this analysis break down in the linear regression setting. 

{\bf Contribution \#3}: We focus on the linear regression setting, and leverage our analysis in the original Horseshoe setting to prove geometric ergodicity of a Gibbs sampler corresponding to \cite{nishimura2019shrinkage}'s regularized variant with the global shrinkage parameter only required to be bounded below by a finite positive constant and to have a finite $(p+\delta)/2^{th}$ moment for some $\delta > 0.00162$. 

The rest of the paper is structured as follows. We introduce the modified version of the JOB Gibbs sampler in Section \ref{Horseshoe:Gibbs:sampler}. Geometric ergodicity of this Gibbs sampler is established in Section \ref{Horseshoe:drift:condition}. The simulation study in Section \ref{simulation:original:Horseshoe} compares the computational time and other metrics for the JOB Gibbs sampler and the proposed modification in a variety of settings. An adaptation of the Horseshoe Gibbs sampler for the regularized Horseshoe posterior is developed in Section \ref{Regularized:Gibbs:sampler}. The geometric 
ergodicity of this regularized Horseshoe Gibbs sampler is established in 
Section \ref{Regularized:drift:minorization}. A related Gibbs sampler for the regularized Horseshoe variant of \cite{nishimura2019shrinkage} is discussed and analyzed in Section \ref{regularized:Horseshoe:Nishimura:Suchard}. Another simulation study in Section \ref{simulation:regularized:Horseshoe} examinies the computational feasibility/scalability of the Gibbs samplers analyzed in 
Sections \ref{Regularized:Gibbs:sampler} and 
\ref{regularized:Horseshoe:Nishimura:Suchard}. The proofs of several 
technical results used in the analysis are contained in an Appendix.

\section{Geometric ergodicity of a Horseshoe Gibbs sampler}\label{original:horseshoe}

\subsection{A modified version of the JOB Gibbs sampler} 
\label{Horseshoe:Gibbs:sampler}

\noindent
In this section, we describe in detail the Horseshoe Gibbs sampler that will be analyzed in subsequent sections. As pointed out in \Myciteauthor{Makalic_2016}, if $\lambda_j^2 \mid \nu_j \sim \mbox{Inverse-Gamma}(1/2, 1/\nu_j)$ and 
$\nu_j \sim \mbox{Inverse-Gamma}(1/2, 1)$, then $\lambda_j \sim C^+ (0,1)$. Using this 
fact, with ${\boldsymbol \nu} = (\nu_1, \nu_2, \cdots, \nu_p)$, the Horseshoe 
prior in 
(\ref{horseshoe_prior}) can be alternatively written as 
\begin{eqnarray}
\left. \bb \right. \mid \bl, \sigma^2,\tau^2\sim \mathcal{N}_p(0,\sigma^2\tau^2 \Bl) \nonumber\\
\left. \lambda_i^2 \right. \mid {\boldsymbol \nu} \sim \mbox{Inverse-Gamma} (1/2,1/\nu_i) \; \text{independently for}\ i=1,2,\cdots,p \nonumber\\
\nu_i \sim \mbox{Inverse-Gamma}(1/2, 1) \; \text{independently for}\ i = 1,2,\cdots, p \nonumber\\
\tau^2\sim \pi_\tau (\cdot), \hspace{0.2in} \sigma^2\sim 
\text{Inverse-Gamma}(a,b) \label{horseshoe_prior_augmented}
\end{eqnarray}

\noindent
Using the prior above and after straightforward calculations, various conditional 
posterior distributions can be derived as follows. 
\begin{eqnarray}\label{posteriors}
& \left.\boldsymbol{\beta}\right|\sigma^2,\tau^2,\boldsymbol{\lambda,\nu,y}\sim\mathcal{N}_p(A^{-1}\X^T\y,\sigma^2A^{-1})\nonumber \\
& \left.\sigma^2\right|\tau^2,\boldsymbol{\lambda,\nu,y}\sim\text{Inverse-Gamma}\left(a+\frac{n}{2},\frac{\y^T\left(I_n-\Tilde{P}_\X\right)\y}{2}+b\right)\nonumber \\
& \left.\lambda_j^2\right|\nu_j,\sigma^2,\tau^2,\beta_j, {\bf y} \sim\text{Inverse-Gamma}\left(1,\frac{1}{\nu_j}+\frac{\beta_j^2}{2\sigma^2\tau^2}\right)\ \text{independently for}\ i=1,2,\cdots,p\nonumber \\
& \left.\nu_j\right|\lambda_j^2, \tau^2, {\bf y} \sim\text{Inverse-Gamma}\left(1,1+\frac{1}{\lambda_j^2}\right)\ \text{independently for}\ i=1,2,\cdots,p\nonumber \\
& \left.\tau^2\right|\boldsymbol{\lambda,y}\sim\pi\left(\left.\tau^2\right|\boldsymbol{\lambda,y}\right)\propto\frac{\left(\frac{\y^T\left(I_n-\Tilde{P}_\X\right)\y}{2}+b\right)^{-\left(a+\frac{n}{2}\right)}}{\sqrt{\left|I_p+\X^T\X.\Bl_*\right|}}\cdot\pi_\tau(\tau^2)
\end{eqnarray}

\noindent
where $\Bl_*=\tau^2\Bl;\ A=\X^T\X+\Bl_*^{-1}$ and $\Tilde{P}_\X=\X 
A^{-1}\X^T$. 

Consider a two-block Gibbs sampling Markov chain with transition kernel 
$K_{aug}$ (with blocks $({\boldsymbol \beta}, \sigma^2, {\boldsymbol \nu}, 
\tau^2)$ and ${\boldsymbol \lambda}$) 
whose one-step transition from $({\boldsymbol \beta}_0, \sigma^2_0, {\boldsymbol \nu}_0, \tau^2_0, {\boldsymbol \lambda}_0)$ to 
$({\boldsymbol \beta}, \sigma^2, {\boldsymbol \nu}, \tau^2, {\boldsymbol \lambda})$ is given as follows. 
\begin{enumerate}
    \item Draw $({\boldsymbol \beta}, \sigma^2, {\boldsymbol \nu}, \tau^2)$ from $\pi({\boldsymbol \beta}, \sigma^2, {\boldsymbol \nu}, \tau^2 \mid {\boldsymbol \lambda}_0, {\bf y})$. This can be 
    done by sequentially drawing ${\boldsymbol \beta}$, then $\sigma^2$, 
    then ${\boldsymbol \nu}$, and then $\tau^2$ from appropriate conditional posterior densities in (\ref{posteriors}). 
    \item Draw ${\boldsymbol \lambda}$ from $\pi({\boldsymbol \lambda} \mid \boldsymbol{\beta,\nu},\sigma^2,\tau^2,\mathbf{\y})$. This can be done by 
    independently drawing the components of ${\boldsymbol \lambda}$ from the 
    appropriate full conditional posterior density in (\ref{posteriors}). 
\end{enumerate}

\noindent
The JOB Gibbs sampler from \Myciteauthor{johndrow2017bayes} is very similar to the above two-block Gibbs sampler $K_{aug}$. The difference is that the latent variables ${\boldsymbol \nu}$ are not used, and the two blocks used in the JOB Gibbs sampler are $({\boldsymbol \beta}, \sigma^2, \tau^2)$ and ${\boldsymbol \lambda}$. While the sampling steps for ${\boldsymbol \beta}, \sigma^2, \tau^2$
are exactly the same as above, the components of ${\boldsymbol \lambda}$ are sampled differently. In particular, each $\lambda_j$ is sampled from the conditional density given $\beta_j, \sigma^2, \tau^2, {\bf y}$ (no conditioning on $\nu_j$). This conditional density is not a standard density, and draws are made using a rejection sampler. To summarize, by considering the latent variables ${\boldsymbol \nu}$, we replace the $p$ rejection sampler based 
draws from a non-standard density in the JOB Gibbs sampler (for components of ${\boldsymbol \lambda}$) with $2p$ draws from standard Inverse-Gamma densities (for components of 
${\boldsymbol \lambda}$ and ${\boldsymbol \nu}$). 

The Gibbs sampler $K_{aug}$ can essentially be considered a hybrid of the JOB Gibbs sampler and the Gibbs sampler in \Myciteauthor{Makalic_2016}, which uses a latent variable $\xi$ (in addition to ${\boldsymbol \nu}$) to replace the draws from the non-standard $\pi(\tau^2 \mid {\boldsymbol \lambda}, {\bf y})$ density with two draws from standard Inverse-Gamma densities. As mentioned in the introduction, the geometric ergodicity result for the JOB Gibbs sampler in 
\cite[Theorem 14]{johndrow2017bayes} has been established by assuming that the local shrinkage parameters in $\bl$ are all bounded above, and the global 
shrinakge parameter $\tau^2$ is bounded above and below. In very recent follow-up work \cite{BBJJ:2020}, the authors establish geoemtric ergodcity for a class of Half-$t$ Gibbs samplers of which the JOB Gibbs sampler is a member. In this work, the truncation assumption on the local shrinkage parameters has been removed, but the global shrinkage parameter is still 
assumed to be truncated above {\it and} below. However, we show below that geometric ergodicity for the hybrid Gibbs sampler $K_{aug}$ can be established with no truncation at all on the local shrinkage parameters in ${\boldsymbol \lambda}$, and only assuming that the global shrinkage parameter $\tau^2$ is truncated below. 

The reasons for this improved analysis of the hybrid chain $K_{aug}$ lie in the intricacies of drift and minorization approach (\cite{10.2307/2291067}), which is the state of the art technique for proving geometric ergodicity for general state space Markov chains. The introduction of the latent variables ${\boldsymbol \nu}$, the resulting Inverse-Gamma posterior conditionals for entries of ${\boldsymbol \lambda}$ and ${\boldsymbol \nu}$, and avoiding the latent variable $\xi$ for the global shrinkage parameter $\tau^2$ provide just the right ingredients for establishing a geometric drift condition in Section \ref{Horseshoe:drift:condition} which is then leveraged to establish geometric ergodicity. Even a minor deviation in the structure of the Markov chain (such as in the JOB Gibbs sampler or the Gibbs sampler of \cite{Makalic_2016}) leads to a breakdown of the intricate argument. 

Before proceeding further, we note that geometric ergodicity of a two-block Gibbs sampler can be established by showing that any of its two marginal chains is geometrically ergodic (see for example \cite{Roberts:Rosenthal:2001}). Hence, we focus on the marginal ${\boldsymbol \lambda}$-chain corresponding to $K_{aug}$. The one-step transition dynamics of this Markov chain from ${\boldsymbol \lambda}_m$ to ${\boldsymbol \lambda}_{m+1}$ is given as follows: 
\begin{enumerate}
    \item Draw $\tau^2$ from $\pi\left(\left.\tau^2\right|\mathbf{\bl}_m,\mathbf{\y}\right)$
    \item Draw $\boldsymbol{\nu}$ from $\pi\left(\left.\boldsymbol{\nu}\right|\mathbf{\bl}_m,\tau^2,\mathbf{\y}\right)=\prod\limits_{j=1}^p\text{Inverse-Gamma}\left(1,1+\frac{1}{\lambda_{j;m}^2}\right)$
    \item Draw $\sigma^2$ from $\pi\left(\left.\sigma^2\right|\tau^2,\mathbf{\bl}_m,{\boldsymbol \nu},\mathbf{\y}\right)=\text{Inverse-Gamma}\left(a+\frac{n}{2},\frac{\mathbf{\y}^T\left(I_n-\Tilde{P}_\X\right)\mathbf{\y}}{2}+b\right)$
    \item Draw $\boldsymbol{\beta}$ from $\pi\left(\left.\boldsymbol{\beta}\right|\sigma^2,\tau^2,\mathbf{\bl}_m,{\boldsymbol \nu},\mathbf{\y}\right)=\mathcal{N}_p(A^{-1}\X^T\y,\sigma^2A^{-1})$
    \item Finally draw $\mathbf{\bl}_{m+1}$ from $\pi\left(\left.\mathbf{\bl}\right|\boldsymbol{\beta,\nu},\sigma^2,\tau^2,\mathbf{\y}\right)=\prod\limits_{j=1}^p\text{Inverse-Gamma}\left(1,\frac{1}{\nu_j}+\frac{\beta_j^2}{2\sigma^2\tau^2}\right)$
\end{enumerate}
The Markov transition density (MTD) corresponding to the marginal ${\boldsymbol \lambda}$-chain is given by 
\begin{eqnarray}\label{MTD}
    k\left(\mathbf{\bl}_0,\mathbf{\bl}\right)= \int_{\mathbb{R}_+}\int_{\mathbb{R}_{+}}\int_{\mathbb{R}^p}\int_{\mathbb{R}_+^p}\pi\left(\left.\mathbf{\bl}\right|\boldsymbol{\beta,\nu},\sigma^2,\tau^2,\mathbf{\y}\right)\pi\left(\left.\boldsymbol{\beta,\nu},\sigma^2,\tau^2\right|\mathbf{\bl}_0,\mathbf{\y}\right)d\boldsymbol{\nu}d\boldsymbol{\beta}d\sigma^2d\tau^2 \nonumber\\
   = \int_{\mathbb{R}_+}\int_{\mathbb{R}_{+}}\int_{\mathbb{R}^p}\int_{\mathbb{R}^{p}_+}\pi\left(\left.\mathbf{\bl}\right|\boldsymbol{\beta,\nu},\sigma^2,\tau^2,\mathbf{\y}\right)\pi\left(\left.\boldsymbol{\beta}\right|\sigma^2,\tau^2,\mathbf{\bl}_0,{\boldsymbol \nu},\mathbf{\y}\right)\nonumber \\
    \times\pi\left(\left.\sigma^2\right|\tau^2,\mathbf{\bl}_0,{\boldsymbol \nu},\mathbf{\y}\right)\pi\left(\left.\boldsymbol{\nu}\right|\mathbf{\bl}_0,\tau^2,\mathbf{\y}\right)\pi\left(\left.\tau^2\right|\mathbf{\bl}_0,\mathbf{\y}\right)d\boldsymbol{\nu}d\boldsymbol{\beta}d\sigma^2d\tau^2 
\end{eqnarray}

\noindent
We now establish a drift condition for the marginal ${\boldsymbol \lambda}$-chain, which will then be used to establish geometric ergodicity for 
the two-block Horseshoe Gibbs sampler $K_{aug}$. 

\subsection{A drift condition for the ${\bl}$-chain}\label{Horseshoe:drift:condition}

\noindent
Consider the function $V:\mathbb{R}_+^p\mapsto[0,\infty)$ given by 
\begin{equation}\label{driftfunction}
    V\left(\mathbf{\bl}\right)=\sum_{j=1}^p\left(\lambda_j^2\right)^{\frac{\delta_0}{2}}+\sum_{j=1}^p\left(\lambda_j^2\right)^{-\frac{\delta_1}{2}}, 
\end{equation}

\noindent
where $\delta_0,\delta_1\in(0,1)$ are some constants. The next result establishes a geometric drift condition for the marginal ${\boldsymbol \lambda}$-chain using the function $V$ with appropriately small values of $\delta_0$ and $\delta_1$. 
\begin{lem}\label{Theorem1}
Suppose the prior density $\pi_\tau$ for the global shrinkage parameter is 
truncated below i.e., $\pi_\tau (u) = 0$ for $ u < T$ for some $T > 0$ and satisfies $$\int_T^{\infty} u^{\delta/2}\pi_\tau (u)du<\infty$$ for some $\delta\in (0.00162,0.22176)$.  
Then, there exist $\delta_0, \delta_1 \in (0,1)$ such that for every 
$\bl_0\in \mathbb{R}_+^p$ we have 
\begin{equation}\label{driftcondition}
    \E\left[\left.V\left(\mathbf{\bl}\right)\right|\mathbf{\bl}_0\right] = 
    \int_{\mathbb{R}_+^p} k\left( \mathbf{\bl}_0,\mathbf{\bl} \right) 
    V(\mathbf{\bl}) d {\mathbf{\bl}} \leq \gamma^* V\left(\mathbf{\bl}_0\right)+b^*
\end{equation}
 
\noindent
with $0<\gamma^*=\gamma^*\left(\delta_0,\delta_1\right)<1$ and $b^* = b^* 
\left( \delta_0, \delta_1 \right) <\infty$. 
\end{lem}

\proof Note that by linearity 
\begin{equation} \label{expectations:linearity}
\E\left[\left.V\left(\mathbf{\bl}\right)\right|\mathbf{\bl}_0\right] = 
\sum_{j=1}^p \E\left[\left. \left(\lambda_j^2\right)^{\frac{\delta_0}{2}} \right|\mathbf{\bl}_0\right] + \sum_{j=1}^p \E\left[\left. \left(\lambda_j^2\right)^{-\frac{\delta_1}{2}} \right|\mathbf{\bl}_0\right]
\end{equation}

\noindent
We first consider terms in the second sum in (\ref{expectations:linearity}). Fix $j\in\left\{1,2,\cdots,p\right\}$ arbitrarily. It follows from the definition of the MTD \eqref{MTD} that 
\begin{eqnarray}
& & \E\left[\left.\left(\lambda_j^2\right)^{-\frac{\delta_1}{2}}\right|\mathbf{\bl}_0\right] \nonumber\\
&=& \E\left[\left.\E\left[\left.\E\left[\left.\E\left[\left.\E\left[\left.
\left(\lambda_j^2\right)^{-\frac{\delta_1}{2}}\right.\mid\boldsymbol{\beta,\nu},\sigma^2,\tau^2,\mathbf{\y}\right]\right.\mid\sigma^2,\tau^2,\mathbf{\bl}_0,{\boldsymbol \nu},\mathbf{\y}\right]\right.\mid\tau^2,\mathbf{\bl}_0,{\boldsymbol \nu},\mathbf{\y}\right]\right.\mid\bl_0,\tau^2,\y\right]\right.\mid\bl_0,\y\right]. \nonumber\\
& & \label{expectations:iterated}
\end{eqnarray}

\noindent
The five iterated expectations correspond to the five conditional densities in (\ref{MTD}). Starting with the innermost expectation, and using the fact that $1/\lambda_j^2$ (conditioned on $\boldsymbol{\beta,\nu},\sigma^2,\tau^2,\mathbf{\y}$) follows a Gamma distribution with shape parameter $1$ and rate parameter $1/\nu_j + \beta_j^2/(2 \sigma^2 \tau^2)$, we obtain 
\begin{eqnarray*}
    \E\left[\left.\left(\lambda_j^2\right)^{-\frac{\delta_1}{2}}\right|\boldsymbol{\beta,\nu},\sigma^2,\tau^2,\mathbf{\y}\right]&=& \Gamma\left(1+\frac{\delta_1}{2}\right)\left(\frac{1}{\nu_j}+\frac{\beta_j^2}{2\sigma^2\tau^2}\right)^{-\frac{\delta_1}{2}} \\
    &=& \Gamma\left(1+\frac{\delta_1}{2}\right)\left(\frac{1}{\nu_j}+\frac{1}{\left\{\frac{\left(2\sigma^2\tau^2\right)^{\frac{\delta_1}{2}}}{\left|\beta_j\right|^{\delta_1}}\right\}^{\frac{2}{\delta_1}}}\right)^{-\frac{\delta_1}{2}} \\
\end{eqnarray*}

\noindent
Note that the function $y\mapsto\left(c+y^{-\frac{2}{\delta_1}}\right)^{-\delta_1/2}$ on $(0,\infty)$ is concave for $c > 0, \delta_1 \in (0,1)$. Applying the second iterated expectation, and using Jensen's inequality, it follows that 
\begin{eqnarray}
&&    \E\left[\left.\E\left[\left.\left(\lambda_j^2\right)^{-\frac{\delta_1}{2}}\right|\boldsymbol{\beta,\nu},\sigma^2,\tau^2,\mathbf{\y}\right]\right|\sigma^2,\tau^2,\mathbf{\bl}_0,\mathbf{\y}\right]\nonumber \\
     &\leq& \Gamma\left(1+\frac{\delta_1}{2}\right)\left(\frac{1}{\nu_j}+\frac{1}{\left\{E\left[\left.\frac{\left(2\sigma^2\tau^2\right)^{\frac{\delta_1}{2}}}{\left|\beta_j\right|^{\delta_1}}\right|\sigma^2,\tau^2,\mathbf{\lambda}_0,{\boldsymbol \nu},\mathbf{y}\right]\right\}^{\frac{2}{\delta_1}}}\right)^{-\frac{\delta_1}{2}} \label{jensen}
     \end{eqnarray}
     
     \noindent
     Note that the conditional distribution of $\beta_j$ given $\sigma^2,\tau^2,\mathbf{\bl}_0,{\boldsymbol \nu},\mathbf{\y}$ is a Gaussian distribution with variance $\sigma_j^2\overset{\text{def}}{=}\sigma^2\boldsymbol{e}_J^TA^{-1}\boldsymbol{e}_j\geq \sigma^2\left(\bar{\omega}+\frac{1}{\tau^2\lambda_{j;0}^2}\right)^{-1}$. 
     Here $\bar{\omega}$ is the maximum eigenvalue of $\X^T\X$ and $\boldsymbol{e}_j$ is the $p \times 1$ vector with $j^{\text{th}}$ entry $1$ and other entries equal to $0$. Using Proposition A1 from \Myciteauthor{pal2014} regarding the negative moments of a Gaussian random variable and choosing 
     $\delta_1 \in (0,1)$, it follows that 
\begin{equation}
E\left[\left.\frac{\left(2\sigma^2\tau^2\right)^{\frac{\delta_1}{2}}}{\left|\beta_j\right|^{\delta_1}} \right| \sigma^2,\tau^2,\mathbf{\bl}_0,{\boldsymbol \nu},\mathbf{\y} \right]\leq \left(2\sigma^2\tau^2\right)^{\frac{\delta_1}{2}}\frac{\Gamma\left(\frac{1-\delta_1}{2}\right)2^{\frac{1-\delta_1}{2}}}{\sqrt{2\pi} \sigma_j^{\delta_1}} \leq \frac{\Gamma\left(\frac{1-\delta_1}{2}\right)}{\sqrt{\pi}}\left(\bar{\omega}\tau^2+\frac{1}{\lambda_{j;0}^2}\right)^{\frac{\delta_1}{2}} \label{Gaussian:moment}\\
\end{equation}

\noindent
Combining (\ref{jensen}) and (\ref{Gaussian:moment}), we obtain 
\begin{eqnarray*}
&&    \E\left[\left.\E\left[\left.\left(\lambda_j^2\right)^{-\frac{\delta_1}{2}}\right|\boldsymbol{\beta,\nu},\sigma^2,\tau^2,\mathbf{\y}\right]\right|\sigma^2,\tau^2,\mathbf{\bl}_0,\mathbf{\y}\right]\\
     &\leq& \Gamma\left(1+\frac{\delta_1}{2}\right)\left(\frac{1}{\nu_j}+\frac{1}{\left\{\frac{\Gamma\left(\frac{1-\delta_1}{2}\right)}{\sqrt{\pi}}\left(\bar{\omega}\tau^2+\frac{1}{\lambda_{j;0}^2}\right)^{\frac{\delta_1}{2}}\right\}^{\frac{2}{\delta_1}}}\right)^{-\frac{\delta_1}{2}}. 
\end{eqnarray*}
     
Using the fact $\left(u+v\right)^{\delta}\leq u^{\delta}+v^{\delta}$ for $\delta\in(0,1)$ and $u,v\geq 0$, it follows that 
\begin{eqnarray}
&&    \E\left[\left.\E\left[\left.\left(\lambda_j^2\right)^{-\frac{\delta_1}{2}}\right|\boldsymbol{\beta,\nu},\sigma^2,\tau^2,\mathbf{\y}\right]\right|\sigma^2,\tau^2,\mathbf{\bl}_0,\mathbf{\y}\right] \nonumber\\
     &\leq& \Gamma\left(1+\frac{\delta_1}{2}\right)\left(\frac{1}{\nu_j}+\frac{1}{\left\{\frac{\Gamma\left(\frac{1-\delta_1}{2}\right)}{\sqrt{\pi}}\left(\bar{\omega}^{\frac{\delta_1}{2}}\left(\tau^2\right)^{\frac{\delta_1}{2}}+\left(\lambda_{j;0}^2\right)^{-\frac{\delta_1}{2}}\right)\right\}^{\frac{2}{\delta_1}}}\right)^{-\frac{\delta_1}{2}}. \label{power:inequality}
\end{eqnarray}

\noindent
Note that the bound in (\ref{power:inequality}) does not depend on $\sigma^2$. Again, using the fact that  $y\mapsto\left(c+y^{-\frac{2}{\delta_1}}\right)^{-\delta_1/2}$ on $(0,\infty)$ is concave for $c > 0, \delta_1 \in (0,1)$, along with  Jensen's inequality, we get 
\begin{eqnarray*}
& & \E\left[\left.\left(\frac{1}{\nu_j}+\frac{1}{\left\{\frac{\Gamma\left(\frac{1-\delta_1}{2}\right)}{\sqrt{\pi}}\left(\bar{\omega}^{\frac{\delta_1}{2}}\left(\tau^2\right)^{\frac{\delta_1}{2}}+\left(\lambda_{j;0}^2\right)^{-\frac{\delta_1}{2}}\right)\right\}^{\frac{2}{\delta_1}}}\right)^{-\frac{\delta_1}{2}}\right|\tau^2,\mathbf{\bl}_0,\mathbf{y}\right]\\
   &\leq& \left(\frac{1}{\left\{E\left[\left.\nu_j^{\frac{\delta_1}{2}}\right|\tau^2,\mathbf{\lambda}_0,\mathbf{y}\right]\right\}^{\frac{2}{\delta_1}}}+\frac{1}{\left\{\frac{\Gamma\left(\frac{1-\delta_1}{2}\right)}{\sqrt{\pi}}\left(\bar{\omega}^{\frac{\delta_1}{2}}\left(\tau^2\right)^{\frac{\delta_1}{2}}+\left(\lambda_{j;0}^2\right)^{-\frac{\delta_1}{2}}\right)\right\}^{\frac{2}{\delta_1}}}\right)^{-\frac{\delta_1}{2}} 
   \end{eqnarray*}

\noindent
Since $\nu_j$ (given $\tau^2,\mathbf{\lambda}_0,\mathbf{y}$) has an 
Inverse-Gamma distribution with shape parameter $1$ and rate parameter $1+1/\lambda_{j;0}^2$, it follows that 

\begin{eqnarray}
   & & \E\left[\left.\left(\frac{1}{\nu_j}+\frac{1}{\left\{\frac{\Gamma\left(\frac{1-\delta_1}{2}\right)}{\sqrt{\pi}}\left(\bar{\omega}^{\frac{\delta_1}{2}}\left(\tau^2\right)^{\frac{\delta_1}{2}}+\left(\lambda_{j;0}^2\right)^{-\frac{\delta_1}{2}}\right)\right\}^{\frac{2}{\delta_1}}}\right)^{-\frac{\delta_1}{2}}\right|\tau^2,\mathbf{\bl}_0,\mathbf{y}\right] 
   \nonumber\\
   &=& \left(\frac{1}{\left\{\Gamma\left(1-\frac{\delta_1}{2}\right)\left(1+\frac{1}{\lambda_{j;0}^2}\right)^{\frac{\delta_1}{2}}\right\}^{\frac{2}{\delta_1}}}+\frac{1}{\left\{\frac{\Gamma\left(\frac{1-\delta_1}{2}\right)}{\sqrt{\pi}}\left(\bar{\omega}^{\frac{\delta_1}{2}}\left(\tau^2\right)^{\frac{\delta_1}{2}}+\left(\lambda_{j;0}^2\right)^{-\frac{\delta_1}{2}}\right)\right\}^{\frac{2}{\delta_1}}}\right)^{-\frac{\delta_1}{2}}. 
   \label{moment}
   \end{eqnarray}

  \noindent
  Let us now take the expectation of the expression in (\ref{moment}) with respect to the conditional distribution of $\tau^2$ given $\mathbf{\bl}_0, \mathbf{y}$. Using for a third time the fact that $y\mapsto\left(c+y^{-\frac{2}{\delta_1}}\right)^{-\delta_1/2}$ on $(0,\infty)$ is concave for $c > 0, \delta_1 \in (0,1)$, along with  Jensen's inequality, we get
  \begin{eqnarray}
   & & \E\left[\left.\left(\frac{1}{\left\{\Gamma\left(1-\frac{\delta_1}{2}\right)\left(1+\frac{1}{\lambda_{j;0}^{\delta_1}}\right)\right\}^{\frac{2}{\delta_1}}}+\frac{1}{\left\{\frac{\Gamma\left(\frac{1-\delta_1}{2}\right)}{\sqrt{\pi}}\left(\bar{\omega}^{\frac{\delta_1}{2}}\left(\tau^2\right)^{\frac{\delta_1}{2}}+\left(\lambda_{j;0}^2\right)^{-\frac{\delta_1}{2}}\right)\right\}^{\frac{2}{\delta_1}}}\right)^{-\frac{\delta_1}{2}}\right|\mathbf{\bl}_0,\mathbf{y}\right] \nonumber \\
   &\leq& \left(\frac{1}{\left\{\Gamma\left(1-\frac{\delta_1}{2}\right)\left(1+\frac{1}{\lambda_{j;0}^{\delta_1}}\right)\right\}^{\frac{2}{\delta_1}}}+\frac{1}{\left\{\frac{\Gamma\left(\frac{1-\delta_1}{2}\right)}{\sqrt{\pi}}\left(\bar{\omega}^{\frac{\delta_1}{2}}\E\left[\left.\left(\tau^2\right)^{\frac{\delta_1}{2}}\right|\mathbf{\bl}_0,\mathbf{y}\right]+\left(\lambda_{j;0}^2\right)^{-\frac{\delta_1}{2}}\right)\right\}^{\frac{2}{\delta_1}}}\right)^{-\frac{\delta_1}{2}} \nonumber \allowdisplaybreaks\\
   &\overset{(\bigstar)}{\leq}&  \left(\frac{1}{\left\{\Gamma\left(1-\frac{\delta_1}{2}\right)\left(1+\frac{1}{\lambda_{j;0}^{\delta_1}}\right)\right\}^{\frac{2}{\delta_1}}}+\frac{1}{\left\{\frac{\Gamma\left(\frac{1-\delta_1}{2}\right)}{\sqrt{\pi}}\left(\bar{\omega}^{\frac{\delta_1}{2}}C_1+\frac{1}{\lambda_{j;0}^{\delta_1}}\right)\right\}^{\frac{2}{\delta_1}}}\right)^{-\frac{\delta_1}{2}} \nonumber \\
  &\leq& \left(C_0+\frac{1}{\lambda_{j;0}^{\delta_1}}\right)\left(\frac{1}{\left\{\Gamma\left(1-\frac{\delta_1}{2}\right)\right\}^{\frac{2}{\delta_1}}}+\frac{\sqrt{\pi}^{\frac{2}{\delta_1}}}{\left\{\Gamma\left(\frac{1-\delta_1}{2}\right)\right\}^{\frac{2}{\delta_1}}}\right)^{-\frac{\delta_1}{2}};\ \ C_0=\text{max}\left\{1,\bar{\omega}^{\frac{\delta_1}{2}}C_1\right\} \label{lower_truncation}
\end{eqnarray}

\noindent
where $(\bigstar)$ follows from Proposition \ref{upperboundontau^2}. Combining (\ref{expectations:iterated}), \eqref{jensen},  (\ref{power:inequality}), 
(\ref{moment}) and (\ref{lower_truncation}), we get 
\begin{eqnarray}\label{drift:1}
   \E\left[\left.\sum_{j=1}^p\left(\lambda_j^2\right)^{-\frac{\delta_1}{2}}\right|\mathbf{\bl}_0\right]\leq \gamma\left(\delta_1\right)\sum_{j=1}^p\left(\lambda_{j;0}^2\right)^{-\frac{\delta_1}{2}}+b_1
\end{eqnarray}

\noindent
where 
$$
\gamma\left(\delta_1\right)=\Gamma\left(1+\frac{\delta_1}{2}\right)\left(\frac{1}{\left\{\Gamma\left(1-\frac{\delta_1}{2}\right)\right\}^{\frac{2}{\delta_1}}}+\frac{\sqrt{\pi}^{\frac{2}{\delta_1}}}{\left\{\Gamma\left(\frac{1-\delta_1}{2}\right)\right\}^{\frac{2}{\delta_1}}}\right)^{-\frac{\delta_1}{2}} 
$$ 

\noindent
and 
$$
b_1=p\cdot C_0\cdot\gamma\left(\delta_1\right). 
$$

\noindent
Next consider $\E\left[\left.\sum\limits_{j=1}^p\left(\lambda_j^2\right)^{\frac{\delta_0}{2}}\right|\mathbf{\bl}_0\right]$. 
Fix a $j\in\left\{1,2,\cdots,p\right\}$ arbitrarily. Since $\delta_0 \in (0,1)$, using the fact that $(u+v)^{\delta_0}\leq u^{\delta_0}+v^{\delta_0}$ for $u,v\geq 0$ we get \begin{align*}
    \E\left[\left.\left(\lambda_j^2\right)^{\frac{\delta_0}{2}}\right|\boldsymbol{\beta,\nu},\sigma^2,\tau^2,\mathbf{y}\right] = {} & \Gamma\left(1-\frac{\delta_0}{2}\right)\left(\frac{1}{\nu_j}+\frac{\beta_j^2}{2\sigma^2\tau^2}\right)^{\frac{\delta_0}{2}} \\
     \leq {} & \Gamma\left(1-\frac{\delta_0}{2}\right)\left(\frac{1}{\nu_j^{\frac{\delta_0}{2}}}+\frac{\left|\beta_j\right|^{\delta_0}}{\left(2\sigma^2\tau^2\right)^{\frac{\delta_0}{2}}}\right). 
\end{align*}

\noindent

For  $j=1,2,\cdots,p$, we denote 
\begin{equation} \label{conditional_mean_definition}
\mu_j = \boldsymbol{e}_j^TA_0^{-1}\X^T\y 
\end{equation}

\noindent

where $A_0 = \X^T\X + (\tau^2 \Bl_0)^{-1}$. 

It follows that 
\begin{eqnarray*}
& &    \E\left[\left.\E\left[\left.\left(\lambda_j^2\right)^{\frac{\delta_0}{2}}\right|\boldsymbol{\beta,\nu},\sigma^2,\tau^2,\mathbf{y}\right]\right|\sigma^2,\tau^2,\mathbf{\bl}_0,\mathbf{y}\right] \\
&\leq& \Gamma\left(1-\frac{\delta_0}{2}\right)\E\left[\left.\left(\frac{1}{\nu_j^{\frac{\delta_0}{2}}}+\frac{\left|\beta_j\right|^{\delta_0}}{\left(2\sigma^2\tau^2\right)^{\frac{\delta_0}{2}}}\right)\right|\sigma^2,\tau^2,\mathbf{\bl}_0,\mathbf{y}\right] \nonumber\\
     &=& \Gamma\left(1-\frac{\delta_0}{2}\right)\left(\E\left[\left.\frac{1}{\nu_j^{\frac{\delta_0}{2}}}\right|\sigma^2,\tau^2,\mathbf{\bl}_0,\mathbf{y}\right]+\E\left[\left.\frac{\left|\beta_j\right|^{\delta_0}}{\left(2\sigma^2\tau^2\right)^{\frac{\delta_0}{2}}}\right|\sigma^2,\tau^2,\mathbf{\bl}_0,\mathbf{y}\right]\right) \nonumber\\
     &\leq& \Gamma\left(1-\frac{\delta_0}{2}\right)\left(\Gamma\left(1+\frac{\delta_0}{2}\right)+\E\left[\left.\frac{\left|\beta_j-\mu_j\right|^{\delta_0}}{\left(2\sigma^2\tau^2\right)^{\frac{\delta_0}{2}}}\right|\sigma^2,\tau^2,\mathbf{\bl}_0,\mathbf{y}\right]+\frac{\left|\mu_j\right|^{\delta_0}}{\left(2\sigma^2\tau^2\right)^{\frac{\delta_0}{2}}}\right); \nonumber \\
     &\leq& \Gamma\left(1-\frac{\delta_0}{2}\right)\left(\Gamma\left(1+\frac{\delta_0}{2}\right)+\frac{\Gamma\left(\frac{1+\delta_0}{2}\right)}{\sqrt{\pi}}\lambda_{j;0}^{\delta_0}+\frac{\left|\mu_j\right|^{\delta_0}}{\left(2\sigma^2\tau^2\right)^{\frac{\delta_0}{2}}}\right) \nonumber \\
     &\overset{(\bigstar\bigstar)}{\leq}&\Gamma\left(1-\frac{\delta_0}{2}\right)\left(\Gamma\left(1+\frac{\delta_0}{2}\right)+\frac{\Gamma\left(\frac{1+\delta_0}{2}\right)}{\sqrt{\pi}}\lambda_{j;0}^{\delta_0}+\frac{T^*}{\left(2\sigma^2\right)^{\frac{\delta_0}{2}}}\right), 
\end{eqnarray*}

\noindent
for some $T^* > 0$. Here $(\bigstar\bigstar)$ follows from 
Proposition $\ref{result:lmUnifBound2}$ (see \hyperref[Appendix A]{Appendix \ref{Appendix A}}) and the fact that $\tau^2$ is supported on $\left[T,\infty\right)$. Hence,
\begin{eqnarray*}
& &   \E\left[\left.\left(\lambda_j^2\right)^{\frac{\delta_0}{2}}\right|\mathbf{\bl}_0\right]=\E\left[\left.\E\left[\left.\E\left[\left.\left(\lambda_j^2\right)^{\frac{\delta_0}{2}}\right|\boldsymbol{\beta,\nu},\sigma^2,\tau^2,\mathbf{\y}\right]\right|\sigma^2,\tau^2,\mathbf{\bl}_0,\mathbf{\y}\right]\right|\mathbf{\bl}_0,\mathbf{\y}\right] \\
   &\leq& \Gamma\left(1-\frac{\delta_0}{2}\right)\left(\Gamma\left(1+\frac{\delta_0}{2}\right)+\frac{\Gamma\left(\frac{1+\delta_0}{2}\right)}{\sqrt{\pi}}\lambda_{j;0}^{\delta_0}+\E\left[\left.\frac{T^*}{\left(2\sigma^2\right)^{\frac{\delta_0}{2}}}\right|\bl_0,\y\right]\right) \nonumber \\
   &=& \Gamma\left(1-\frac{\delta_0}{2}\right)\left(\Gamma\left(1+\frac{\delta_0}{2}\right)+\frac{\Gamma\left(\frac{1+\delta_0}{2}\right)}{\sqrt{\pi}}\lambda_{j;0}^{\delta_0}+\E\left[\left.\E\left[\left.\frac{T^*}{\left(2\sigma^2\right)^{\frac{\delta_0}{2}}}\right|\tau^2,\bl_0,\y\right]\right|\bl_0,\y\right]\right) \nonumber \\
   &\leq& \Gamma\left(1-\frac{\delta_0}{2}\right)\left(\Gamma\left(1+\frac{\delta_0}{2}\right)+\frac{\Gamma\left(\frac{1+\delta_0}{2}\right)}{\sqrt{\pi}}\lambda_{j;0}^{\delta_0}+\frac{T^*}{\left(2b\right)^{\frac{\delta_0}{2}}}\cdot\frac{\Gamma\left(a+\frac{n+\delta_0}{2}\right)}{\Gamma\left(a+\frac{n}{2}\right)}\right) 
\end{eqnarray*}

\noindent
It follows that 
\begin{eqnarray}\label{drift:2}
   \E\left[\left.\sum_{j=1}^p\left(\lambda_j^2\right)^{\frac{\delta_0}{2}}\right|\mathbf{\bl}_0\right]\leq \gamma\left(\delta_0\right)\sum_{j=1}^p\left(\lambda_{j;0}^2\right)^{\frac{\delta_0}{2}}+b_2
\end{eqnarray}

\noindent
where 
$$
\gamma\left(\delta_0\right)=\Gamma\left(1-\frac{\delta_0}{2}\right)\frac{\Gamma\left(\frac{1+\delta_0}{2}\right)}{\sqrt{\pi}}
$$

\noindent
and 
$$
b_2=p\cdot\Gamma\left(1-\frac{\delta_0}{2}\right)\Gamma\left(1+\frac{\delta_0}{2}\right)\frac{T^*}{\left(2b\right)^{\frac{\delta_0}{2}}}\cdot\frac{\Gamma\left(a+\frac{n+\delta_0}{2}\right)}{\Gamma\left(a+\frac{n}{2}\right)}. 
$$

\noindent
The result follows by combining $\eqref{drift:1}$ and $\eqref{drift:2}$ with $$\gamma^*=\text{max}\left\{\gamma\left(\delta_0\right),\gamma\left(\delta_1\right)\right\}$$ and $$b^*=b_1+b_2.$$ 
Note that $\gamma^*=\text{max}\left\{\gamma\left(\delta_0\right),\gamma\left(\delta_1\right)\right\} < 1$ for small enough choices of $\delta_0$ and 
$\delta_1$, for example $\delta_0,\delta_1\in (0.00162, 0.22176)$. \qed

\begin{remark} \label{Horseshoe:negative:moment}
Note that the only place in the proof of Lemma \ref{Theorem1} where we 
need $\tau^2$ to be truncated below is to show that ${\bf E} \left[ 
\left( \tau^2 \right)^{-\frac{\delta_0}{2}} \mid \bl_0, {\bf y} \right]$ 
is uniformly bounded in $\bl_0$. In 
Proposition \ref{upperboundontau^2:negative:moment}, we show this follows by assuming the weaker condition that the {\it prior} negative $(p+\delta_0)/2^{th}$ moment for $\tau^2$ is finite. 
\end{remark}

\noindent
We now explain why the geometric drift condition established in 
Theorem \ref{Theorem1} for the marginal $\bl$-chain implies geometric 
ergodicity of the two-block Horseshoe Gibbs sampler $K_{aug}$. Note that 
for every $d \in \mathbb{R}$, the set 
$$
B \left( V,d \right) = \left\{\bl \in \mathbb{R}_+^p: V(\bl)=\sum_{j=1}^p\left(\lambda_{j}^2\right)^{\frac{\delta_0}{2}}+\sum_{j=1}^p\left(\lambda_{j}^2\right)^{-\frac{\delta_1}{2}} \leq d\right\}
$$ \label{originalhorseshoeminorizationset}

\noindent
is a compact set. Since $k(\bl_0, \bl)$ is continuous in $\bl_0$, a 
standard argument using Fatou's lemma along with Theorem 6.0.1 of \Myciteauthor{meyn1993markov} can be used to establish that 
the marginal $\bl$-chain is {\it unbounded off petite sets}. 
 Lemma 15.2.8 of \Myciteauthor{meyn1993markov} then implies 
geometric ergodicity of the marginal $\lambda$-chain. Using   Lemma 2.4 in \Myciteauthor{diaconis2008} now 
gives the following result. 
\begin{thm} \label{Theorem 2.1}
Suppose the prior density $\pi_\tau$ for the global shrinkage parameter 
truncated below i.e., $\pi_\tau (u) = 0$ for $ u < T$ for some $T > 0$ and satisfies $$\int_0^{\infty}u^{\delta/2}\pi_\tau\left(u\right)du<\infty$$ for some $\delta\in (0.00162, 0.22176)$. 
Then the two-block Horseshoe Gibbs sampler with transition kernel $K_{aug}$
is geometrically ergodic. The assumption of truncation below (i.e., 
$T > 0$) can be replaced by the weaker assumption that $T = 0$ and that 
the {\it prior} negative $(p+\delta)/2^{th}$ moment for $\tau^2$ is 
finite for some $\delta > 0.00162$. 
\end{thm}

\noindent
Note that the above result establishes geometric ergodicity, which as 
described earlier, helps rigorously establish the asymptotic validity of 
Markov chain CLT based standard error estimates. However, if quantitative 
bounds on the distance to stationarity are needed, then an additional 
{\it minorization condition} needs to be established. For the sake of 
completeness, we derive such a condition in \hyperref[Appendix C]{Appendix \ref{Appendix C}} (see Lemma \ref{Theorem2}).

\subsection{A simulation study} \label{simulation:original:Horseshoe}

\noindent
 The objective of this study is to examine the practical feasibility/scalability  of the  Gibbs sampler described and analyzed in Section~\ref{Horseshoe:Gibbs:sampler} by comparing its computational performance with the JOB Gibbs sampler. We have considered two different simulation settings. For the first simulation setting, we fix the sample size $n$ to be $500$ and the number of predictors $p$ to be $1000$. The  first $20$ entries of the ``true" regression coefficient vector  ${\boldsymbol \beta}^0 := (\beta^0_1, \ldots , \beta^0_{p})$ are specified as  $ \beta^0_j=2^{s_j}$ where $s_j's$ are a sequence of equally spaced values in the interval $(-3.5,3)$, and the 
 other entries are set to zero. The entries of the design matrix $\X$ are generated independently from $\mathcal{N}(0,1)$. Then, we generate the response vector ${\y}$ from the model ${\y}= \X {\boldsymbol \beta}^0 + {\boldsymbol \epsilon}$ where the error vector ${\boldsymbol \epsilon}$ has i.i.d. normal entries with mean $0$ and standard deviation $0.1$. For the second simulation setting, the same procedure described 
 above is followed with $n = 750$ and $p = 1500$. 

\begin{figure}[hbt!]
     \centering
     \begin{subfigure}[b]{0.8\textwidth}
         \centering
         \includegraphics[width=\textwidth]{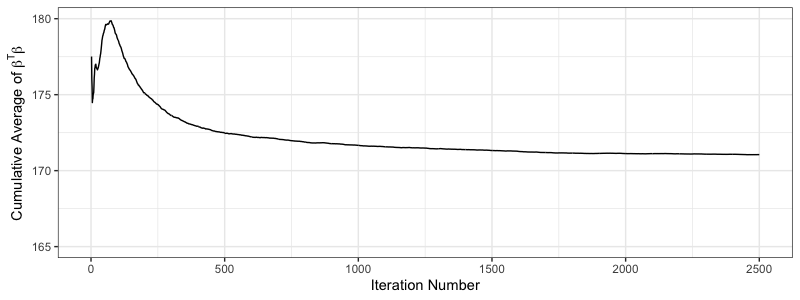}
         \caption{Cumulative Average plot  for JOB  sampler.}
         \label{fig:JOB_p_1000}
     \end{subfigure}\\
     \vspace{.1in}
     \begin{subfigure}[b]{0.8\textwidth}
         \centering
         \includegraphics[width=\textwidth]{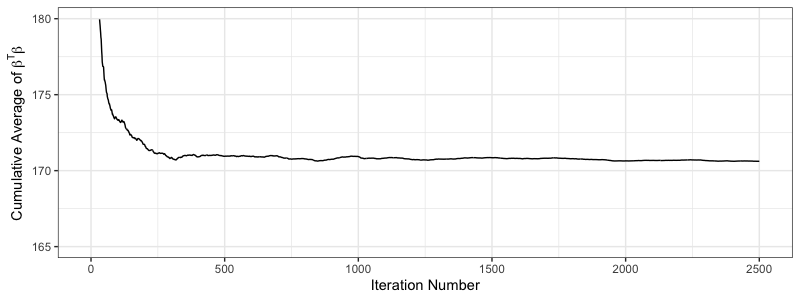}
         \caption{Cumulative Average for proposed sampler.}
         \label{fig:prop_p_1000}
     \end{subfigure}
     \hfill
     \caption{Cumulative average plots for the function $\boldsymbol{\beta}^T\boldsymbol{\beta}$ corresponding to a randomly 
     selected data set in the $n=500, p=1000$ simulation setting.}
     \label{fig:p_1000}
  \end{figure}

\begin{figure}[hbt!]
     \centering
     \begin{subfigure}[hbt!]{0.8\textwidth}
         \centering
         \includegraphics[width=\textwidth]{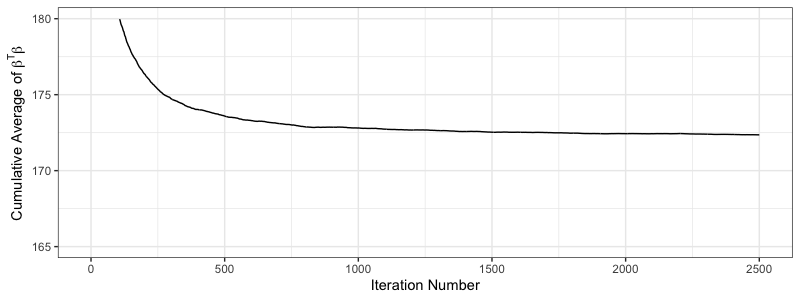}
         \caption{Cumulative Average plot  for JOB  sampler.}
         \label{fig:JOB_p_1500}
     \end{subfigure}\\
     \vspace{.1in}
     \begin{subfigure}[hbt!]{0.8\textwidth}
         \centering
         \includegraphics[width=\textwidth]{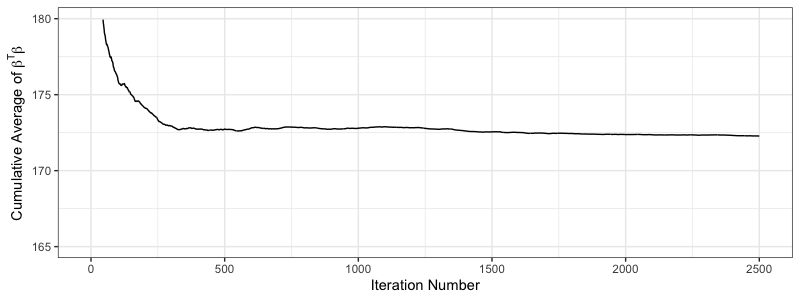}
         \caption{Cumulative Average for proposed sampler.}
         \label{fig:prop_p_1500}
     \end{subfigure}
     \hfill
     \caption{Cumulative average plots for the function $\boldsymbol{\beta}^T\boldsymbol{\beta}$ corresponding to a randomly 
     selected data set in the $n=750, p=1500$ simulation setting.}
       \label{fig:p_1500}
  \end{figure}

We generate $20$ data sets each from each of the two simulation settings, and run both the Gibbs samplers on each of the $40$ data sets. For a fair comparison, both algorithms were implemented in $R$. The simulations were run on a machine with a 64 bit Windows 7 operating system, 8 GB RAM and a 
3.4 GHz processor. We provide the run-times for generating $2500$ iterations from each of the Gibbs sampler in Table~\ref{tab:table1}.  In the case $n=750,\  p=1500$,  the average CPU time required for JOB sampler is $ 9213.84 $ seconds compared to the average of $9114.19$ seconds for the proposed  sampler.  In  the case $ n=500,\  p=1000$ ,  the average required times are $2446.9$ and $2391.1$ seconds respectively for the JOB sampler and the proposed sampler. In order to check the convergence of the MCMC chains,  we considered the cumulative average plots of the function $\boldsymbol{\beta}^T\boldsymbol{\beta}$. These plots for a randomly selected data set from each of the two simulation settings are provided in  Figure~\ref{fig:p_1000} and   Figure~\ref{fig:p_1500}.  The plots for all the other Markov chains are similar to the ones presented here. 

It is evident from the above results that the Gibbs sampler analysed in this paper has comparable (slightly better) computational performance than the JOB Gibbs sampler in the above settings, and hence is practically useful. The geoemtric ergodicity result in Theorem \ref{Theorem 2.1} therefore helps provide asymptotically valid standard error estimates for corresponding MCMC approximations, under weaker assumptions compared to the JOB Gibbs sampler. 

In \cite{johndrow2017bayes}, the authors discuss a time-inhomogeneous approximation/modification to the JOB Gibbs sampler, termed as the approximate Gibbs sampler, for faster and more scalable computation. We would like to mention that the exact same modifications can be used for the Gibbs sampler described in 
Section \ref{Horseshoe:Gibbs:sampler} to obtain a corresponding approximate faster 
and time-inhomogeneous version. 

\begin{table}[hbt!]
\begin{subtable}{0.45\textwidth}
\centering
\begin{tabular}{|@{}cccc@{}|}
\toprule
 & \multicolumn{2}{l}{Simulation setting: n=500, p=1000} &  \\
 \midrule
 & JOB sampler   & Proposed Sampler  &  \\
 & 2300.2        & 2236.7            &  \\
 & 2303.5        & 2247.33           &  \\
 & 2662.55       & 2312.96           &  \\
 & 2700.01       & 2669.57           &  \\
 & 2593.41       & 2600.05           &  \\
 & 2675.34       & 2569.93           &  \\
 & 2755.46       & 2763.22           &  \\
 & 2528.5        & 2447.42           &  \\
 & 2641.8        & 2580.94           &  \\
 & 2664.14       & 2593.7            &  \\
 & 2360.96       & 2593.69           &  \\
 & 2313.37       & 2248.14           &  \\
 & 2301.66       & 2240.99           &  \\
 & 2298.01       & 2248.02           &  \\
 & 2305.22       & 2254.25           &  \\
 & 2294.81       & 2242.13           &  \\
 & 2313.77       & 2244.53           &  \\
 & 2290.15       & 2240.65           &  \\
 & 2307.51       & 2242.81           &  \\
 & 2327.64       & 2244.93           & \\
 \bottomrule
\end{tabular}
\caption{}
\label{tab:table1_a}
\end{subtable}%
\hfill
\begin{subtable}{0.45\textwidth}
\centering 
\begin{tabular}{|@{}cccc@{}|}
\toprule
 & \multicolumn{2}{l}{Simulation setting: n=750, p=1500} &  \\ \midrule
 & JOB sampler   & Proposed Sampler  &  \\
 & 9211.42       & 9092.32           &  \\
 & 9192.64       & 9101.27           &  \\
 & 9236.63       & 9105.14           &  \\
 & 9190.44       & 9125.47           &  \\
 & 9241.26       & 9128.65           &  \\
 & 9204.89       & 9108.78           &  \\
 & 9222.37       & 9100.23           &  \\
 & 9225.12       & 9137.84           &  \\
 & 9196.45       & 9110.89           &  \\
 & 9197.96       & 9135.35           &  \\
 & 9204.08       & 9110.05           &  \\
 & 9213.73       & 9126.43           &  \\
 & 9215.79       & 9117.76           &  \\
 & 9216.99       & 9108.28           &  \\
 & 9212.47       & 9098.57           &  \\
 & 9220          & 9115.62           &  \\
 & 9232.6        & 9125.46           &  \\
 & 9196.04       & 9112.85           &  \\
 & 9233.48       & 9122.75           &  \\
 & 9212.5        & 9100.18           &  \\ \bottomrule
\end{tabular}
\caption{}
\label{tab:table1_d}
\end{subtable}%
\caption{\footnotesize The run-times (in seconds) required to generate 2500 MCMC samples using the JOB sampler and the proposed sampler for $40$ simulated datasets are tabulated. The left  sub-table (\subref{tab:table1_a}) corresponds to the simulation setting $n=500, p=1000$,  while the right sub-table (\subref{tab:table1_d}) corresponds to the simulation setting $n=750, p=1500$. }
\label{tab:table1}
\end{table}

\section{Geometric ergodicity for regularized Horseshoe Gibbs samplers}

\subsection{A Gibbs sampler for the regularized Horseshoe }
\label{Regularized:Gibbs:sampler}

\noindent
Recall from the introduction that the regularized Horseshoe prior  
developed in \Myciteauthor{piironen2017} is given by 
\begin{eqnarray}
\left.\beta_i \right. \mid \lambda_i^2, \sigma^2,\tau^2\sim \mathcal{N}_p \left( 0, \left( \frac{1}{c^2} + \frac{1}{\lambda_i^2 \tau^2} \right)^{-1} \sigma^2 \right) \text{independently for}\ i=1,2,\cdots,p \nonumber\\
\left.\lambda_i\right.\sim C^+ (0,1) \text{independently for}\ i=1,2,\cdots,p \nonumber\\
\tau^2\sim \pi_\tau (\cdot) \hspace{0.2in} \sigma^2\sim 
\text{Inverse-Gamma}(a,b) \label{regularized_prior}
\end{eqnarray}

\noindent
The only difference between this prior and the original Horseshoe prior 
in (\ref{horseshoe_prior}) is the additional regularization introduced in 
the the prior conditional variance of the $\beta_i$s through the constant $c$. As $c \rightarrow \infty$ in (\ref{regularized_prior}), then one reverts back to the 
original Horseshoe specification in (\ref{horseshoe_prior}). 

Note that one of the salient features of the Horseshoe prior is the lack of shrinkage/regularization of parameter values that are far away from zero. The authors in \cite{piironen2017} argue that while this feature is one of the key strengths of the Horseshoe prior in many situations, it can be a drawback in settings where the parameters are weakly identified. We refer the reader to \cite{piironen2017} for a thorough motivation and discussion of the properties and performance of this prior vis-a-vis the Horseshoe prior. Our focus in this paper is to look at Markov chains to sample from the resulting intractable regularized Horseshoe posterior, and investigate properties such as geometric ergodicity. 

The authors in \cite{piironen2017} use Hamiltonian Monte Carlo (HMC) to generate samples from the posterior distribution. Geometric ergodicity of this HMC chain, however, is not established. In recent work \cite{livingstone2019}, sufficient conditions for geometric ergodicity (or lack thereof) for general HMC chains have been provided. However, these conditions, namely Assumptions A1, A2, A3 in \cite{livingstone2019}, are rather complex and intricate, and at least to the best of our understanding it is unclear and hard to verify if these conditions are satisfied by the HMC chain in \cite{piironen2017}. 

Given the host of Gibbs samplers available in the literature for the original Horseshoe posterior, it is natural to consider a Gibbs sampler to sample from the regularized Horseshoe posterior as well. In fact, after introducing the augmented variables $\left\{\nu_j\right\}_{j=1}^p$, the following conditional posterior distributions can be obtained after straightforward computations:
\begin{eqnarray}\label{reghorseshoeposteriors}
\left.\bb\right|\sigma^2,\tau^2,\bl,\y\sim\mathcal{N}\left(A_c^{-1}\X^T\y,\sigma^2A_c^{-1}\right)\nonumber \\
\left.\sigma^2\right|\tau^2,\bl,\y\sim \text{Inverse-Gamma}\left(a+\frac{n}{2},\frac{\y^T\left(I_n-\X A_c^{-1}\X^T\right)\y}{2}+b\right)\nonumber \\
\left.\nu_j\right|\lambda_j^2,\y\sim\text{Inverse-Gamma}\left(1,1+\frac{1}{\lambda_j^2}\right),\ \text{independently for}\ j=1,2,\cdots,p\nonumber \\
 \pi \left(\left.\mathbf{\bl}\right|\boldsymbol{\beta,\nu},\sigma^2,\tau^2,\mathbf{\y}\right)=\prod\limits_{j=1}^p g\left(\left.\lambda_j^2\right|\nu_j,\beta_j,\sigma^2,\tau^2,\y\right) 
 \nonumber\\
\left.\tau^2\right|\bl,\y\sim\pi\left(\left.\tau^2\right|\bl,\y\right)
\label{regularized_posterior}
\end{eqnarray}

\noindent
where 
$$
g\left(\left.\lambda_j^2\right|\nu_j,\beta_j,\sigma^2,\tau^2,\y\right) \propto \left(\frac{1}{c^2}+\frac{1}{\tau^2\lambda_j^2}\right)^{\frac{1}{2}}\left(\lambda_j^2\right)^{-\frac{3}{2}}\exp{\left[-\frac{1}{\lambda_j^2}\left(\frac{1}{\nu_j}+\frac{\beta_j^2}{2\sigma^2\tau^2}\right)\right]}
$$
   
\noindent
for $j=1,2,\cdots,p$, 
$$
\pi\left(\left.\tau^2\right|\bl,\y\right)\propto\left|A_c\right|^{-\frac{1}{2}}\prod\limits_{j=1}^p\left\{\left(\frac{1}{c^2}+\frac{1}{\tau^2\lambda_j^2}\right)^{\frac{1}{2}}\right\}\left(\frac{\y^T\left(I_n-\X A_c^{-1}\X^T\right)\y}{2}+b\right)^{-\left(a+\frac{n}{2}\right)}\pi_\tau \left(\tau^2\right)
$$

\noindent
and $A_c=\X^T\X+\left(\tau^2\Bl\right)^{-1}+c^{-2}I_p$. Most of the above densities are standard and can be easily sampled from. Efficient rejection/Metropolis samplers can be used to sample from the one-dimensional non-standard densities $g\left(\left.\lambda_j^2\right|\nu_j,\beta_j,\sigma^2,\tau^2,\y\right)$ and $\pi\left(\left.\tau^2\right|\bl,\y\right)$ (see in Appendix \ref{Appendix D}). Hence, a two-block Gibbs sampler, whose one step-transition from $({\boldsymbol \beta}_0, \sigma^2_0, {\boldsymbol \nu}_0, \tau^2_0, {\boldsymbol \lambda}_0)$ to 
$({\boldsymbol \beta}, \sigma^2, {\boldsymbol \nu}, \tau^2, {\boldsymbol \lambda})$ is given by sampling sequentially from $\pi({\boldsymbol \beta}, \sigma^2, {\boldsymbol \nu}, \tau^2 \mid {\boldsymbol \lambda}_0, \y)$ and 
$\pi({\boldsymbol \lambda} \mid {\boldsymbol \beta}, \sigma^2, {\boldsymbol \nu}, \tau^2, \y)$, can be used to generate approximate samples from the regularized Horseshoe posterior. We will denote the transition kernel of this two-block Gibbs sampler by $K_{aug,reg}$ (analogous to $K_{aug}$ in the original Horseshoe setting). 

Our goal now is to establish geometric ergodicity for $K_{aug,reg}$. We 
will achieve this by focusing on the marginal $\bl$-chain corresponding to 
$K_{aug,reg}$. The one-step transition dynamics of this Markov chain from 
$\bl_m$ to $\bl_{m+1}$ is given as follows:\\
\begin{enumerate}
    \item Draw $\tau^2$ from $\pi\left(\left.\tau^2\right|\mathbf{\bl}_m,\mathbf{\y}\right)$
    \item Draw $\boldsymbol{\nu}$ from $\pi\left(\left.\boldsymbol{\nu}\right|\mathbf{\bl}_m,\mathbf{\y}\right)=\prod\limits_{j=1}^p\text{Inverse-Gamma}\left(1,1+\frac{1}{\lambda_{j;m}^2}\right)$
    \item Draw $\sigma^2$ from $\pi\left(\left.\sigma^2\right|\tau^2,\mathbf{\bl}_m,\mathbf{\y}\right)=\text{Inverse-Gamma}\left(a+\frac{n}{2},\frac{\mathbf{\y}^T\left(I_n-\X A_c^{-1}X^T\right)\mathbf{\y}}{2}+b\right)$
    \item Draw $\boldsymbol{\beta}$ from $\pi\left(\left.\boldsymbol{\beta}\right|\sigma^2,\tau^2,\mathbf{\bl}_m,\mathbf{\y}\right)=\mathcal{N}_p(A_c^{-1}\X^T\y,\sigma^2A_c^{-1})$
    \item Finally draw $\mathbf{\bl}_{m+1}$ from $\pi\left(\left.\mathbf{\bl}\right|\boldsymbol{\beta,\nu},\sigma^2,\tau^2,\mathbf{\y}\right)=\prod\limits_{j=1}^p g\left(\left.\lambda_j^2\right|\nu_j,\beta_j,\sigma^2,\tau^2,\y\right)$. 
\end{enumerate}

The Markov transition density (MTD) corresponding to the marginal $\bl$-chain is given by 
\begin{eqnarray}\label{MTDreghorseshoe}
     k\left(\mathbf{\bl}_0,\mathbf{\bl}\right) &=& \int_{\mathbb{R}_+}\int_{\mathbb{R}_{+}}\int_{\mathbb{R}^p}\int_{\mathbb{R}_+^p}\pi\left(\left.\mathbf{\bl}\right|\boldsymbol{\beta,\nu},\sigma^2,\tau^2,\mathbf{\y}\right)\pi\left(\left.\boldsymbol{\beta,\nu},\sigma^2,\tau^2\right|\mathbf{\bl}_0,\mathbf{\y}\right)d\boldsymbol{\nu}d\boldsymbol{\beta}d\sigma^2d\tau^2 \nonumber\\
   &=& \int_{\mathbb{R}_+}\int_{\mathbb{R}_{+}}\int_{\mathbb{R}^p}\int_{\mathbb{R}^{p}_+}\pi\left(\left.\mathbf{\bl}\right|\boldsymbol{\beta,\nu},\sigma^2,\tau^2,\mathbf{\y}\right)\pi\left(\left.\boldsymbol{\beta}\right|\sigma^2,\tau^2,\mathbf{\bl}_0,\mathbf{\nu,\y}\right)\nonumber \\
    &\qquad \qquad \times& \pi\left(\left.\sigma^2\right|\tau^2,\mathbf{\bl}_0,\mathbf{\nu,\y}\right)\pi\left(\left.\boldsymbol{\nu}\right|\mathbf{\bl}_0,\tau^2,\mathbf{\y}\right)\pi\left(\left.\tau^2\right|\mathbf{\bl}_0,\mathbf{\y}\right)d\boldsymbol{\nu}d\boldsymbol{\beta}d\sigma^2d\tau^2
\end{eqnarray}

\subsection{Drift and minorization analysis for the regularized Horseshoe 
$\bl$-chain} \label{Regularized:drift:minorization}

\noindent
The geometric ergodicity of the $\bl$-chain will be established using a 
drift and minorization analysis. However, given the modifications in the 
regularized Horseshoe posterior, the drift function $V(\bl)$ (see 
(\ref{driftfunction})) used for the original Horseshoe does not work in 
this case. We will instead use another drift function $\Tilde{V}(\bl)$ 
defined by 
\begin{equation}\label{driftfunctionreghorseshoe}
    \Tilde{V}\left(\mathbf{\bl}\right)=\sum\limits_{j=1}^p\left(\lambda_j^2\right)^{-\frac{\delta}{2}};\ \text{for some constant}\ \delta\in(0,1).
    \end{equation}

\noindent
As discussed previously, the function $V(\bl)$ is unbounded off petite sets and the $V$-based drift condition in Lemma \ref{Theorem1} is enough to guarantee geometric ergodicity for the original Horseshoe Gibbs sampler 
$K_{aug}$. A minorization condition is only needed if one also wants to get quantitative convergence bounds for distance to stationarity. The function 
$\Tilde{V}$ however, is {\it not} unbounded off petite sets since 
$$
B \left(\Tilde{V}, d \right) = \left\{ \bl \in \mathbb{R}_+^p: 
\Tilde{V}(\bl) \leq d \right\}
$$

\noindent
is not a compact subset of $\mathbb{R}_+^p$ for $d > 0$. Hence, a drift 
condition with $\Tilde{V}$ needs to be complemented with a minorization condition in order to establish geometric ergodicity (Theorem \ref{thm:regularized:Horseshoe:Piironen:Vehtari}). We establish these 
two conditions respectively in Sections \ref{driftreghorseshoe} and \ref{reghorseshoeminorization} below. As opposed to the original Horseshoe setting, 
we do not require that the prior density $\pi_\tau$ is truncated below away from zero. Only the existence of the $\delta/2^{th}$-moment is assumed for some 
$\delta \in (0.00162, 0.22176)$: a very mild condition, satisfied for example by the commonly used half-Cauchy density. 

\subsubsection{Drift condition}\label{driftreghorseshoe}

\begin{lem}\label{Theorem1reghorseshoe}
Suppose $\int_{\mathbb{R}_+} u^{\delta/2} \pi_\tau (u) du < \infty$ for some $\delta \in (0.00162, 0.22176)$. Then, there exist constants $0<\gamma^*=\gamma^*\left(\delta\right)<1$ and $b^*<\infty$ such that
\begin{equation}\label{driftconditionreghorseshoe}
    \E\left[\left.\Tilde{V}\left(\mathbf{\bl}\right)\right|\mathbf{\bl}_0\right]\leq \gamma^* \Tilde{V}\left(\mathbf{\bl}_0\right)+b^*
\end{equation}
for every $\bl_0\in \mathbb{R}_+^p$. \\
\end{lem}

\proof Note that by linearity 
\begin{equation} \label{reghorseshoeexpectations:linearity}
\E\left[\left.\Tilde{V}\left(\mathbf{\bl}\right)\right|\mathbf{\bl}_0\right] = \sum_{j=1}^p \E\left[\left. \left(\lambda_j^2\right)^{-\frac{\delta}{2}} \right|\mathbf{\bl}_0\right]
\end{equation}

\noindent
Fix $j\in\left\{1,2,\cdots,p\right\}$ arbitrarily. It follows from the definition of the MTD \eqref{MTDreghorseshoe} that 
\begin{eqnarray} \label{reghorseshoeexpectations:iterated}
& & \E\left[\left.\left(\bl_j^2\right)^{-\frac{\delta}{2}}\right|\mathbf{\bl}_0\right] \nonumber\\
&=& \E\left[\left.\E\left[\left.\E\left[\left.\E\left[\left.\E\left[\left.\left(\bl_j^2\right)^{-\frac{\delta}{2}}\right.\mid\boldsymbol{\beta,\nu},\sigma^2,\tau^2,\mathbf{\y}\right]\right.\mid\sigma^2,\tau^2,\mathbf{\bl}_0,\mathbf{\nu,\y}\right]\right.\mid\tau^2,\bl_0,\mathbf{\nu},\y\right]\right.\mid\tau^2\mathbf{\bl}_0,\mathbf{\y}\right]\right.\mid\bl_0,\y\right]. \nonumber\\
& & 
\end{eqnarray}

\noindent
We begin by evaluating the innermost expectation. It follows by using $\sqrt{u+v}\leq\sqrt{u}+\sqrt{v}$ for $u,v\geq 0$ that \\

\begin{eqnarray}\label{driftreghorseshoeproof}
& & \E\left[\left.\left(\lambda_j^2\right)^{-\frac{\delta}{2}}\right|\boldsymbol{\bb,\nu},\sigma^2,\tau^2,\y\right] \nonumber\\
&=& \frac{\int\limits_0^\infty\left(\lambda_j^2\right)^{-\frac{\delta}{2}}\left(\frac{1}{c^2}+\frac{1}{\tau^2\lambda_j^2}\right)^{\frac{1}{2}}\left(\lambda_j^2\right)^{-\frac{3}{2}}\exp{\left[-\frac{1}{\lambda_j^2}\left(\frac{1}{\nu_j}+\frac{\beta_j^2}{2\sigma^2\tau^2}\right)\right]}d\lambda_j^2}{\int\limits_0^\infty\left(\frac{1}{c^2}+\frac{1}{\tau^2\lambda_j^2}\right)^{\frac{1}{2}}\left(\lambda_j^2\right)^{-\frac{3}{2}}\exp{\left[-\frac{1}{\lambda_j^2}\left(\frac{1}{\nu_j}+\frac{\beta_j^2}{2\sigma^2\tau^2}\right)\right]}d\lambda_j^2}\nonumber \\
&\leq& \frac{\int\limits_0^\infty\left(\lambda_j^2\right)^{-\frac{\delta}{2}}\left(\frac{1}{|c|}+\frac{1}{\sqrt{\tau^2\lambda_j^2}}\right)\left(\lambda_j^2\right)^{-\frac{3}{2}}\exp{\left[-\frac{1}{\lambda_j^2}\left(\frac{1}{\nu_j}+\frac{\beta_j^2}{2\sigma^2\tau^2}\right)\right]}d\lambda_j^2}{\int\limits_0^\infty\left(\frac{1}{c^2}+\frac{1}{\tau^2\lambda_j^2}\right)^{\frac{1}{2}}\left(\lambda_j^2\right)^{-\frac{3}{2}}\exp{\left[-\frac{1}{\lambda_j^2}\left(\frac{1}{\nu_j}+\frac{\beta_j^2}{2\sigma^2\tau^2}\right)\right]}d\lambda_j^2}\nonumber \\
&\leq& \frac{\frac{\left(\tau^2\right)^{\frac{\delta}{2}}}{|c|}\int\limits_0^\infty\left(\tau^2\lambda_j^2\right)^{-\frac{\delta}{2}}\left(\lambda_j^2\right)^{-\frac{3}{2}}\exp{\left[-\frac{1}{\lambda_j^2}\left(\frac{1}{\nu_j}+\frac{\beta_j^2}{2\sigma^2\tau^2}\right)\right]}d\lambda_j^2}{\int\limits_0^\infty\left(\frac{1}{c^2}+\frac{1}{\tau^2\lambda_j^2}\right)^{\frac{1}{2}}\left(\lambda_j^2\right)^{-\frac{3}{2}}\exp{\left[-\frac{1}{\lambda_j^2}\left(\frac{1}{\nu_j}+\frac{\beta_j^2}{2\sigma^2\tau^2}\right)\right]}d\lambda_j^2}\nonumber \\
&\qquad \qquad \qquad +& \frac{\int\limits_0^\infty\left(\lambda_j^2\right)^{-\frac{\delta}{2}}\left(\lambda_j^2\right)^{-2}\exp{\left[-\frac{1}{\lambda_j^2}\left(\frac{1}{\nu_j}+\frac{\beta_j^2}{2\sigma^2\tau^2}\right)\right]}d\lambda_j^2}{\int\limits_0^\infty\left(\lambda_j^2\right)^{-2}\exp{\left[-\frac{1}{\lambda_j^2}\left(\frac{1}{\nu_j}+\frac{\beta_j^2}{2\sigma^2\tau^2}\right)\right]}d\lambda_j^2}. 
\end{eqnarray}

\noindent
The first term in the last inequality of (\ref{driftreghorseshoeproof}) can be expressed as 
$$
\frac{\left(\tau^2\right)^{\frac{\delta}{2}}}{|c|}\frac{\E\left[X^\delta\right]}{\E\left[\sqrt{\frac{1}{c^2}+X^2}\right]}
$$

\noindent
where $\tau^2X^2\sim\text{Gamma}\left(\frac{1}{2},\frac{1}{\nu_j}+\frac{\beta_j^2}{2\sigma^2\tau^2}\right)$. Using Young's inequality, it follows that the first term is bounded above by 
$$
\frac{\text{max}\left\{1,|c|\right\}}{|c|}\sqrt{\delta}\left(\tau^2\right)^{\frac{\delta}{2}}. 
$$

\noindent
The second term in the last inequality of (\ref{driftreghorseshoeproof}) is basically an Inverse-Gamma expectation, and is exactly equal to 
$$
\frac{\Gamma\left(1+\frac{\delta}{2}\right)}{\left(\frac{1}{\nu_j}+\frac{\beta_j^2}{2\sigma^2\tau^2}\right)^{\frac{\delta}{2}}}. 
$$

\noindent
Hence, we get 
$$
\E\left[\left.\left(\lambda_j^2\right)^{-\frac{\delta}{2}}\right|\boldsymbol{\bb,\nu},\sigma^2,\tau^2,\y\right] \leq \frac{\text{max}\left\{1,|c|\right\}}{|c|}\sqrt{\delta}\left(\tau^2\right)^{\frac{\delta}{2}} + \frac{\Gamma\left(1+\frac{\delta}{2}\right)}{\left(\frac{1}{\nu_j}+\frac{\beta_j^2}{2\sigma^2\tau^2}\right)^{\frac{\delta}{2}}}. 
$$

\noindent
Note that the conditional distribution of $\beta_j$ given $\sigma^2,\tau^2,\mathbf{\bl}_0,{\boldsymbol \nu},\mathbf{\y}$ is a Gaussian distribution with variance $\sigma_j^2\overset{\text{def}}{=}\sigma^2\boldsymbol{e}_J^TA_c^{-1}\boldsymbol{e}_j\geq \sigma^2\left(\bar{\omega}+\frac{1}{c^2} +\frac{1}{\tau^2\lambda_{j;0}^2}\right)^{-1}$. Here $\bar{\omega}$ is the maximum eigenvalue of $\X^T\X$. Now, proceeding exactly with the analysis 
from (\ref{Gaussian:moment}) to (\ref{drift:1}) in the proof of 
Lemma \ref{Theorem1} with $\bar{\omega}$ replaced by $\bar{\omega} + c^{-2}$ and using Proposition \ref{reghsupperboundontau^2} instead of 
Proposition \ref{upperboundontau^2} yields 
$$
\E\left[\left.\sum_{j=1}^p\left(\lambda_j^2\right)^{-\frac{\delta_1}{2}}\right|\mathbf{\bl}_0\right]\leq \gamma^* \left(\delta\right)\sum_{j=1}^p\left(\lambda_{j;0}^2\right)^{-\frac{\delta_1}{2}}+b^*
$$

\noindent
with 
$$
\gamma^*\left(\delta\right)=\Gamma\left(1+\frac{\delta}{2}\right)\left(\frac{1}{\left\{\Gamma\left(1-\frac{\delta}{2}\right)\right\}^{\frac{2}{\delta}}}+\frac{\sqrt{\pi}^{\frac{2}{\delta}}}{\left\{\Gamma\left(\frac{1-\delta}{2}\right)\right\}^{\frac{2}{\delta}}}\right)^{-\frac{\delta}{2}} 
$$

\noindent
and 
$$
b^*=p\frac{\text{max}\left\{1,|c|\right\}}{|c|}\sqrt{\delta}C_2+p\cdot\text{max}\left\{1,(\bar{\omega}+c^{-2})^{\frac{\delta}{2}}C_2\right\}\cdot\gamma^*\left(\delta\right). 
$$ 

\noindent
Here $C_2$ is as in Proposition \ref{reghsupperboundontau^2}. It can be shown that $\gamma^* \left( \delta \right) < 1$ for $\delta\in (0.00162, 0.22176) $. \\
Hence, the required geometric drift condition has been established. \qed 
\subsubsection{Minorization condition}\label{reghorseshoeminorization}

\noindent
As discussed previously, the drift function $\Tilde{V}$ is not unbounded off compact sets, and the drift condition in Lemma \ref{Theorem1reghorseshoe} needs to be complemented by an associated minorization condition to establish geometric ergodicity. Fix a $d>0$. Define
\begin{equation}\label{reghorseshoeminorizationset}
    B \left( \Tilde{V},d \right) = \left\{\bl\in\mathbb{R}_+^p:\Tilde{V}\left(\bl\right)\leq d\right\}
\end{equation}
We now establish the following minorization condition associated to the geometric drift condition in Lemma \ref{Theorem1reghorseshoe}.
\begin{lem}\label{Theorem2reghorseshoe}
There exists a constant $\epsilon^*=\epsilon^*\left(\Tilde{V},d\right)>0$ and a density function $h$ on $\mathbb{R}_+^p$ such that
\begin{equation}\label{reghorseshoeminorizationcondition}
   k\left(\bl_0,\bl\right)\geq \epsilon^* h\left(\bl\right) 
\end{equation}
for every $\bl_0\in B\left(\Tilde{V},d\right)$.
\end{lem}
\proof Fix a $\bl_0\in B\left(\Tilde{V},d\right)$ arbitrarily. In order to prove \eqref{reghorseshoeminorizationcondition} we will demonstrate appropriate lower bounds for the conditional densities appearing in \eqref{MTDreghorseshoe}. From \eqref{reghorseshoeposteriors} we have the following:

\begin{eqnarray*}
 \pi\left(\left.\tau^2\right|\bl_0,\y\right)&\geq& b^{a+\frac{n}{2}} \omega_*^{-\frac{p}{2}}\left(1+\frac{1}{\tau^2}\right)^{-\frac{p}{2}}|c|^{-p}\left(\y^T\y+b\right)^{-\left(a+\frac{n}{2}\right)}\pi_\tau\left(\tau^2\right);
\end{eqnarray*}

\noindent
where $\omega_*=\text{max}\left\{\bar{\omega}+c^{-2},d^{\frac{2}{\delta}}\right\}$; (recall that $\bar \omega$ denotes the maximum eigenvalue of $\X^T\X$),
\begin{eqnarray*}
 \pi\left(\left.\boldsymbol{\beta}\right|\sigma^2,\tau^2,\bl_0,\y\right)&\geq& \left(2\pi\sigma^2\right)^{-\frac{p}{2}}\left|c\right|^{-p}\\
 &\times&\exp{\left[-\frac{\left(\boldsymbol{\beta}-\Omega^{-1}\X^T\y\right)^T\Omega\left(\boldsymbol{\beta}-\Omega^{-1}\X^T\y\right)+\y^T\X\left(c^2I_p-\Omega^{-1}\right)\X^T\y}{2\sigma^2}\right]}; 
\end{eqnarray*}

\noindent
where $\Omega=\omega_*\left(1+\frac{1}{\tau^2}\right)I_p$, 
\begin{eqnarray*}
  \pi\left(\left.\bl\right|\boldsymbol{\beta,\nu},\sigma^2,\tau^2,\y\right)&=& \prod_{j=1}^p\left\{\frac{\left(\frac{1}{c^2}+\frac{1}{\tau^2\lambda_j^2}\right)^{\frac{1}{2}}\left(\lambda_j^2\right)^{-\frac{3}{2}}\exp{\left[-\frac{1}{\lambda_j^2}\left(\frac{1}{\nu_j}+\frac{\beta_j^2}{2\sigma^2\tau^2}\right)\right]}}{\int\limits_0^\infty\left(\frac{1}{c^2}+\frac{1}{\tau^2\lambda_j^2}\right)^{\frac{1}{2}}\left(\lambda_j^2\right)^{-\frac{3}{2}}\exp{\left[-\frac{1}{\lambda_j^2}\left(\frac{1}{\nu_j}+\frac{\beta_j^2}{2\sigma^2\tau^2}\right)\right]}d\lambda_j^2}\right\}\\
  &\geq& \prod\limits_{j=1}^p\left\{\frac{\left(\tau^2\right)^{-\frac{1}{2}}\left(\lambda_j^2\right)^{-2}\exp{\left[-\frac{1}{\lambda_j^2}\left(\frac{1}{\nu_j}+\frac{\beta_j^2}{2\sigma^2\tau^2}\right)\right]}}{k^*\left(\sqrt{\nu_j}+\frac{\sigma^2\sqrt{\tau^2}}{\beta_j^2}\right)}\right\};
\end{eqnarray*}

\noindent
where $k^*=\text{max}\left\{\sqrt{\pi}|c|^{-1},2\right\}$, and 
\begin{eqnarray} \label{lowerboundsreghorseshoe}
   \pi\left(\left.\boldsymbol{\nu}\right|\bl_0,\y\right)&\geq& \prod\limits_{j=1}^p\left\{\nu_j^{-2}\exp{\left[-\frac{1}{\nu_j}\left(1+d^{\frac{2}{\delta}}\right)\right]}\right\}\nonumber \\
\pi\left(\left.\sigma^2\right|\tau^2,\bl_0,\y\right)&\geq& \frac{b^{a+\frac{n}{2}}}{\Gamma\left(a+\frac{n}{2}\right)}\left(\sigma^2\right)^{-\left(a+\frac{n}{2}\right)-1}\exp{\left[-\frac{1}{\sigma^2}\left(\frac{\y^T\y}{2}+b\right)\right]}. 
\end{eqnarray}

\noindent
Combining all the lower bounds provided above, it follows from \eqref{MTDreghorseshoe} that
\begin{align*}\label{minorizationproofreghorseshoe}
& \quad k\left(\bl_0,\bl\right)\nonumber \\
& \geq \frac{\left(2\pi\right)^{-\frac{p}{2}}b^{2\left(a+\frac{n}{2}\right)}\left(\sqrt{\omega_*}k^*c^2\right)^{-p}}{\left(\y^T\y+b\right)^{a+\frac{n}{2}}\Gamma\left(a+\frac{n}{2}\right)}\int\limits_{\mathbb{R}_+}\int\limits_{\mathbb{R}_+}\int\limits_{\mathbb{R}^p}\int\limits_{\mathbb{R}^p_+}\left(\sigma^2\right)^{-\left(a+\frac{n+p}{2}\right)-1}\prod\limits_{j=1}^p\left\{\frac{\nu_j^{-2}\exp{\left[-\frac{1}{\nu_j}\left(1+d^{\frac{2}{\delta}}+\frac{1}{\lambda_j^2}\right)\right]}}{\sqrt{\nu_j}+\frac{\sigma^2\sqrt{\tau^2}}{\beta_j^2}}\right\}\nonumber \\ 
& \qquad \qquad \exp{\left[-\frac{\left(\bb-\Omega^{-1}\X^T\y\right)^T\Omega\left(\bb-\Omega^{-1}\X^T\y\right)+\bb^T\left(\tau^2\Bl\right)^{-1}\bb}{2\sigma^2}\right]}\prod\limits_{j=1}^p\left\{\left(\lambda_j^2\right)^{-2}\right\}\nonumber \\
& \qquad \qquad \exp{\left[-\frac{1}{\sigma^2}\left(\frac{\y^T\y+\y^T\X\left(c^2I_p-\Omega^{-1}\right)\X^T\y}{2}+b\right)\right]}\left(1+\tau^2\right)^{-\frac{p}{2}}\pi_\tau\left(\tau^2\right)d\boldsymbol{\nu}d\bb d\sigma^2 d\tau^2\nonumber \\
\end{align*}

\noindent
Now for the inner most integral wrt $\boldsymbol \nu$, substituting the lower bounds given in Proposition \ref{integralwrtnureghorseshoe}, induce the following lower bound on $k\left(\bl_0,\bl\right)$:
\begin{align*}
& \quad k\left(\bl_0,\bl\right)\nonumber \\
& \geq \frac{\left(2\pi\right)^{-\frac{p}{2}}b^{2\left(a+\frac{n}{2}\right)}\alpha^p\left(\sqrt{\omega_*}k^*c^2\right)^{-p}}{\left(\y^T\y+b\right)^{a+\frac{n}{2}}\Gamma\left(a+\frac{n}{2}\right)}\int\limits_{\mathbb{R}_+}\int\limits_{\mathbb{R}_+}\int\limits_{\mathbb{R}^p}\left(\sigma^2\right)^{-\left(a+\frac{n+p}{2}\right)-1}\prod\limits_{j=1}^p\left\{\left(\lambda_j^2\right)^{-2}\right\}\nonumber \\ 
& \qquad \qquad \exp{\left[-\frac{\left(\bb-\Omega^{-1}\X^T\y\right)^T \Omega\left(\bb-\Omega^{-1}\X^T\y\right)+\bb^T\left(\tau^2\Bl\right)^{-1}\bb}{2\sigma^2}\right]}\prod\limits_{j=1}^p\left\{\frac{\left(1+\frac{1}{\lambda_j^2}\right)^{-2}}{1+\frac{\sigma^2\sqrt{\tau^2}}{\beta_j^2}}\right\}\nonumber \\
& \qquad \qquad \exp{\left[-\frac{1}{\sigma^2}\left(\frac{\y^T\y+\y^T\X\left(c^2I_p-\Omega^{-1}\right)\X^T\y}{2}+b\right)\right]}\left(1+\tau^2\right)^{-\frac{p}{2}}\pi_\tau\left(\tau^2\right)d\bb d\sigma^2 d\tau^2 
\end{align*}

\noindent
where $\alpha$ is some positive constant (see Proposition \ref{integralwrtnureghorseshoe}). For the inner most integral wrt $\bb$ we use the lower bound in Proposition \ref{integralwrtbetareghorseshoe} and get the following:
\begin{align*}
& \quad k\left(\bl_0,\bl\right)\nonumber \\
& \geq \frac{b^{2\left(a+\frac{n}{2}\right)}\alpha^p\left(\sqrt{\omega_*}k^*|c|^3\right)^{-p}}{\left(\y^T\y+b\right)^{a+\frac{n}{2}}\Gamma\left(a+\frac{n}{2}\right)}\nonumber \\
& \qquad \qquad \int\limits_{\mathbb{R}_+}\int\limits_{\mathbb{R}_+}\left(\sigma^2\right)^{-\left(a+\frac{n}{2}\right)-1}\exp{\left[-\frac{1}{\sigma^2}\left(\frac{\y^T\y+2\y^T\X\left(c^2I_p-M_{\tau^2}^{-1}\right)\X^T\y}{2}+b\right)\right]}\nonumber \\
& \qquad \qquad \prod\limits_{j=1}^p\left(1+\lambda_j^2\right)^{-2}\times \left|M_{\tau^2}\right|^{-1}\left(1+\frac{\sqrt{\tau^2}}{c^2}\right)^{-p}\left(1+\tau^2\right)^{-\frac{p}{2}}\pi_\tau\left(\tau^2\right)d\sigma^2d\tau^2; 
\end{align*}

\noindent
where $M_{\tau^2}$ is as in Proposition \ref{integralwrtbetareghorseshoe}. It follows that 
\begin{align*}
& \quad k\left(\bl_0,\bl\right)\\
& \geq \frac{b^{2\left(a+\frac{n}{2}\right)}\alpha^p\left(\sqrt{\omega_*}k^*|c|^3\right)^{-p}}{\left(\y^T\y+b\right)^{a+\frac{n}{2}}\Gamma\left(a+\frac{n}{2}\right)}\int\limits_{\mathbb{R}_+}\int\limits_{\mathbb{R}_+}\left(\sigma^2\right)^{-\left(a+\frac{n}{2}\right)-1}\exp{\left[-\frac{1}{\sigma^2}\left(\frac{\y^T\y+2c^2\y^T\X\X^T\y}{2}+b\right)\right]}\nonumber \\
&\qquad \qquad \prod\limits_{j=1}^p\left(1+\lambda_j^2\right)^{-2}\times \left|M_{\tau^2}\right|^{-1}\left(1+\frac{\sqrt{\tau^2}}{c^2}\right)^{-p}\left(1+\tau^2\right)^{-\frac{p}{2}}\pi_\tau\left(\tau^2\right)d\sigma^2d\tau^2, 
\end{align*}

\noindent
since $\y^T\X M_{\tau^2}^{-1}\X^T\y/\sigma^2\geq 0.$ Next by virtue of the 
inverse-gamma integral, we have 

$$\int\limits_{\mathbb{R}_+}\left(\sigma^2\right)^{-\left(a+\frac{n}{2}\right)-1}\exp{\left[-\frac{1}{\sigma^2}\left(\frac{\y^T\y+2c^2\y^T\X\X^T\y}{2}+b\right)\right]}d\sigma^2=\frac{\Gamma\left(a+\frac{n}{2}\right)}{\left(\frac{\y^T\y+2c^2\y^T\X\X^T\y}{2}+b\right)^{a+\frac{n}{2}}}$$

\noindent
This together with the fact that $\mid M_{\tau^2}\mid\leq \left(1+\frac{1}{\tau^2}\right)^p \prod\limits_{j=1}^p\left(\omega_*+\frac{1}{\lambda_j^2}\right)$ gives the following lower bound:
\begin{align*} 
& \quad k\left(\bl_0,\bl\right)\nonumber \\
& \geq \frac{b^{2\left(a+\frac{n}{2}\right)}\alpha^p\left(\sqrt{\omega_*}k^*|c|^3\right)^{-p}}{\left(\y^T\y+b\right)^{a+\frac{n}{2}}\left(\frac{\y^T\y+2c^2\y^T\X\X^T\y}{2}+b\right)^{\left(a+\frac{n}{2}\right)}}\prod\limits_{j=1}^p\left\{\left(1+\lambda_j^2\right)^{-2}\left(\omega_*+\frac{1}{\lambda_j^2}\right)^{-1}\right\}\nonumber \\
& \qquad \qquad \int\limits_{0}^{\infty}\left(1+\frac{1}{\tau^2}\right)^{-p} \left(1+\frac{\sqrt{\tau^2}}{c^2}\right)^{-p}\left(1+\tau^2\right)^{-\frac{p}{2}}\pi_\tau\left(\tau^2\right)d\tau^2 
\end{align*}

\noindent
Further denoting $\eta=\text{max}\left\{1,\omega_*\right\}$ we get
\begin{align*} 
 k\left(\bl_0,\bl\right) = \epsilon^* h\left(\bl\right)
\end{align*}

\noindent
where
$$\epsilon^*= \frac{b^{2\left(a+\frac{n}{2}\right)}\alpha^p\left(2\eta\sqrt{\omega_*}k^*|c|^3\right)^{-p}}{\left(\y^T\y+b\right)^{a+\frac{n}{2}}\left(\frac{\y^T\y+2c^2\y^T\X\X^T\y}{2}+b\right)^{\left(a+\frac{n}{2}\right)}} E_{\pi_\tau\left(\tau^2\right)}\left[\left(1+\frac{1}{\tau^2}\right)^{-p} \left(1+\frac{\sqrt{\tau^2}}{c^2}\right)^{-p}\left(1+\tau^2\right)^{-\frac{p}{2}}\right],$$
and $h$ is a probability density on $\mathbb{R}^p_+$ given by $$h\left(\bl\right)= \prod\limits_{j=1}^p\left\{\frac{2\eta\lambda_j^2}{\left(1+\eta\lambda_j^2\right)^3}I_{(0,\infty)}\left(\lambda_j^2\right)\right\},$$ 

\noindent
and this completes the proof of minorization condition for the MTD $k$ corresponding to the regularized Horseshoe $\bl$-chain. \qed \\
The drift and minorization conditions in Lemma \ref{Theorem1reghorseshoe} and Lemma \ref{Theorem2reghorseshoe} can be combined with Theorem 12 of \Myciteauthor{10.2307/2291067} to establish geometric ergodicity of the regularized Horseshoe Gibbs sampler which is stated as follows:

\begin{thm} \label{thm:regularized:Horseshoe:Piironen:Vehtari}
Suppose the prior density $\pi_\tau (\cdot)$ for the global shrinkage parameter satisfies 
$$\int_{\mathbb{R}_+} u^{\delta/2} \pi_\tau (u) du < \infty$$ for some $\delta \in (0.00162, 0.22176)$. 
Then, the regularized Horseshoe Gibbs sampler with transition kernel $K_{aug,reg}$
is geometrically ergodic.
\end{thm}

\subsection{Geometric ergodicity of a Gibbs sampler for the regularized Horseshoe 
variant in \Myciteauthor{nishimura2019shrinkage}} \label{regularized:Horseshoe:Nishimura:Suchard}
The following variant of regularized Horseshoe shrinkage prior has been introduced in \Myciteauthor{nishimura2019shrinkage}. 
\begin{eqnarray}\label{Nishimurahorseshoepriors}
\pi\left(\beta_j, \lambda_j \mid \tau^2,\sigma^2\right)\propto \frac{1}{\sqrt{\tau^2\lambda_j^2}}\exp{\left[-\frac{\beta_j^2}{2\sigma^2}\left(\frac{1}{c^2}+\frac{1}{\tau^2\lambda_j^2}\right)\right]}\pi_\ell\left(\lambda_j\right) \nonumber\\
& & \text{independently for}\ j=1,2,\cdots,p \nonumber \\
\sigma^2\sim\text{Inverse-Gamma}\left(a,b\right);\ \tau^2\sim\pi_\tau\left(\cdot\right) 
\end{eqnarray}

\noindent
where $\pi_\ell$ and $\pi_\tau$ are probability densities. Note that based on the 
above specification 
$$
\beta_j \mid \lambda_j^2, \tau^2,\sigma^2 \sim N \left( 0, \left( \frac{1}{c^2}+\frac{1}{\tau^2\lambda_j^2} \right) \right)
$$

\noindent
identical to the specification in \cite{piironen2017}. The difference is that instead $\bl, \tau^2$ and $\sigma^2$ having independent priors, we now have 
\begin{equation} \label{nishimura:conditional:prior}
\pi(\lambda_j \mid \tau^2, \sigma^2) = c(\tau^2) \left( 1 + \frac{\tau^2 \lambda_j^2}{c^2} \right)^{-1/2} \pi_\ell (\lambda_j), 
\end{equation}

\noindent
where 
\begin{equation} \label{nishimura:normalizing:constant}
1/c(\tau^2) = \int_0^\infty \left( 1 + \frac{\tau^2 \lambda_j^2}{c^2} \right)^{-1/2} \pi_\ell (\lambda_j) d \lambda_j. 
\end{equation}

\noindent
The principal motivation for the algebraic modification of the prior as compared to 
that of \cite{piironen2017} is the resulting simplification of the posterior 
computation, although an alternative interpretation using fictitious data is also 
discussed in \cite{nishimura2019shrinkage}. In fact, using $\pi_\ell$ to be the half-Cauchy density 
and using its representation in terms of a mixture of Inverse-Gamma densities (\cite{Makalic_2016}) the following conditional 
posterior distributions can be obtained from straightforward computations after augmenting the latent 
variables $\left\{\nu_j\right\}_{j=1}^p$. 
\begin{eqnarray}\label{Nishimurahorseshoeposteriors}
\left.\bb\right|\sigma^2,\tau^2,\bl,\y\sim\mathcal{N}\left(A_c^{-1}\X^T\y,\sigma^2A_c^{-1}\right)\nonumber \\
\left.\sigma^2\right|\tau^2,\bl,\y\sim \text{Inverse-Gamma}\left(a+\frac{n}{2},\frac{\y^T\left(I_n-\X A_c^{-1}\X^T\right)\y}{2}+b\right)\nonumber \\
\left.\nu_j\right|\lambda_j^2,\y\sim\text{Inverse-Gamma}\left(1,1+\frac{1}{\lambda_j^2}\right),\ \text{independently for}\ j=1,2,\cdots,p\nonumber \\
\left.\tau^2\right|\bl,\y\sim\pi\left(\left.\tau^2\right|\bl,\y\right)\nonumber \\
\left.\lambda_j^2\right|\nu_j,\beta_j,\sigma^2,\tau^2,\y\sim \text{Inverse-Gamma}\left(1,\frac{1}{\nu_j}+\frac{\beta_j^2}{2\sigma^2\tau^2}\right),\ \text{independently for}\ j=1,2,\cdots,p\nonumber \\
\end{eqnarray}

\noindent
where 
$$
\pi\left(\left.\tau^2\right|\bl,\y\right)\propto\left|\tau^2 A_c\right|^{-\frac{1}{2}}\left(\frac{\y^T\left(I_n-\X A_c^{-1}\X^T\right)\y}{2}+b\right)^{-\left(a+\frac{n}{2}\right)}\pi_\tau\left(\tau^2\right) c(\tau^2)^p 
$$

\noindent
and $A_c=\X^T\X+\left(\tau^2\Bl\right)^{-1}+c^{-2}I_p$. Most of the above conditional posterior densities, including that for the local shrinkage parameters $\left\{\lambda_j^2\right\}_{j=1}^p$ are standard probability distributions (as opposed to the non-standard ones in the regularized Horseshoe posterior in \eqref{reghorseshoeposteriors}) and can be easily sampled from. An efficient Metropolis sampler for the non-standard (one-dimensional) density  $\pi\left(\left.\tau^2\right|\bl,\y\right)$ can be constructed similar to the one provided in Appendix \ref{Appendix D}. Hence, a two-block Gibbs sampler, whose one step-transition from $({\boldsymbol \beta}_0, \sigma^2_0, {\boldsymbol \nu}_0, \tau^2_0, {\boldsymbol \lambda}_0)$ to 
$({\boldsymbol \beta}, \sigma^2, {\boldsymbol \nu}, \tau^2, {\boldsymbol \lambda})$ is given by sampling sequentially from $\pi({\boldsymbol \beta}, \sigma^2, {\boldsymbol \nu}, \tau^2 \mid {\boldsymbol \lambda}_0, \y)$ and 
$\pi({\boldsymbol \lambda} \mid {\boldsymbol \beta}, \sigma^2, {\boldsymbol \nu}, \tau^2, \y)$, can be used to generate approximate samples from the regularized Horseshoe posterior in \eqref{Nishimurahorseshoeposteriors}. We will denote the Markov transition kernel of this two-block Gibbs sampler by $\tilde K_{aug,reg}$ (analogous to $K_{aug,reg}$ in the regularized Horseshoe setting). The transition 
density can be obtained by substituting the appropriate conditional posterior  densities in the expression (\ref{MTDreghorseshoe}). 

Note that the above conditional posterior distributions are very similar to that for the original Horseshoe Gibbs sampler $K_{aug}$ given by \eqref{posteriors} in Section \ref{original:horseshoe}. The only differences are 
\begin{enumerate}
    \item the matrix $A$ appearing in \eqref{posteriors} has been replaced by $A_c$ (which is the matrix $A$ plus the added regularization introduced in the prior conditional variance of $\bb$ through the constant c) in \eqref{Nishimurahorseshoeposteriors}, and
    \item the form of the posterior conditional density of the global shrinkage parameter, namely, $\pi\left(\left.\tau^2\right|\bl,\y\right)$ is different due to the additional term $c(\tau^2)^p$. 
\end{enumerate}

\begin{thm}\label{TheoremNishimurahorseshoe}
Suppose the prior density of the global shrinkage parameter is truncated below away 
from zero; that is, $\pi_\tau(u)=0$ for $u<T$ for some $T>0$ and satisfies 
$$
\int_T^{\infty}u^{\frac{p+\delta}{2}}\pi_\tau(u) du<\infty
$$

\noindent
for some $\delta\in (0.00162,0.22176)$. Then, the regularized Horseshoe Gibbs sampler corresponding to the transition kernel $\tilde K_{aug,reg}$ is geometrically ergodic. \end{thm}

\noindent
The above theorem can be proved by essentially following verbatim the proof 
of Lemma \ref{Theorem1}  (which establishes geometric ergodicity for 
$K_{aug}$) with the same geometric drift function as in 
Lemma \ref{Theorem1}, and replacing the matrix $A$ by the matrix $A_c$ at 
relevant places. However, appropriate modifications are needed using the following two facts.  \begin{enumerate}
    \item In the original Horseshoe setting, a uniform upper bound for the conditional posterior means of $\beta_j'$s (see  \eqref{conditional_mean_definition} for definition) was established in Proposition \ref{result:lmUnifBound2} in Appendix \ref{Appendix A}. However, in the current context, the added regularization of $c^{-2} I_p$ in $A_c$ immediately provides the uniform upper bound without need for additional analysis. 
    \item The conditional posterior density $\pi\left(\left.\tau^2\right|\bl_0,\y\right)$ is different from the original Horseshoe setting. Hence, the upper bound for the $\delta_0/2^{\text{th}}$ moment of this density for some $\delta_0 \in (0.00162,0.22176)$ (see \eqref{lower_truncation}) needs to be independently established. We have provided this bound in Proposition \ref{Nishimuraupperboundontau^2} of Appendix \ref{Appendix B}. Due to the presence of the additional term $c(\tau^2)^p$ in the conditional density, a stronger assumption of the existence of $(p+\delta_0)/2^{th}$
    moment is required (as compared to the $\delta_0/2^{th}$ moment in 
    Theorem \ref{Theorem 2.1} and Theorem \ref{Theorem1reghorseshoe}). 
\end{enumerate}

\begin{remark} \label{remregular}
In \cite{nishimura2019shrinkage}, the authors focus on Bayesian logistic regression for their geometric ergodicity analysis. They use the regularized Horseshoe prior in (\ref{Nishimurahorseshoepriors}) without the parameter $\sigma^2$ as their is no need for an error variance parameter for the Binomial likelihood. However, for computational purposes, additional parameters ${\boldsymbol \omega} = \left\{\omega_j\right\}_{j=1}^p$ with Polya-Gamma prior distributions are introduced. A two-block Gibbs sampler with blocks $(\bb, \bl)$ and $({\boldsymbol \omega}, \tau^2)$ is then constructed and its geometric ergodicity is then established assuming that the global shrinkage parameter $\tau^2$ is bounded away from zero and infinity 
\cite[Theorem 4.6]{nishimura2019shrinkage}. 

Many details of this analysis break down when translating to the Bayesian linear 
regression framework considered in our paper. The parameters 
${\boldsymbol \omega}$ are now replaced by the error variance parameter $\sigma^2$. 
One can still construct a two-block Gibbs sampler with blocks $(\bb, \bl)$ and 
$(\sigma^2, \tau^2)$, but many conditional independence and other algebraic niceties 
involving ${\boldsymbol \omega}$ which are crucial in establishing the 
minorization condition in the logistic regression context, do not hold analogously 
with $\sigma^2$ in the linear regression context. The structural differences also 
imply that the drift condition with the function $\sum\limits_{j=1}^p |\beta_j|^{-\delta}$ 
does not work out in the linear regression setting. 

\end{remark}

\begin{remark}
The geometric ergodicity result (Theorem \ref{TheoremNishimurahorseshoe}) 
corresponding to the regularized Horseshoe variant in \cite{nishimura2019shrinkage} 
requires truncation of the global shrinkage parameter $\tau^2$ below away from zero. 
Such an assumption is not required for the geometric ergodicity result 
(Theorem \ref{Theorem1reghorseshoe}) corresponding to the regularized Horseshoe of 
\cite{piironen2017}. Also, due to the presence of the additional term 
$(c(\tau^2))^p$ in $\pi\left(\left.\tau^2\right|\bl_0,\y\right)$, a stronger moment assumption is required for Theorem \ref{TheoremNishimurahorseshoe} as compared to Theorem \ref{Theorem1reghorseshoe}. 
\end{remark}

\subsection{A simulation study} \label{simulation:regularized:Horseshoe}

\noindent
The primary objective of this study is to examine the practical 
feasibility/scalability of the two regularized Horseshoe Gibbs samplers 
described in Sections \ref{Regularized:Gibbs:sampler} and 
\ref{regularized:Horseshoe:Nishimura:Suchard}. We consider a simulation 
setting with $n = 500$ samples and $p = 1000$ variables. We generate $10$ 
replicated datasets following exactly the same procedure as outlined in 
Section \ref{simulation:original:Horseshoe}. For each of these 10 datasets, 
we run four Gibbs samplers each: the Gibbs sampler $K_{aug,reg}$ for the 
regularized Horseshoe in \cite{piironen2017} with $c = 1$ and $c = 100$, 
and the Gibbs sampler $\tilde K_{aug,reg}$ for the regularized Horseshoe 
variant in \cite{nishimura2019shrinkage} with $c = 1$ and $c = 100$. 

Again, both algorithms were implemented in {\it R}. Due to maintenance 
issues, an older machine albeit with the same OS/RAM/processor 
specifications was used for these experiments as compared to the one used in
Section \ref{simulation:original:Horseshoe}. The run-times for $2500$ 
iterations of each Gibbs sampler for each replication and each value of $c$ 
are provided in Table \ref{table_RegularizedHS}. Cumulative average plots 
for the function $\bb^T \bb$ were used to monitor and confirm sufficient 
mixing of all the Markov chains. In all the settings, and across all the 
replications, the Gibbs samplers roughly needed 5500 seconds to complete the 
required 2500 iterations. 

We also tried to use the Hamiltonian Monte Carlo based algorithm for the  
regularized Horseshoe in \cite{piironen2017}, as implemented in the {\it R} 
package {\it hsstan}. However, the maximum treedepth (set to $10$) is exceeded in 
{\it all} of the 2500 iterations. This issue persists even after warming up for up 
to 7000 iterations, and then running for 2500 more iterations. As we understand, this
indicates poor adaptation, and raises questions about adequate posterior exploration 
and  mixing of the Markov chain. A proposed remedy in this setting (Chapter 15.2 of 
the Stan reference manual on {\it mc-stan.org}) is to increase the tree depth. The 
{\it hsstan} function, however, did not allow us to pass the {\it max\_treedepth} or 
{\it max\_depth} as a parameter and change its value. Anyway, from the point of view 
of scalability, the time taken per iteration with maximum treedepth $10$ was roughly 
one-fourth as compared to the various Gibbs samplers. When the maximum treedepth is 
increased appropriately to resolve the issue pointed out above, it is very likely 
that the time taken per iteration will be around the same or more than those of the 
Gibbs samplers (increasing the tree-depth by $1$ in the No U-turn HMC sampler 
effectively doubles the computation time). 

To conclude, the Gibbs samplers described in 
Sections \ref{Regularized:Gibbs:sampler} and 
\ref{regularized:Horseshoe:Nishimura:Suchard} provide practically feasible 
approaches which are computationally competitive with the HMC based approach. The 
geometric ergodicity results in 
Theorems \ref{Theorem1reghorseshoe} and \ref{TheoremNishimurahorseshoe} help provide the practitioner with asymptotically valid standard error estimates for corresponding MCMC based approximations to posterior quantities of interest. 
\begin{table}[hbt!]
    \begin{subtable}[h]{0.45\textwidth}
        \centering
        \begin{tabular}{| l | l | l |}
        \hline
        Replication \# & $c=100$ & $c=1$ \\
        \hline \hline
        1 & 5505.93 & 5783.77\\
        2 & 5704.37 & 5755.31\\
        3 & 5563.22 & 5756.45\\
        4 & 5644.05 & 5807.36\\
        5 & 5534.64 & 5917.48\\
        6 & 5491.00 & 5692.67\\
        7 & 5487.85 & 5629.29\\
        8 & 5589.35 & 5812.13\\
        9 & 5634.18 & 5819.75\\
       10 & 5679.75 & 5715.75\\
        \hline
       \end{tabular}
       \caption{Gibbs sampler $K_{aug, reg}$ for the regularized Horseshoe posterior in \cite{piironen2017}}
       \label{tab:week1}
    \end{subtable}
    \hfill
    \begin{subtable}[h]{0.45\textwidth}
        \centering
        \begin{tabular}{| l | l | l |}
        \hline
        Replication \# & $c=100$ & $c=1$ \\
        \hline \hline
        1 & 5422.53 & 5289.71\\
        2 & 5488.5  & 5461.39\\
        3 & 5459.94 & 5425.37\\
        4 & 5479.54 & 5272.11\\
        5 & 5486.41 & 5373.16\\
        6 & 5332.89 & 5430.52\\
        7 & 5549.28 & 5180.04\\
        8 & 5343.74 & 5102.16\\
        9 & 5462.57 & 5446.23\\
       10 & 5252.53 & 5103.76\\
        \hline
        \end{tabular}
        \caption{Gibbs sampler $\tilde K_{aug,reg}$ for the regularized Horseshoe variant in \cite{nishimura2019shrinkage}}
        \label{tab:week2}
     \end{subtable}
     \caption{The run-times (in seconds) required to generate 2500 MCMC samples for Gibbs samplers corresponding to the regularized Horseshoe in \cite{piironen2017} and the regularized Horseshoe variant in \cite{nishimura2019shrinkage} for 10 simulated data sets with 
     $n=500$ and $p=1000$.}
     \label{table_RegularizedHS}
\end{table}

\appendix

\section{Uniform bound on $\mu_j$} \label{Appendix A}

\noindent
The goal of this subsection is to show that $\mu_j = 
\boldsymbol{e}_j^TA_0^{-1}\X^T\y$ defined in 
(\ref{conditional_mean_definition}) is uniformly bounded in $\bl_0$ (even 
when $n < p$). This result will be established through a sequence of five 
propositions. 
\begin{prop} \label{result:lm1}
Let $\lmLambda\in \R^{p\times p}$ be any diagonal matrix with positive diagonal elements $\lmlambda_1, \lmlambda_2, \ldots, \lmlambda_p $.  Let $\lmX\in \R^{n\times p}$ be any matrix with rank $r$.  Let  the singular value decomposition of $X$ is $\lmX=\lmU\lmD\lmV^T$ where  $ \lmD\in \R^{r\times r}$  is diagonal matrix with positive diagonal elements $\lmd_1, \ldots \lmd_r$ while    $\lmU\in \R^{n\times r}$ and  $\lmV\in \R^{p\times r}$ are such that  $\lmU^T\lmU=I$ and  $\lmV^T\lmV=I$.   If $\lmlambda\leq \min\{\lmlambda_1, \ldots , \lmlambda_p\}$ be any positive number then for arbitrary $\lmy\in \R^n$
$$   \lmy^T\lmX (\lmX^T\lmX +\lmLambda)^{-1}\lmX^T\lmy\leq \lmy^T\lmy-  \norm{P_{\lmU^\perp}\lmy}^2-\sum_{i=1}^{r}\frac{\lmlambda \tilde{\lmu}^2_i}{\lmd_i^2+\lmlambda},$$
where $\tilde{\lmu}_i$ is the $i^{\text{th}}$ component of the vector $\tilde{\lmu}=\lmU^T\lmy$ and $ P_{\lmU^\perp}  $ is the orthogonal projection matrix for the orthogonal complement of the column space of $\lmU$.
\end{prop}

\proof
Without loss of generality we assume that the matrix $\lmLambda$ is diagonal matrix with diagonal elements $\lmlambda_1, \lmlambda_2, \ldots, \lmlambda_p $ where $0<\lmlambda_1\leq \lmlambda_2\leq  \ldots \leq  \lmlambda_p. $ According to the condition of the  result $\lmlambda\leq \lmlambda_1$. Now we define a set of diagonal matrices $\{\lmLambda_{(1)}, \ldots, \lmLambda_{(j)} , \ldots, \lmLambda_{(p)} \}$ in the following manner. For $j=1, \ldots (p-1)$, the matrix $\lmLambda_{(j)} $ has first $j$ diagonal elements to be identical and equal to $\lmlambda$ while rest of the $(p-j)$ diagonal elements are identical as that of the matrix $\lmLambda$. Also let  $\lmLambda_{(p)}=\lmlambda I_{p\times p} $ and  $\lmLambda_{(0)}=\lmLambda$. The above set of matrices satisfy the following relation
$$ \lmLambda_{(j-1)}=\lmLambda_{(j)}+ (\lmlambda_j-\lmlambda) \lme_j\lme_j^T  \text{ for } j=1 , \ldots, p $$
where $\lme_j\in\R^p$ denotes the $j^{\text{th}}$ elementary vector.  Now using the Sherman-Woodbury formula for inverting matrices we get that 
\begin{eqnarray}
& & \left( \lmX^T\lmX+\lmLambda_{(j-1)} \right)^{-1} \nonumber\\
& = &  \left( \lmX^T\lmX+\lmLambda_{(j)}+  (\lmlambda_j-\lmlambda) \lme_j\lme_j^T  \right)^{-1} \nonumber\\
& = &  \left( \lmX^T\lmX+\lmLambda_{(j)}\right)^{-1}- \frac{(\lmlambda_j-\lmlambda)   \left( \lmX^T\lmX+\lmLambda_{(j)}\right)^{-1} \lme_j\lme_j^T \left( \lmX^T\lmX+\lmLambda_{(j)}\right)^{-1}   }{ 1- (\lmlambda_j-\lmlambda) \lme_j^T \left( \lmX^T\lmX+\lmLambda_{(j)}\right)^{-1} \lme_j    } .\nonumber
\end{eqnarray}
Consequently, 
\begin{eqnarray}\label{eq:lmSWInversion}
\left( \lmX^T\lmX+\lmLambda_{(j-1)} \right)^{-1}-  \left( \lmX^T\lmX+\lmLambda_{(j)}\right)^{-1} = - \frac{(\lmlambda_j-\lmlambda)   \left( \lmX^T\lmX+\lmLambda_{(j)}\right)^{-1} \lme_j\lme_j^T \left( \lmX^T\lmX+\lmLambda_{(j)}\right)^{-1}   }{ 1- (\lmlambda_j-\lmlambda) \lme_j^T \left( \lmX^T\lmX+\lmLambda_{(j)}\right)^{-1} \lme_j    },
\end{eqnarray}
Aggregating the equations \ref{eq:lmSWInversion}  over $j=1, \ldots p$,   we get that  

\begin{eqnarray}\label{eq:lmMatrixRecursion}
\left( \lmX^T\lmX+\lmLambda \right)^{-1}= \left( \lmX^T\lmX+\lmlambda I\right)^{-1}  - \sum_{j=1}^{p}\frac{(\lmlambda_j-\lmlambda)   \left( \lmX^T\lmX+\lmLambda_{(j)}\right)^{-1} \lme_j\lme_j^T \left( \lmX^T\lmX+\lmLambda_{(j)}\right)^{-1}   }{ 1- (\lmlambda_j-\lmlambda) \lme_j^T \left( \lmX^T\lmX+\lmLambda_{(j)}\right)^{-1} \lme_j    },
\end{eqnarray}
where we have used the fact that  $\lmLambda_{(0)}=\lmLambda$ and $\lmLambda_{(p)}=\lmlambda I $. If $\lmy\in \R^n$ be arbitrary vector then  it follows from \ref{eq:lmMatrixRecursion} that 

\begin{eqnarray}\label{eq:lmIneqDiagIdentity}
& & \lmy^T \lmX \left( \lmX^T\lmX+\lmLambda \right)^{-1} \lmX^T \lmy \nonumber\\
& = & \lmy^T \lmX \left( \lmX^T\lmX+\lmlambda I\right)^{-1} \lmX^T \lmy  - \sum_{j=1}^{p}\frac{(\lmlambda_j-\lmlambda)   \lmy^T \lmX\left( \lmX^T\lmX+\lmLambda_{(j)}\right)^{-1} \lme_j\lme_j^T \left( \lmX^T\lmX+\lmLambda_{(j)}\right)^{-1} \lmX^T \lmy  }{ 1- (\lmlambda_j-\lmlambda) \lme_j^T \left( \lmX^T\lmX+\lmLambda_{(j)}\right)^{-1} \lme_j    }\nonumber\\
& = & \lmy^T \lmX \left( \lmX^T\lmX+\lmlambda I\right)^{-1} \lmX^T \lmy  - \sum_{j=1}^{p}\frac{(\lmlambda_j-\lmlambda)   \norm{\lme_j^T \left( \lmX^T\lmX+\lmLambda_{(j)}\right)^{-1} \lmX^T \lmy }^2 }{ 1- (\lmlambda_j-\lmlambda) \lme_j^T \left( \lmX^T\lmX+\lmLambda_{(j)}\right)^{-1} \lme_j    } \nonumber\\
& \leq & \lmy^T \lmX \left( \lmX^T\lmX+\lmlambda I\right)^{-1} \lmX^T \lmy, 
\end{eqnarray}
because $\lmlambda_j\geq \lmlambda$.  
Now consider the singular value decomposition of the matrix $\lmX=\lmU \lmD\lmV^T$ where $\lmD$ is a diagonal matrix with positive diagonal elements $d_1, \ldots , d_r$, $\lmU\in \R^{n\times r}, \lmV\in R^{p\times r}$ such that  $\lmU^T\lmU =\lmV^T\lmV=I$. Also let $\lmV^{\perp} \in \R^{p\times (p-r)}$ be  such that the matrix  $\left[ \lmV , \lmV^{\perp}\right] \in \R^{p\times p}$ is a orthogonal matrix, i.e. the columns of $\lmV^{\perp}$ constitutes a orthonormal basis for the orthogonal complement of the column space of $\lmV$.  Now consider the fact that 
\begin{eqnarray}
(\lmX^T\lmX+\lmlambda I )^{-1}
 & = & \left[ \lmV\lmD\lmU^T \lmU\lmD \lmV^T +\lmlambda I \right]^{-1} \nonumber\\
 & = & \left[ \lmV\lmD^2 \lmV^T +\lmlambda \lmV\lmV^T +\lmlambda \lmV^{\perp} (\lmV^{\perp})^T \right]^{-1} \nonumber\\
 & = & \left[ \lmV(\lmD^2+ \lmlambda I ) \lmV^T +\lmlambda \lmV^{\perp} (\lmV^{\perp})^T \right]^{-1} \nonumber\\
 & = & \lmV(\lmD^2+ \lmlambda I ) ^{-1}\lmV^T +\frac{1}{\lmlambda} \lmV^{\perp} (\lmV^{\perp})^T .\nonumber
\end{eqnarray}

Note that  $(\lmV^{\perp})^T \lmV=0$. Thus
\begin{eqnarray}\label{eq:lmSum1}
\lmy^T \lmX \left( \lmX^T\lmX+\lmlambda I\right)^{-1} \lmX^T \lmy 
& =&  \lmy^T \lmX \left[ \lmV(\lmD^2+ \lmlambda I ) ^{-1}\lmV^T +\frac{1}{\lmlambda} \lmV^{\perp} (\lmV^{\perp})^T\right] \lmV\lmD\lmU^T \lmy \nonumber\\
& =&  \lmy^T (\lmU\lmD\lmV^T)  \lmV(\lmD^2+ \lmlambda I ) ^{-1}\lmV^T  \lmV\lmD\lmU^T \lmy \nonumber\\
& =&  \lmy^T\lmU\lmD(\lmD^2+ \lmlambda I )^{-1}\lmD\lmU^T \lmy \nonumber\\
& =& \sum_{i=1}^{r} \frac{d_i^2 \tilde{\lmu}^2_i}{d_i^2+\lmlambda},
\end{eqnarray}
where $d_1, \ldots d_r>0$ are the diagonal elements of the matrix $\lmD$ and $ \tilde{\lmu}_i$ is the $i^{\text{th}}$ entry of the vector $\lmU^T\lmy$.  Let  $\lmU^{\perp}$ refers to the orthogonal completion of the matrix $\lmU$ then 
\begin{eqnarray}\label{eq:lmNorm1}
\sum_{i=1}^{r}  \tilde{\lmu}_i^2= \lmy^T\lmU\lmU^T\lmy=\lmy^T\lmy -\lmy^T \lmU^{\perp} (\lmU^{\perp})^T\lmy = \lmy^T\lmy -\norm{ P_{\lmU^{\perp}}\lmy}^2,
\end{eqnarray}
where $P_{\lmU^{\perp}}$ denotes the orthogonal projection for the column space of ${\lmU^{\perp}}$.  Finally, it follows from \ref{eq:lmIneqDiagIdentity},  \ref{eq:lmSum1} and \ref{eq:lmNorm1}  that 
\begin{eqnarray}
\lmy^T \lmX \left( \lmX^T\lmX+\lmLambda \right)^{-1} \lmX^T \lmy \leq  \sum_{i=1}^{r}  \tilde{\lmu}_i^2 - \sum_{i=1}^{r} \frac{\lmlambda\tilde{\lmu}^2_i}{d_i^2+\lmlambda} = \lmy^T\lmy -\norm{ P_{\lmU^{\perp}}\lmy}^2 - \sum_{i=1}^{r} \frac{\lmlambda\tilde{\lmu}^2_i}{d_i^2+\lmlambda} .\nonumber
\end{eqnarray}
Note that $\norm{ P_{\lmU^{\perp}}\lmy}^2 + \sum_{i=1}^{r} \frac{\lmlambda\tilde{\lmu}^2_i}{d_i^2+\lmlambda} >0$ beacuse $ \sum_{i=1}^{r}  \tilde{\lmu}_i^2+\norm{ P_{\lmU^{\perp}}\lmy}^2=\lmy^T\lmy>0$.
\qed

\vspace{.25in}

%
%

\begin{prop} \label{result:lm2} 
Let  $\lmX=\left[ \lmxbf_1, \ldots , \lmxbf_p   \right]\in \R^{n\times p}$ and $\lmX_2=\left[ \lmxbf_2, \ldots , \lmxbf_p   \right]$.
Assume  $\lmX_{2}=UDV^T$ be the singular value decomposition  where $d_1, \ldots d_r>0$ are the diagonal elements of $\lmD$.  Let $\lmdelta_1>0$ and $\lmDelta_{p-1}$ be any diagonal matrix with positive diagonal elements $\lmdelta_2, \ldots, \lmdelta_p$.  If  
$T:= (\norm{\lmxbf_1}^2+\lmdelta_1)- \lmxbf_1^T \lmX_2 (\lmX_{2}^T\lmX_2+\lmDelta_{p-1})^{-1}\lmX_2\lmxbf_1$ then  for any $0<\lmdelta\leq \min\{ \lmdelta_1, \ldots, \lmdelta_p \},$
\begin{eqnarray}
T\geq \delta_1+\norm{P_{U^\perp}\lmxbf_1}^2+\sum_{i=1}^{r} \frac{\lmdelta\tilde{\lmu}^2_i}{d_i^2+\lmdelta},\nonumber
\end{eqnarray}
%
where $\tilde{\lmu}_i$ is the $i^{\text{th}}$ component of the vector $\tilde{\lmu}=\lmU^T\lmxbf_{1}$ and $ P_{\lmU^\perp}  $ is the orthogonal projection matrix for the orthogonal complement of the column space of $\lmU$.
\end{prop}

\proof
 Using Result~\ref{result:lm1},  for any $\lmdelta\leq min\{ \lmdelta_1, \ldots, \lmdelta_p \}$, we get that  
\begin{eqnarray}
 \lmxbf_1^T \lmX_2 (\lmX_{2}^T\lmX_2+\lmDelta_{p-1})^{-1}\lmX_2\lmxbf_1
& \leq & \lmxbf_1^T\lmxbf_1-\norm{P_{U^\perp}\lmxbf_1}^2-\sum_{i=1}^{r} \frac{\lmdelta\tilde{\lmu}^2_i}{d_i^2+\lmdelta}\nonumber
\end{eqnarray}
Consequently $T= (\lmxbf_1^T\lmxbf_1+\lmdelta_1)- \lmxbf_1^T \lmX_2 (\lmX_{2}^T\lmX_2+\lmDelta_{p-1})^{-1}\lmX_2\lmxbf_1\geq \lmdelta_1+\norm{P_{U^\perp}\lmxbf_1}^2+\sum_{i=1}^{r} \frac{\lmdelta\tilde{\lmu}^2_i}{d_i^2+\lmdelta}.$
\qed.
 \vspace{.25in}

\begin{prop} \label{result:lm2.1}
Let $\lmX\in \R^{n\times p}$ be arbitrary matrix and $\lmDelta_p\in \R^{p\times p}$ be any diagonal matrix with positive diagonal elements $\lmdelta_1, \ldots, \lmdelta_p$.  Consider the following partition of the matrix  
$$\lmX^T\lmX+\lmDelta_p=\left[
\begin{array}{c|c}
\norm{\lmxbf_{1}}^2+\lmdelta_1 &     \lmxbf_1^T\lmX_2 \\ \hline
  \lmX_2^T\lmxbf_1 &\lmX_2^T\lmX_2+\lmDelta_{p-1}
\end{array}\right],$$
 where $\lmDelta_{p-1}$ is the diagonal matrix with diagonal elements $\lmdelta_2, \ldots, \lmdelta_p$. If  $(\lmX_2^T\lmX_2+\lmDelta_{p-1})^{-1}\lmDelta_{p-1}$ is uniformly bounded, then the first column of the matrix 
$$  \left( \lmX^T\lmX+\lmDelta_p \right)^{-1}\lmDelta_p $$ is uniformly bounded. The notations $\lmxbf_1$ and $ \lmX_2$ are as they are defined  in the Result~\ref{result:lm2}.
\end{prop}  
 
 \proof
%
If we  consider the partition of 
$$
\lmX^T\lmX+\lmDelta_p=\left[
\begin{array}{c|c}
\norm{\lmxbf_{1}}^2+\lmdelta_1 &     \lmxbf_1^T\lmX_2 \\ \hline
  \lmX_2^T\lmxbf_1 &\lmX_2^T\lmX_2+\lmDelta_{p-1}
\end{array}\right],
$$
then the Schur complement of the first block of the matrix is given as 
$$T= (\norm{\lmxbf_{1}}^2+\lmdelta_1)- \lmxbf_1^T\lmX_2 (\lmX_2^T\lmX_2+\lmDelta_{p-1})^{-1}\lmX_2^T\lmxbf_1.$$
 Employing the inversion formula of the block matrices  \citep{lu2002inverses}, we get that 
 \begin{eqnarray}\label{eq:lmStep1}
 & & (\lmX^T\lmX+\lmDelta_p)^{-1} \nonumber\\
 & = &  \left[\begin{array}{c|c}
\frac{1}{T}      &  \frac{\lmxbf_1^T\lmX_2(\lmX_2^T\lmX_2+\lmDelta_{p-1})^{-1}}{T}      \\ \hline
 \frac{(\lmX_2^T\lmX_2+\lmDelta_{p-1})^{-1} \lmX_2^T\lmxbf_1}{T}    &   (\lmX_2^T\lmX_2+\lmDelta_{p-1})^{-1}+\frac{(\lmX_2^T\lmX_2+\lmDelta_{p-1})^{-1}   \lmX_2^T\lmxbf_1\lmxbf_1^T\lmX_2  (\lmX_2^T\lmX_2+\lmDelta_{p-1})^{-1}}{T} \nonumber
\end{array}\right].
\end{eqnarray}  
 In the next two bullet points, we are going to show if  $(\lmX_2^T\lmX_2+\lmDelta_{p-1})^{-1}\lmDelta_{p-1}$ is uniformly bounded then so is all the entries of the vector $$ \left[\begin{array}{cc}
\frac{\lmdelta_1}{T}  & \frac{ \lmdelta_1 (\lmX_2^T\lmX_2+\lmDelta_{p-1})^{-1} \lmX_2^T\lmxbf_1}{T}            
\end{array}\right]^T,  $$
which is the first column of the matrix $(\lmX^T\lmX+\lmDelta_p)^{-1}\lmDelta_p $.

\begin{itemize}
\item { \bf To show $0<\frac{\lmdelta_1}{T}\leq 1:$  }\\
It is evident from the Result~\ref{result:lm2} that  $T>0$ .   Therefore $ \frac{\lmdelta_1}{T}>0$ as $\delta_1>0$ as well.  On the other hand, a direct implication of Result~\ref{result:lm2} is that for $0<\lmdelta\leq \min\{\lmdelta_1, \ldots , \lmdelta_p\}$, 
\begin{eqnarray}\label{eq:lm1byT}
\frac{\lmdelta_1}{T}
& \leq  &\frac{\lmdelta_1}{ \lmdelta_1+\norm{P_{U^\perp}\lmxbf_1}^2+\sum_{i=1}^{r} \frac{\lmdelta\tilde{\lmu}^2_i}{d_i^2+\lmdelta}    }\leq  1.
\end{eqnarray}
where the details about the notations  $\tilde{\lmu}_i$, $ P_{\lmU^\perp} , d_i $ can be found in Result~\ref{result:lm2}. \\

\item {\bf To show $ \frac{ \lmdelta_1 (\lmX_2^T\lmX_2+\lmDelta_{p-1})^{-1} \lmX_2^T\lmxbf_1}{T}  $ uniformly bounded:       }\\
Let $\lmxv_1, \lmxv_2$ be such that $\lmxbf_1=\lmxv_1+ \lmxv_2$ where $\lmxv_1$ belongs to the column space of $\lmX_2$ and $\lmxv_2$ belongs to the orthogonal complement of the column space of $\lmX_2$. Therefore $\lmxv_1=\lmX_2 l $ for some vector  $l\in \R^{p-1}$. Consequently, 
\begin{eqnarray}\label{eq:lmUNIFpartRow1}
  (\lmX_2^T\lmX_2+\lmDelta_{p-1})^{-1} \lmX_2^T\lmxbf_1
& =& (\lmX_2^T\lmX_2+\lmDelta_{p-1})^{-1}\lmX_2^T( \lmX_2 l +\lmxv_2 ) \nonumber\\
& =&(\lmX_2^T\lmX_2+\lmDelta_{p-1})^{-1}\lmX_2^T\lmX_2 l \nonumber\\
& =&  \left [ I-  (\lmX_2^T\lmX_2+\lmDelta_{p-1})^{-1}\lmDelta_{p-1}\right] l , \nonumber
\end{eqnarray}
is uniformly bounded as  we are assuming that the matrix  $(\lmX_2^T\lmX_2+\lmDelta_{p-1})^{-1}\lmDelta_{p-1}$ is uniformly bounded. Combining this fact along with \ref{eq:lm1byT}, we conclude that all the entries of the  vector  $\frac{ \lmdelta_1 (\lmX_2^T\lmX_2+\lmDelta_{p-1})^{-1} \lmX_2^T\lmxbf_1}{T}  $  are also uniformly bounded. 
\end{itemize}
\qed

\begin{prop}\label{result:lmUnifBound1}
Let $\lmX\in \R^{n\times p}$ be arbitrary matrix and $\lmDelta_p\in \R^{p\times p}$ be any diagonal matrix with positive diagonal elements $\lmdelta_1, \ldots, \lmdelta_p$. Then for arbitrary $p$ and $n$ the matrix $ \left( \lmX^T\lmX+\lmDelta_p \right)^{-1}\lmDelta_p  $ is uniformly bounded. Specifically 
$$ \sup_{\lmdelta_1, \ldots, \lmdelta_p>0}\left\vert e_{i}^T  \left( \lmX^T\lmX+\lmDelta_p \right)^{-1}\lmDelta_p e_j \right\vert<C  $$
where $C$ is a finite constant that does not depend on $\lmdelta_1, \ldots, \lmdelta_p$. 
\end{prop}

 \proof
 We will show the result by induction on the integer $k$ where the hypothesis of induction is as follows,
  
$  \mathcal{H}(k): $ Let $n$ be arbitrary positive integer. Then for any positive integer $k$,  the matrix   $\left( \lmX^T\lmX+\lmDelta_k \right)^{-1}\lmDelta_k$ is uniformly bounded for all $X\in \R^{n\times k}$  and arbitrary diagonal matrix with positive diagonal elements $\lmDelta_k\in \R^{k\times k}$.  \\

 {\bf Initial step: }  The hypothesis  trivially  holds for $k=1$. We will show that $\mathcal{H}(k)$ is true for the case  $k=2.$ Let $\lmX=\left[\lmxbf_1, \lmxbf_2  \right] \in \R^{n\times 2}$ and  $ \lmDelta_2=\begin{bmatrix}
 \lmdelta_{1} & 0\\
0 & \lmdelta_{2}
\end{bmatrix}$, $\lmdelta_{1}, \lmdelta_{2}>0$ be arbitrary.  Define 
 $\lmA:=\begin{bmatrix}
 \lma_{1,1} & \lma_{1,2}\\
\lma_{2,1} & \lma_{2,2}
\end{bmatrix} :=\lmX^T\lmX \text{ and }   \text{ then }$

\begin{eqnarray}
\left( \lmA+\lmDelta_2 \right)^{-1}\lmDelta_2
=\frac{1}{(\lma_{1,1}\lma_{2,2}-\lma_{1,2}\lma_{2,1})+ \lmdelta_1\lma_{2,2}+\lmdelta_2\lma_{1,1}+ \lmdelta_1\lmdelta_2}
  \begin{bmatrix}
 \lmdelta_1(\lma_{2,2}+\lmdelta_2) & -\lmdelta_2\lma_{2,1}\\
-\lmdelta_1\lma_{1,2} & \lmdelta_2(\lma_{1,1}+\lmdelta_1)
\end{bmatrix}  \nonumber
\end{eqnarray}

\begin{itemize}
\item Note that
\begin{eqnarray}
 \sup_{\lmdelta_1, \lmdelta_2>0}\left\vert e_{1}^T  \left( \lmA+\lmDelta_2 \right)^{-1}\lmDelta_2 e_1 \right\vert
 & = &  \sup_{\lmdelta_1, \lmdelta_2>0}\frac{    \lmdelta_1(\lma_{2,2}+\lmdelta_2)   }{(\lma_{1,1}\lma_{2,2}-\lma_{1,2}\lma_{2,1})+ \lmdelta_1(\lma_{2,2} + \lmdelta_2) +\lmdelta_2\lma_{1,1}}
\leq 1 \nonumber
\end{eqnarray}
because $(\lma_{1,1}\lma_{2,2}-\lma_{1,2}\lma_{2,1})\geq 0$ as $\lmX^T\lmX$ is nonnegative definte matrix. Additionally  $\lma_{1,1}=\norm{\lmxbf_1}^2\geq 0$ and $\lma_{2,2}=\norm{\lmxbf_2}^2\geq 0$, where $\lmX=[\lmxbf_1, \lmxbf_2].$

 \item  
\begin{eqnarray}\label{eq:lmUpper}
 & & \sup_{\lmdelta_1, \lmdelta_2>0}\left\vert e_{1}^T  \left( \lmA+\lmDelta_2 \right)^{-1}\lmDelta_2 e_2 \right\vert\nonumber\\
 & = &  \sup_{\lmdelta_1, \lmdelta_2>0}\frac{  \vert \lma_{2,1}\lmdelta_2  \vert }{(\lma_{1,1}\lma_{2,2}-\lma_{1,2}\lma_{2,1})+ \lmdelta_1(\lma_{2,2} + \lmdelta_2) +\lmdelta_2\lma_{1,1}}\\
& \leq &  \frac{\vert\lma_{2,1}\vert}{\lma_{1,1}}\nonumber
\end{eqnarray} 
 for the case when $\lma_{1,1}\neq 0$. On the contrary, if   $\lma_{1,1}=0$ then it follows from \ref{eq:lmUpper} that 
 \begin{eqnarray}
\sup_{\lmdelta_1, \lmdelta_2>0}\left\vert e_{1}^T  \left( \lmA+\lmDelta_2 \right)^{-1}\lmDelta_2 e_2 \right\vert
& = &0,\end{eqnarray}
because $\lma_{1,1}=\norm{\lmxbf_1}^2=0$ implies that $\lma_{2,1}=\lmxbf_2^T\lmxbf_1=0$. 
\end{itemize}
 In a similar fashion we can show that  the absolute value of the other two entries of the matrix $ \left( \lmX^T\lmX+\lmDelta_2 \right)^{-1}\lmDelta_2$ can be bounded above by numbers that does not depend on $\lmdelta_1, \lmdelta_2$. Consequently $\mathcal{H}(k)$ holds for $k=2.$

\vspace{.2in}
{\bf Induction step:}
Let $\mathcal{H}(k)$ holds for  $k=1, 2, \ldots , (p-1)$.  We will show that the result holds for $k=p$ as well. Let $\lmX\in \R^{n\times p}$ be arbitrary matrix and  $\lmDelta_{p}$ for diagonal matrix with positive diagonal elements $\lmdelta_1, \ldots , \lmdelta_p$.  Consider the partition of the matrices $\lmX^T\lmX+\lmDelta_p$ as follows 
\begin{eqnarray}\label{eq:lmPartition}
\lmX^T\lmX+\lmDelta_p=\left[
\begin{array}{c|c}
\norm{\lmxbf_{1}}^2+\lmdelta_1 &     \lmxbf_1^T\lmX_2 \\ \hline
  \lmX_2^T\lmxbf_1 &\lmX_2^T\lmX_2+\lmDelta_{p-1}
\end{array}\right],
\end{eqnarray}
where $\lmX_{2}, \lmxbf_1, \lmDelta_{p-1}$ are as it is in Result~\ref{result:lm2}.  As   it  satisfies the conditions of the induction hypothesis $\mathcal{H}{(p-1)}$, the matrix $(\lmX_{2}^T\lmX_{2}+\lmDelta_{p-1})^{-1}\lmDelta_{p-1}$  is uniformly bounded.  Therefore using Result~\ref{result:lm2.1}, the first column of the matrix $(\lmX^T\lmX+\lmDelta_p)^{-1}\lmDelta_p $ is uniformly bounded.

  In remaining of the proof, we show that  the $m^{\text{th}}$ column of  $(\lmX^T\lmX+\lmDelta_p)^{-1}\lmDelta_p$ is uniformly bounded for any  $m>1$.  Consider the permutation matrix $P_{1, m}=\left[ e_m,e_2 \ldots ,e_{m-1},e_{1} \ldots, e_p \right]$.  Note that $P_{1,m}$ can be generated by exchanging the $1^{\text{st}}$ and $m^{\text{th}}$ columns of an identity matrix.   $P_{1, m}$ is a symmetric and orthogonal matrix, i.e. $P_{1,m}^T=P_{1,m}$ and $P_{1,m}^TP_{1,m}=P_{1,m}^2=I$. Now consider 
\begin{eqnarray}\label{eq:lmPerm}
P_{1,m} (\lmX^T\lmX+\lmDelta_p)^{-1}\lmDelta_p P_{1,m}  
& = &  P_{1,m}(\lmX^T\lmX+\lmDelta_p)^{-1} P_{1,m} P_{1,m}\lmDelta_p P_{1,m}\nonumber\\
& = &  (P_{1,m}^T \lmX^T\lmX P_{1,m} +P_{1,m} \lmDelta_p P_{1,m} )^{-1} P_{1,m}\lmDelta_p P_{1,m}\nonumber\\
& = &  ( {\lmX^{*}}^T \lmX^{*}  + \lmDelta_p^{*} )^{-1} \lmDelta_p^{*},
\end{eqnarray}
where the $ \lmX^{*}:=\lmX P_{1,m} $ is obtained by exchanging the first and $m^{\text{th}}$ columns of $\lmX$  while $ \lmDelta_p^{*}: =P_{1,m}\lmDelta_p P_{1,m}$ is the diagonal matrix where the first and the $m^{\text{th}}$ diagonal elements of $\lmDelta_p$ are exchanged.  We can represent ${\lmX^*}^T\lmX^*+\lmDelta^*_p$ as
$$\left[
\begin{array}{c|c}
\norm{\lmxbf^*_{1}}^2+\lmdelta^*_1 &     {\lmxbf^*_1}^T\lmX^{*}_2 \\ \hline
  {\lmX_2^*}^T\lmxbf^*_1 &{\lmX^*_2}^T\lmX_2+\lmDelta^*_{p-1}
\end{array}\right], $$
\noindent
where the notations are equivalent to that of the ones in \ref{eq:lmPartition}. The matrix $({\lmX^*_{2}}^T\lmX^*_{2}+\lmDelta^*_{p-1})^{-1}\lmDelta^*_{p-1}$  satisfies the conditions of the induction hypothesis $\mathcal{H}{(p-1)}$, thus  it is uniformly bounded.  Therefore using Result~\ref{result:lm2.1}, the first column of the matrix $({\lmX^*}^T\lmX^*+\lmDelta^*_p)^{-1}\lmDelta^*_p $ is uniformly bounded as well.
 It follows from \ref{eq:lmPerm} that the permuted version of the  first column of  $( {\lmX^{*}}^T \lmX^{*}  + \lmDelta_p^{*} )^{-1} \lmDelta_p^{*} $ is
\begin{eqnarray}
P_{1,m}\left[( {\lmX^{*}}^T \lmX^{*}  + \lmDelta_p^{*} )^{-1} \lmDelta_p^{*} e_1 \right]
& =&  P_{1,m}P_{1,m} (\lmX^T\lmX+\lmDelta_p)^{-1}\lmDelta_p P_{1,m}   e_1 = \left[  (\lmX^T\lmX+\lmDelta_p)^{-1}\lmDelta_p \right]  e_m, \nonumber
\end{eqnarray}
 which is the $m^{\text{th}}$ column of $ (\lmX^T\lmX+\lmDelta_p)^{-1}\lmDelta_p$. 
Therefore, we infer that  all the columns of the matrix  $(\lmX^T\lmX+\lmDelta_p)^{-1}\lmDelta_p$ are uniformly bounded and conclude that  $\mathcal{H}(k)$ holds for the case $k=p$.
 
 
 \qed
\vspace{.2in }

\begin{prop}\label{result:lmUnifBound2}
Let $\lmX\in \R^{n\times p}$ be arbitrary matrix and $\lmDelta_p\in \R^{p\times p}$ be any diagonal matrix with positive diagonal elements $\lmdelta_1, \ldots, \lmdelta_p$. Then for arbitrary $p$ and $n$ 
\begin{enumerate}
\item The matrix $ \left( \lmX^T\lmX+\lmDelta_p \right)^{-1} \lmX^T\lmX  $ is uniformly bounded. Specifically 
$$ \sup_{\lmdelta_1, \ldots, \lmdelta_p>0}\left\vert e_{i}^T  \left( \lmX^T\lmX+\lmDelta_p \right)^{-1} \lmX^T\lmX e_j \right\vert<C  $$
where $C$ is a finite constant that does not depend on $\lmdelta_1, \ldots, \lmdelta_p$. 
\item The vector $ \left( \lmX^T\lmX+\lmDelta_p \right)^{-1} \lmX^T\lmy$ is uniformly bounded.
\end{enumerate}
\end{prop}

\proof

part(1):
Note that $ \left(\lmX^T\lmX+\lmDelta_p \right)^{-1} \lmX^T\lmX= I-\left(\lmX^T\lmX+\lmDelta_p \right)^{-1}  \lmDelta_p$. Using Result\ref{result:lmUnifBound1} we know that the matrix $\left(\lmX^T\lmX+\lmDelta_p \right)^{-1}  \lmDelta_p$ is uniformly bounded. Consequently $\left(\lmX^T\lmX+\lmDelta_p \right)^{-1} \lmX^T\lmX$ is also uniformly bounded. \\

part(2):

Let $\lmy=\lmxv_1+ \lmxv_2$ where $\lmxv_1$ belongs to the column space of $\lmX$ and $\lmxv_2$ belongs to the perpendicular to the column space of $\lmX$. Therefore $\lmxv_1=\lmX l $ for some vector  $l\in \R^{p-1}$. Consequently, 

\begin{eqnarray}\label{eq:lmUNIFpartRow}
(\lmX^T\lmX+\lmDelta_{p})^{-1}\lmX^T\lmy
& = &  (\lmX^T\lmX+\lmDelta_{p})^{-1}\lmX^T( \lmxv_1+\lmxv_2  ) \nonumber\\
& =&   (\lmX^T\lmX+\lmDelta_{p})^{-1}\lmX^T( \lmX l +\lmxv_2 ) \nonumber\\
& =&   \left[(\lmX^T\lmX+\lmDelta_{p})^{-1}\lmX^T\lmX \right]l .\nonumber
\end{eqnarray}
Therefore part(a) of the result ensures that the $(\lmX^T\lmX+\lmDelta_{p})^{-1}\lmX^T\lmy$ is uniformly bounded. 

\qed

\section{Other technical results} \label{Appendix B}

\begin{prop}\label{upperboundontau^2}
Let $\delta$ be chosen as in Lemma \ref{Theorem1}. Then for any 
$\epsilon>0$ there exists $C_1>0$ (not depending on $\bl_0$) such that  $$
\E\left[\left.\left(\tau^2\right)^{\frac{\delta}{2}}\right|\bl_0,\y\right]
\leq C_1 
$$.
\end{prop}

\proof For any $\epsilon > 0$, note that 
\begin{eqnarray}\label{upperboundtau^2:1}
 \E\left[\left.\left(\tau^2\right)^{\frac{\delta}{2}}\right| \bl_0,\y\right]&=&\E\left[\left.\left(\tau^2\right)^{\frac{\delta}{2}}I_{\left[\tau^2<\epsilon\right]}\right| \bl_0,\y\right]+\E\left[\left.\left(\tau^2\right)^{\frac{\delta}{2}}I_{\left[\tau^2\geq\epsilon\right]}\right| \bl_0,\y\right]\nonumber \\
&\leq & \epsilon^{\frac{\delta}{2}}+\E\left[\left.\left(\tau^2\right)^{\frac{\delta}{2}}I_{\left[\tau^2\geq\epsilon\right]}\right| \bl_0,\y\right] 
\end{eqnarray}

\noindent
Next we demonstrate an upper bound to the second term in \eqref{upperboundtau^2:1}.
\begin{eqnarray}
 \E\left[\left.\left(\tau^2\right)^{\frac{\delta}{2}}I_{\left[\tau^2\geq\epsilon\right]}\right| \bl_0,\y\right]&=&\int\limits_{\epsilon}^\infty\left(\tau^2\right)^{\frac{\delta}{2}}\pi\left(\left.\tau^2\right|\bl_0,\y\right)d\tau^2\nonumber \\
 &=& \frac{\int\limits_{\epsilon}^\infty\left(\tau^2\right)^{\frac{\delta}{2}}\left(\frac{\y^T\left(I_n-\X A_0^{-1}\X^T\right)\y}{2}+b\right)^{-\left(a+\frac{n}{2}\right)}\frac{\pi_\tau\left(\tau^2\right)}{\left|I_p+\tau^2\X^T\X\cdot\Bl_0\right|^{\frac{1}{2}}}d\tau^2}{\int\limits_{0}^\infty\left(\frac{\y^T\left(I_n-\X A_0^{-1}\X^T\right)\y}{2}+b\right)^{-\left(a+\frac{n}{2}\right)}\frac{\pi_\tau\left(\tau^2\right)}{\left|I_p+\tau^2\X^T\X\cdot\Bl_0\right|^{\frac{1}{2}}}d\tau^2}\nonumber \\
 &\leq& \left(1+\frac{\y^T\y}{b}\right)^{a+\frac{n}{2}}\frac{\int\limits_{\epsilon}^\infty\left(\tau^2\right)^{\frac{\delta}{2}}\frac{\pi_\tau\left(\tau^2\right)}{\left|I_p+\tau^2\X^T\X\cdot\Bl_0\right|^{\frac{1}{2}}}d\tau^2}{\int\limits_{0}^\epsilon\frac{\pi_\tau\left(\tau^2\right)}{\left|I_p+\tau^2\X^T\X\cdot\Bl_0\right|^{\frac{1}{2}}}d\tau^2}\nonumber \\
 &\leq& \left(1+\frac{\y^T\y}{b}\right)^{a+\frac{n}{2}}\frac{\int\limits_{\epsilon}^\infty\left(\tau^2\right)^{\frac{\delta}{2}}\frac{\pi_\tau\left(\tau^2\right)}{\left|I_p+\epsilon\X^T\X\cdot\Bl_0\right|^{\frac{1}{2}}}d\tau^2}{\int\limits_0^\epsilon\frac{\pi_\tau\left(\tau^2\right)}{\left|I_p+\epsilon\X^T\X\cdot\Bl_0\right|^{\frac{1}{2}}}d\tau^2}\nonumber \\
 &\leq&\frac{\left(1+\frac{\y^T\y}{b}\right)^{a+\frac{n}{2}}}{\int\limits_0^\epsilon\pi_\tau\left(\tau^2\right)d\tau^2}\int\limits_\epsilon^\infty\left(\tau^2\right)^{\frac{\delta}{2}}\pi_\tau\left(\tau^2\right)d\tau^2\nonumber \\
 &<& \infty. 
\end{eqnarray}

\noindent
This completes the proof with $C_1=\epsilon^{\frac{\delta}{2}}+\frac{\left(1+\frac{\y^T\y}{b}\right)^{a+\frac{n}{2}}}{\int\limits_{0}^\epsilon\pi_\tau\left(\tau^2\right)d\tau^2}\int\limits_{\epsilon}^\infty\left(\tau^2\right)^{\frac{\delta}{2}}\pi_\tau\left(\tau^2\right)d\tau^2$. \qed 

\begin{prop}\label{upperboundontau^2:negative:moment}
Suppose there exists a $\delta > 0.00162$ such that 
$$
\int_0^\infty (\tau^2)^{-\frac{p+\delta}{2}} \pi_\tau (\tau^2) d \tau^2 < 
\infty. 
$$

\noindent
Then for any $\epsilon>0$ there exists $\tilde{C}_1>0$ (not depending on $\bl_0$) 
such that  
$$
\E\left[\left.\left(\tau^2\right)^{\frac{-\delta}{2}}\right|\bl_0,\y\right] \leq \tilde{C}_1. 
$$
\end{prop}

\proof For any $\epsilon > 0$, note that 
\begin{eqnarray}\label{upperboundtau^2:negative:moment}
 \E\left[\left.\left(\tau^2\right)^{-\frac{\delta}{2}}\right| \bl_0,\y\right]&=&\E\left[\left.\left(\tau^2\right)^{-\frac{\delta}{2}}I_{\left[\tau^2<\epsilon\right]}\right| \bl_0,\y\right]+\E\left[\left.\left(\tau^2\right)^{\frac{\delta}{2}}I_{\left[\tau^2\geq\epsilon\right]}\right| \bl_0,\y\right]\nonumber \\
&\leq & \epsilon^{-\frac{\delta}{2}}+\E\left[\left.\left(\tau^2\right)^{-\frac{\delta}{2}}I_{\left[\tau^2\leq\epsilon\right]}\right| \bl_0,\y\right] \end{eqnarray}

\noindent
Next we demonstrate an upper bound to the second term in \eqref{upperboundtau^2:negative:moment}.
\begin{eqnarray}
 \E\left[\left.\left(\tau^2\right)^{-\frac{\delta}{2}}I_{\left[\tau^2\leq\epsilon\right]}\right| \bl_0,\y\right]&=&\int\limits_0^{\epsilon}\left(\tau^2\right)^{-\frac{\delta}{2}}\pi\left(\left.\tau^2\right|\bl_0,\y\right)d\tau^2 \nonumber\\
 &=& \frac{\int\limits_0^{\epsilon}\left(\tau^2\right)^{-\frac{\delta}{2}}\left(\frac{\y^T\left(I_n-\X A_0^{-1}\X^T\right)\y}{2}+b\right)^{-\left(a+\frac{n}{2}\right)}\frac{\pi_\tau\left(\tau^2\right)}{\left|I_p+\tau^2\X^T\X\cdot\Bl_0\right|^{\frac{1}{2}}}d\tau^2}{\int\limits_{0}^\infty\left(\frac{\y^T\left(I_n-\X A_0^{-1}\X^T\right)\y}{2}+b\right)^{-\left(a+\frac{n}{2}\right)}\frac{\pi_\tau\left(\tau^2\right)}{\left|I_p+\tau^2\X^T\X\cdot\Bl_0\right|^{\frac{1}{2}}}d\tau^2} \nonumber\\
 &\leq& \left(1+\frac{\y^T\y}{b}\right)^{a+\frac{n}{2}}\frac{\int\limits_0^{\epsilon}\left(\tau^2\right)^{-\frac{\delta}{2}}\frac{\pi_\tau\left(\tau^2\right)}{\left|I_p+\tau^2\X^T\X\cdot\Bl_0\right|^{\frac{1}{2}}}d\tau^2}{\int\limits_\epsilon^\infty\frac{\pi_\tau\left(\tau^2\right)}{\left|I_p+\tau^2\X^T\X\cdot\Bl_0\right|^{\frac{1}{2}}}d\tau^2} \nonumber\\
 &=& \left(1+\frac{\y^T\y}{b}\right)^{a+\frac{n}{2}}\frac{\int\limits_0^{\epsilon}\left(\tau^2\right)^{-\frac{p+\delta}{2}}\frac{\pi_\tau\left(\tau^2\right)}{\left|\tau^{-2}I_p+\X^T\X\cdot\Bl_0\right|^{\frac{1}{2}}}d\tau^2}{\int\limits_\epsilon^\infty\left(\tau^2\right)^{-\frac{p}{2}}\frac{\pi_\tau\left(\tau^2\right)}{\left|\tau^{-2}I_p+\X^T\X\cdot\Bl_0\right|^{\frac{1}{2}}}d\tau^2} \nonumber\\
 &\leq& \left(1+\frac{\y^T\y}{b}\right)^{a+\frac{n}{2}}\frac{\int\limits_0^{\epsilon}\left(\tau^2\right)^{-\frac{p+\delta}{2}}\frac{\pi_\tau\left(\tau^2\right)}{\left|\epsilon^{-1}I_p+\X^T\X\cdot\Bl_0\right|^{\frac{1}{2}}}d\tau^2}{\int\limits_\epsilon^\infty\left(\tau^2\right)^{-\frac{p}{2}} \frac{\pi_\tau\left(\tau^2\right)}{\left|\epsilon^{-1}I_p+\X^T\X\cdot\Bl_0\right|^{\frac{1}{2}}}d\tau^2} \nonumber\\
 &\leq& \left(1+\frac{\y^T\y}{b}\right)^{a+\frac{n}{2}}\frac{\int\limits_0^{\epsilon}\left(\tau^2\right)^{-\frac{p+\delta}{2}}\pi_\tau\left(\tau^2\right)}{\int\limits_\epsilon^\infty\left(\tau^2\right)^{-\frac{p}{2}} \pi_\tau\left(\tau^2\right)d\tau^2} \nonumber\\
&<& \infty. 
\end{eqnarray}

\noindent
\qed

\begin{prop}\label{reghsupperboundontau^2}
Refer to \ref{reghorseshoeposteriors} and Lemma \ref{Theorem1reghorseshoe}. Then for any $\epsilon>0$ and for any $\delta>0$, there exists $C_2>0$ such that  $$\E\left[\left.\left(\tau^2\right)^{\frac{\delta}{2}}\right|\bl_0,\y\right]\leq C_2$$.
\end{prop}

\proof Fix an $\epsilon>0$ and a $\delta>0$. \\
\begin{eqnarray}\label{reghsupperboundtau^2:1}
 \E\left[\left.\left(\tau^2\right)^{\frac{\delta}{2}}\right| \bl_0,\y\right]&=&\E\left[\left.\left(\tau^2\right)^{\frac{\delta}{2}}I_{\left[\tau^2<\epsilon\right]}\right| \bl_0,\y\right]+\E\left[\left.\left(\tau^2\right)^{\frac{\delta}{2}}I_{\left[\tau^2\geq\epsilon\right]}\right| \bl_0,\y\right]\nonumber \\
&\leq & \epsilon^{\frac{\delta}{2}}+\E\left[\left.\left(\tau^2\right)^{\frac{\delta}{2}}I_{\left[\tau^2\geq\epsilon\right]}\right| \bl_0,\y\right] \nonumber \\
\end{eqnarray}

\noindent
Next we demonstrate an upper bound to the second term in \eqref{reghsupperboundtau^2:1}.
\begin{eqnarray}
 & & \E\left[\left.\left(\tau^2\right)^{\frac{\delta}{2}}I_{\left[\tau^2\geq\epsilon\right]}\right| \bl_0,\y\right] \nonumber\\
 &=&\int\limits_{\epsilon}^\infty\left(\tau^2\right)^{\frac{\delta}{2}}\pi\left(\left.\tau^2\right|\bl_0,\y\right)d\tau^2\nonumber \\
 &=& \frac{\int\limits_{\epsilon}^\infty\left(\tau^2\right)^{\frac{\delta}{2}}\left(\frac{\y^T\left(I_n-\X A_0^{-1}\X^T\right)\y}{2}+b\right)^{-\left(a+\frac{n}{2}\right)}\frac{\left|c^{-2}I_p+\left(\tau^2\Bl_0\right)^{-1}\right|^{1/2}}{\left|\X^T\X+c^{-2}I_p+\left(\tau^2\Bl_0\right)^{-1}\right|^{1/2}}\pi_\tau\left(\tau^2\right)d\tau^2}{\int\limits_{0}^\infty\left(\frac{\y^T\left(I_n-\X A_0^{-1}\X^T\right)\y}{2}+b\right)^{-\left(a+\frac{n}{2}\right)}\frac{\left|c^{-2}I_p+\left(\tau^2\Bl_0\right)^{-1}\right|^{1/2}}{\left|\X^T\X+c^{-2}I_p+\left(\tau^2\Bl_0\right)^{-1}\right|^{1/2}}\pi_\tau\left(\tau^2\right)d\tau^2}\nonumber \\
 &\leq& \left(1+\frac{\y^T\y}{b}\right)^{a+\frac{n}{2}}\frac{\int\limits_{\epsilon}^\infty\left(\tau^2\right)^{\frac{\delta}{2}}\frac{\left|c^{-2}I_p+\left(\tau^2\Bl_0\right)^{-1}\right|^{1/2}}{\left|\X^T\X+c^{-2}I_p+\left(\tau^2\Bl_0\right)^{-1}\right|^{1/2}}\pi_\tau\left(\tau^2\right)d\tau^2}{\int\limits_0^\epsilon\frac{\left|c^{-2}I_p+\left(\tau^2\Bl_0\right)^{-1}\right|^{1/2}}{\left|\X^T\X+c^{-2}I_p+\left(\tau^2\Bl_0\right)^{-1}\right|^{1/2}}\pi_\tau\left(\tau^2\right)d\tau^2} \label{det_ratio}\\
 &\leq&\frac{\left(1+\frac{\y^T\y}{b}\right)^{a+\frac{n}{2}}}{\int\limits_{0}^\epsilon\pi_\tau\left(\tau^2\right)d\tau^2}\int\limits_{0}^\infty\left(\tau^2\right)^{\frac{\delta}{2}}\pi_\tau\left(\tau^2\right)d\tau^2\nonumber \\
 &<& \infty \nonumber,
\end{eqnarray}

\noindent
Note that the ratio of two determinants inside the integral in the numerator and denominator in \eqref{det_ratio} can be represented in the form 
$$
\frac{|B_1 + \tau^{-2} I_p|}{|B_1 + B_2 + \tau^{-2} I_p|} = \frac{\prod_{k=1}^p 
(s_k (B_1) + \tau^{-2})}{\prod_{k=1}^p (s_k (B_1 + B_2) + \tau^{-2})} 
$$

\noindent
for appropriate symmetric non-negative definite matrices $B_1$ and $B_2$, and their 
respective eigenvalues denoted by $s_k (\cdot)$. Since every eigenvalue of $B_1$ is 
bounded above by the corresponding eigenvalue of $B_1 + B_2$, it follows that the 
ratio of determinants is a decreasing function of $\tau^2$, and can be replaced by 
the value at $\tau^2 = \epsilon$ in both places with the inequality going in the 
right direction. This completes the proof with 
$$
C_2=\epsilon^{\frac{\delta}{2}}+\frac{\left(1+\frac{\y^T\y}{b}\right)^{a+\frac{n}{2}}}{\int\limits_{0}^\epsilon\pi_\tau\left(\tau^2\right)d\tau^2}\int\limits_{0}^\infty\left(\tau^2\right)^{\frac{\delta}{2}}\pi_\tau\left(\tau^2\right)d\tau^2
$$. 

\noindent
\qed

\begin{prop}\label{Nishimuraupperboundontau^2}
Let $\delta$ be chosen as in Theorem \ref{TheoremNishimurahorseshoe}. Then 
for any $\epsilon>0$ there exists $C_3>0$ (not depending on $\bl_0$) such 
that  
$$
\E\left[\left.\left(\tau^2\right)^{\frac{\delta}{2}}\right|\bl_0,\y\right]
\leq C_3. 
$$
\end{prop}

\proof For any $\epsilon > 0$ note that 
\begin{eqnarray}\label{Nishimuraupperboundtau^2:1}
 \E\left[\left.\left(\tau^2\right)^{\frac{\delta}{2}}\right| \bl_0,\y\right]&=&\E\left[\left.\left(\tau^2\right)^{\frac{\delta}{2}}I_{\left[\tau^2<\epsilon\right]}\right| \bl_0,\y\right]+\E\left[\left.\left(\tau^2\right)^{\frac{\delta}{2}}I_{\left[\tau^2\geq\epsilon\right]}\right| \bl_0,\y\right]\nonumber \\
&\leq & \epsilon^{\frac{\delta}{2}}+\E\left[\left.\left(\tau^2\right)^{\frac{\delta}{2}}I_{\left[\tau^2\geq\epsilon\right]}\right| \bl_0,\y\right]. 
\end{eqnarray}

\noindent
Next we demonstrate an upper bound to the second term in \eqref{Nishimuraupperboundtau^2:1}.
\begin{eqnarray}
 \E\left[\left.\left(\tau^2\right)^{\frac{\delta}{2}}I_{\left[\tau^2\geq\epsilon\right]}\right| \bl_0,\y\right]&=&\int\limits_{\epsilon}^\infty\left(\tau^2\right)^{\frac{\delta}{2}}\pi\left(\left.\tau^2\right|\bl_0,\y\right)d\tau^2\nonumber \\
 &=& \frac{\int\limits_{\epsilon}^\infty\left(\tau^2\right)^{\frac{\delta}{2}}\left(\frac{\y^T\left(I_n-\X A_0^{-1}\X^T\right)\y}{2}+b\right)^{-\left(a+\frac{n}{2}\right)}\frac{\pi_\tau\left(\tau^2\right) c(\tau^2)^p}{\left|\tau^2\left(\X^T\X+c^{-2}I_p\right)+\Bl_0^{-1}\right|^{\frac{1}{2}}}d\tau^2}{\int\limits_{T}^\infty\left(\frac{\y^T\left(I_n-\X A_0^{-1}\X^T\right)\y}{2}+b\right)^{-\left(a+\frac{n}{2}\right)}\frac{\pi_\tau\left(\tau^2\right) c(\tau^2)^p}{\left|\tau^2\left(\X^T\X+c^{-2}I_p\right)+\Bl_0^{-1}\right|^{\frac{1}{2}}}d\tau^2}\nonumber \\
 &\leq& \left(1+\frac{\y^T\y}{b}\right)^{a+\frac{n}{2}}\frac{\int\limits_{\epsilon}^\infty\left(\tau^2\right)^{\frac{\delta}{2}}\frac{\pi_\tau\left(\tau^2\right) c(\tau^2)^p}{\left|\tau^2\left(\X^T\X+c^{-2}I_p\right)+\Bl_0^{-1}\right|^{\frac{1}{2}}}d\tau^2}{\int\limits_{T}^\epsilon\frac{\pi_\tau\left(\tau^2\right)c(\tau^2)^p}{\left|\tau^2\left(\X^T\X+c^{-2}I_p\right)+\Bl_0^{-1}\right|^{\frac{1}{2}}}d\tau^2}\nonumber \\
 &\leq& \left(1+\frac{\y^T\y}{b}\right)^{a+\frac{n}{2}}\frac{\int\limits_{\epsilon}^\infty\left(\tau^2\right)^{\frac{\delta}{2}}\frac{\pi_\tau\left(\tau^2\right) c(\tau^2)^p}{\left|\epsilon\left(\X^T\X+c^{-2}I_p\right)+\Bl_0^{-1}\right|^{\frac{1}{2}}}d\tau^2}{\int\limits_{T}^\epsilon\frac{\pi_\tau\left(\tau^2\right)c(\tau^2)^p}{\left|\epsilon\left(\X^T\X+c^{-2}I_p\right)+\Bl_0^{-1}\right|^{\frac{1}{2}}}d\tau^2}\nonumber \\
 &\leq& \tilde{C}^{\frac{p}{2}}\frac{\left(1+\frac{\y^T\y}{b}\right)^{a+\frac{n}{2}}}{\int\limits_{T}^\epsilon\pi_\tau\left(\tau^2\right)d\tau^2}\int\limits_{T}^\infty\left(\tau^2\right)^{\frac{\delta}{2}}\left(1+\tau^2\right)^{\frac{p}{2}}\pi_\tau\left(\tau^2\right)d\tau^2.  
\end{eqnarray}

\noindent
The last inequality follows from the fact that 
$$
1 \leq c(\tau^2) \leq \tilde{C} \sqrt{1 + \tau^2} 
$$

\noindent
for an appropriate constant $\tilde{C}$. \qed

\begin{prop}\label{cauchyschwarzinequality}
$f:\mathbb{R}^p\mapsto [0,\infty)$ and $g:\mathbb{R}^p\mapsto (0,\infty)$ be two functions such that $\int\limits_{\mathbb{R}^p} f(\boldsymbol{x})d\boldsymbol{x}<\infty$ and  $0<\int\limits_{\mathbb{R}^p} f(\boldsymbol{x})g(\boldsymbol{x})d\boldsymbol{x}<\infty$. Then $\int\limits_{\mathbb{R}^p}\frac{f(\boldsymbol{x})}{g(\boldsymbol{x})}d\boldsymbol{x}\geq \frac{\left(\int\limits_{\mathbb{R}^p} f(\boldsymbol{x})d\boldsymbol{x}\right)^2}{\int\limits_{\mathbb{R}^p} f(\boldsymbol{x})g(\boldsymbol{x})d\boldsymbol{x}}$.
\end{prop}

\proof Follows from Cauchy-Schwarz inequality. \qed

\begin{prop}\label{integralwrtnureghorseshoe}
For any $j\in\left\{1,2,\cdots,p\right\}$ and any $d>0$, there exists some $\alpha>0$ such that
$$\displaystyle\int\limits_0^\infty\ \frac{\nu_j^{-2}\exp{\left[-\frac{1}{\nu_j}\left(1+d^{\frac{2}{\delta}}+\frac{1}{\lambda_j^2}\right)\right]}}{\sqrt{\nu_j}+\frac{\sigma^2\sqrt{\tau^2}}{\beta_j^2}}d\nu_j\geq \alpha\frac{\left(1+\frac{1}{\lambda_j^2}\right)^{-2}}{\left(1+\frac{\sigma^2\sqrt{\tau^2}}{\beta_j^2}\right)}.$$ 
\end{prop}

\proof Follows from Proposition \ref{cauchyschwarzinequality} with \\

$f\left(\nu_j\right)=\nu_j^{-2}\exp{\left[-\frac{1}{\nu_j}\left(1+d^{\frac{2}{\delta}}+\frac{1}{\lambda_j^2}\right)\right]}$, \\
$g\left(\nu_j\right)=\sqrt{\nu_j}+\frac{\sigma^2\sqrt{\tau^2}}{\beta_j^2}$, \\
and $\alpha=\left(1+d^{2/\delta}\right)^{-2}$ \qed

\begin{prop} \label{integralwrtbetareghorseshoe} There exists a positive definite matrix $M_{\tau^2}$ such that \\
\begin{align*}
\qquad \displaystyle\int\limits_{\mathbb{R}^p}\ \frac{\exp{\left[-\frac{\left(\bb-\Omega^{-1}\X^T\y\right)^T \Omega\left(\bb-\Omega^{-1}\X^T\y\right)+\bb^T\left(\tau^2\Bl\right)^{-1}\bb}{2\sigma^2}\right]}}{\prod\limits_{j=1}^p\left(1+\frac{\sigma^2\sqrt{\tau^2}}{\beta_j^2}\right)}d\bb &\geq&  \left(2\pi\sigma^2\right)^{\frac{p}{2}}|c|^{-p}\mid M_{\tau^2}\mid^{-1}\left(1+\frac{\sqrt{\tau^2}}{c^2}\right)^{-p} \\
 \qquad \qquad &\times&\exp{\left[-\frac{\y^T\X\left(c^2I_p+\Omega^{-1}-2M_{\tau^2}^{-1}\right)\X^T\y}{2\sigma^2}\right]}
\end{align*}
\end{prop}

\proof Follows from Proposition \ref{cauchyschwarzinequality} with \\

$f\left(\bb\right)=\exp{\left[-\frac{\left(\bb-\Omega^{-1}\X^T\y\right)^T \Omega\left(\bb-\Omega^{-1}\X^T\y\right)+\bb^T\left(\tau^2\Bl\right)^{-1}\bb}{2\sigma^2}\right]}$, \\ $g\left(\bb\right)=\prod\limits_{j=1}^p\left(1+\frac{\sigma^2\sqrt{\tau^2}}{\beta_j^2}\right)$, \\ and $M_{\tau^2}=\Omega+\left(\tau^2\Bl\right)^{-1}$ \qed

\section{Minorization condition for Horseshoe Gibbs sampler} \label{Appendix C}

\begin{lem}\label{Theorem2}
For every $d > 0$, there exists a constant $\epsilon^*=\epsilon^*\left(V,d\right)>0$ and a density function $h$ on $\mathbb{R}_+^p$ such that
\begin{equation}\label{minorizationcondition}
   k\left(\bl_0,\bl\right)\geq \epsilon^* h\left(\bl\right) 
\end{equation}
for every $\bl_0\in B\left(V,d\right)$ (see Section \ref{originalhorseshoeminorizationset} for definition).
\end{lem}

\noindent
{\it Proof}: Fix a $\bl_0\in B\left(V,d\right)$. In order to prove \eqref{minorizationcondition} we will demonstrate appropriate lower bounds to the conditional densities appearing in \eqref{MTD}. From \eqref{posteriors} we have the following:

\begin{eqnarray*}
  \pi\left(\left.\tau^2\right|\bl_0,\y\right) &\geq&  \left(\frac{b}{\y^T\y/2+b}\right)^{a+\frac{n}{2}}\omega_*^{-p/2} \left(1+\tau^2\right)^{-p/2}\pi_\tau\left(\tau^2\right)
\end{eqnarray*}

\noindent

where $\omega_*=\text{max}\left\{1,\bar \omega\cdot d^{2/\delta_0}\right\}$ (recall that $\bar \omega$ is the maximum eigenvalue of $\X^T\X$ and that the prior density $\pi_\tau$ is truncated below at some $T>0$).
\begin{eqnarray*}
  \pi\left(\left.\boldsymbol{\nu}\right|\bl_0,\y\right) &\geq&  \prod_{j=1}^p\left\{\nu_j^{-2}\exp{\left[-\frac{1}{\nu_j}\left(1+d^{\frac{2}{\delta_1}}\right)\right]}\right\}\nonumber \\
  \pi\left(\left.\sigma^2\right|\tau^2,\bl_0,\y\right) &\geq& \frac{b^{a+\frac{n}{2}}}{\Gamma\left(a+\frac{n}{2}\right)}\left(\sigma^2\right)^{-\left(a+\frac{n}{2}\right)-1}\exp{\left[-\frac{1}{\sigma^2}\left(\frac{\y^T\y}{2}+b\right)\right]}\nonumber \\
  \pi\left(\left.\boldsymbol{\beta}\right|\sigma^2,\tau^2,\bl_0,\y\right) &\geq& \left(2\pi\sigma^2\right)^{-\frac{p}{2}}d^{-p/\delta_0}\left(\tau^2\right)^{-p/2}\nonumber \\
\qquad \qquad &\times&  \exp{\left[-\frac{\left(\boldsymbol{\beta}-M^{-1}\X^T\y\right)^T M\left(\boldsymbol{\beta}-M^{-1}\X^T\y\right)+\y^T\left(I-\X M^{-1}\X\right)\y}{2\sigma^2}\right]}\nonumber \\
\end{eqnarray*}

\noindent
since, \begin{eqnarray*}
    & & \left(\boldsymbol{\beta}-A_0^{-1}\X^T\y\right)^T A_0\left(\boldsymbol{\beta}-A_0^{-1}\X^T\y\right) \nonumber \\
    &=& \boldsymbol{\beta}^T A_0\boldsymbol{\beta}-2\boldsymbol{\beta}^T\X^T\y+\y^T\X A_0^{-1}\X^T\y\nonumber \\
    &\leq& \boldsymbol{\beta}^T M\boldsymbol{\beta}-2\boldsymbol{\beta}^T\X^T\y+\y^T\y\nonumber \\
    &=&\left(\boldsymbol{\beta}-M^{-1}\X^T\y\right)^T M\left(\boldsymbol{\beta}-M^{-1}\X^T\y\right)+\y^T\left(I-\X M^{-1}\X^T\right)\y\nonumber \\
\end{eqnarray*}

\noindent
where $M=\omega^*\left(1+\frac{1}{\tau^2}\right)I_p$ and  $\omega^*=\text{max}\left\{\bar{\omega},d^{2/\delta_1}\right\}$.
Finally, 
\begin{eqnarray}\label{lowerbounds}
    \pi\left(\left.\bl\right|\boldsymbol{\beta,\nu},\sigma^2,\tau^2,\y\right) &\geq&  \prod_{j=1}^p\left\{\frac{\beta_j^2}{2\sigma^2\tau^2}\left(\lambda_j^2\right)^{-2}\exp{\left[-\frac{1}{\lambda_j^2}\left(\frac{1}{\nu_j}+\frac{\beta_j^2}{2\sigma^2\tau^2}\right)\right]}\right\} 
\end{eqnarray}
  
\noindent

Putting all lower bounds in \eqref{lowerbounds} in the equation of MTD \eqref{MTD} we have:

\begin{eqnarray*}\label{minorizationproof}
k\left(\bl_0,\bl\right) 
&\geq& \left(2\pi\right)^{-\frac{p}{2}}\left(\omega_*\right)^{-p/2}d^{-p/\delta_0}\left(\frac{b}{\y^T\y/2+b}\right)^{a+\frac{n}{2}}\frac{b^{a+\frac{n}{2}}}{\Gamma\left(a+\frac{n}{2}\right)} \nonumber \\ 
&\times& \int_{[T,\infty)}\int_{\mathbb{R}_+}\int_{\mathbb{R}^p}\int_{\mathbb{R}_+^p}\prod_{j=1}^p\left\{\nu_j^{-2}\exp{\left[-\frac{1}{\nu_j}\left(1+d^{\frac{2}{\delta_1}}+\frac{1}{\lambda_j^2}\right)\right]}\right\}\nonumber \\
&\times& \exp{\left[-\frac{\left(\boldsymbol{\beta}-M^{-1}\X^T\y\right)^T M\left(\boldsymbol{\beta}-M^{-1}\X^T\y\right)+\boldsymbol{\beta}^T\left(\tau^2\Bl\right)^{-1}\boldsymbol{\beta}}{2\sigma^2}\right]}\nonumber \\
&\times&\exp{\left[-\frac{\y^T\left(I-\X M^{-1}\X^T\right)\y+\y^T\y+2b}{2\sigma^2}\right]} \nonumber \\
&\times& \left(\sigma^2\right)^{-\left(a+\frac{n+p}{2}\right)-1}\prod_{j=1}^p\left\{\frac{\beta_j^2}{2\sigma^2\tau^2}\left(\lambda_j^2\right)^{-2}\right\}\left(1+\tau^2\right)^{-p/2}\left(\tau^2\right)^{-p/2}\pi_\tau\left(\tau^2\right)d\boldsymbol{\nu}d\boldsymbol{\beta}d\sigma^2d\tau^2\nonumber \\
\end{eqnarray*}

\noindent
Next we perform the inner integral wrt $\boldsymbol \nu$ and noting that $1+d^{\frac{2}{\delta_1}}+\frac{1}{\lambda_j^2}\leq \left(1+d^{\frac{2}{\delta_1}}\right)\left(1+\frac{1}{\lambda_j^2}\right)$ we have:

\begin{eqnarray*}
k\left(\bl_0,\bl\right) &\geq& \left(2\pi\right)^{-\frac{p}{2}}\left(\omega_*\right)^{-p/2}d^{-p/\delta_0}\left(1+d^{\frac{2}{\delta_1}}\right)^{-p}\left(\frac{b}{\y^T\y/2+b}\right)^{a+\frac{n}{2}}\frac{b^{a+\frac{n}{2}}}{\Gamma\left(a+\frac{n}{2}\right)} \nonumber \\
&\times& \int_{[T,\infty)}\int_{\mathbb{R}_+}\int_{\mathbb{R}^p}\prod_{j=1}^p\left\{\frac{\beta_j^2}{2\sigma^2}\right\}\exp{\left[-\frac{\left(\boldsymbol{\beta}-M^{-1}\X^T\y\right)^T M\left(\boldsymbol{\beta}-M^{-1}\X^T\y\right)+\boldsymbol{\beta}^T\left(\tau^2\Bl\right)^{-1}\boldsymbol{\beta}}{2\sigma^2}\right]}\nonumber \\
&\times& \left(\sigma^2\right)^{-\left(a+\frac{n+p}{2}\right)-1}\exp{\left[-\frac{\y^T\left(I-\X M^{-1}\X^T\right)\y+\y^T\y+2b}{2\sigma^2}\right]} \nonumber \\
&\times& \prod_{j=1}^p\left\{\left(1+\frac{1}{\lambda_j^2}\right)^{-1}\left(\lambda_j^2\right)^{-2}\right\}\left(1+\tau^2\right)^{-p/2}\left(\tau^2\right)^{-3p/2}\pi_\tau\left(\tau^2\right)d\boldsymbol{\beta}d\sigma^2d\tau^2\nonumber \\
\end{eqnarray*}

\noindent
Now recall that $M=\omega^*\left(1+\frac{1}{\tau^2}\right)I_p$. Hence 
\begin{eqnarray*}
& \left(\boldsymbol{\beta}-M^{-1}\X^T\y\right)^T M\left(\boldsymbol{\beta}-M^{-1}\X^T\y\right)+\boldsymbol{\beta}^T\left(\tau^2\Bl\right)^{-1}\boldsymbol{\beta} \\
&= \boldsymbol{\beta}^T\left(M+\left(\tau^2\Bl\right)^{-1}\right)\boldsymbol{\beta}-2\boldsymbol{\beta}^T\X^T\y+\y^T\X M^{-1}\X^T\y\nonumber \\
&\leq \boldsymbol{\beta}^T Q\boldsymbol{\beta}-2\boldsymbol{\beta}^T\X^T\y+\y^T\X M^{-1}\X^T\y\nonumber \\
&= \left(\boldsymbol{\beta}-Q^{-1}\X^T\y\right)^T Q\left(\boldsymbol{\beta}-Q^{-1}\X^T\y\right)+\y^T\X\left(M^{-1}- Q^{-1}\right)\X^T\y\nonumber \\
\end{eqnarray*}

\noindent

where $Q=\omega^*\left(1+\frac{1}{\tau^2}\right)\left(I_p+\Bl^{-1}\right).$ Hence it follows that

\begin{eqnarray*}
k\left(\bl_0,\bl\right) &\geq&
\left(\omega_*\right)^{-p/2}d^{-p/\delta_0}\left(1+d^{\frac{2}{\delta_1}}\right)^{-p}\left(\frac{b}{\y^T\y/2+b}\right)^{a+\frac{n}{2}}\frac{b^{a+\frac{n}{2}}}{\Gamma\left(a+\frac{n}{2}\right)} \nonumber \\
&\times& \int_{[T,\infty)}\int_{\mathbb{R}_+}\int_{\mathbb{R}^p}\prod_{j=1}^p\left\{\frac{\beta_j^2}{2\sigma^2}\right\}\left(2\pi\sigma^2\right)^{-\frac{p}{2}}\mid Q\mid ^{1/2}\exp{\left[-\frac{\left(\boldsymbol{\beta}-Q^{-1}\X^T\y\right)^T Q\left(\boldsymbol{\beta}-Q^{-1}\X^T\y\right)}{2\sigma^2}\right]}\nonumber \\
&\times& \mid Q\mid ^{-1/2}\left(\sigma^2\right)^{-\left(a+\frac{n}{2}\right)-1}\exp{\left[-\frac{\y^T\left(I-\X Q^{-1}\X^T\right)\y+\y^T\y+2b}{2\sigma^2}\right]} \nonumber \\&\times& \prod_{j=1}^p\left\{\left(1+\frac{1}{\lambda_j^2}\right)^{-1}\left(\lambda_j^2\right)^{-2}\right\}\left(1+\tau^2\right)^{-p/2}\left(\tau^2\right)^{-3p/2}\pi_\tau\left(\tau^2\right)d\boldsymbol{\beta}d\sigma^2d\tau^2\nonumber \\
\end{eqnarray*}

\noindent
Note that if $\bb\sim \mathcal{N}\left(Q^{-1}\X^T\y,\sigma^2Q^{-1}\right)$ then the inner most integral wrt $\bb$ is equal to 
\begin{eqnarray*}
    \E\left[\prod\limits_{j=1}^p\left\{\frac{\beta_j^2}{2\sigma^2}\right\}\right] &=& \left(2\sigma^2\right)^{-p}\prod\limits_{j=1}^p\left\{ E\left[\beta_j^2\right]\right\};\ \text{since $Q$ is a diagonal matrix, $\beta_j'$s are indep.}\nonumber \\
    &\geq& \left(2\sigma^2\right)^{-p}\prod\limits_{j=1}^p\left\{ \text{Var}\left[\beta_j\right]\right\} \nonumber \\
    &=& \left(2\omega^*\right)^{-p}\left(1+\frac{1}{\tau^2}\right)^{-p}\prod\limits_{j=1}^p\left\{\left(1+\frac{1}{\lambda_j^2}\right)^{-1}\right\} \nonumber \\
\end{eqnarray*} 

\noindent
Also, noting that $\mid Q\mid=(\omega^*)^p\left(1+\frac{1}{\tau^2}\right)^p\prod\limits_{j=1}^p\left(1+\frac{1}{\lambda_j^2}\right)$ we have the following lower bound:

\begin{eqnarray*}
k\left(\bl_0,\bl\right) &\geq& 2^{-p}\left(\omega_*\right)^{-p/2}\left(\omega^*\right)^{-3p/2}d^{-p/\delta_0}\left(1+d^{\frac{2}{\delta_1}}\right)^{-p}\left(\frac{b}{\y^T\y/2+b}\right)^{a+\frac{n}{2}}\frac{b^{a+\frac{n}{2}}}{\Gamma\left(a+\frac{n}{2}\right)} \nonumber \\
&\times& \int_{[T,\infty)}\int_{\mathbb{R}_+}
\left(\sigma^2\right)^{-\left(a+\frac{n}{2}\right)-1}\exp{\left[-\frac{\y^T\left(I-\X Q^{-1}\X^T\right)\y+\y^T\y+2b}{2\sigma^2}\right]} \nonumber \\
&\times& \prod_{j=1}^p\left\{\left(1+\frac{1}{\lambda_j^2}\right)^{-5/2}\left(\lambda_j^2\right)^{-2}\right\}\left(1+\tau^2\right)^{-2p}\pi_\tau\left(\tau^2\right)d\sigma^2d\tau^2\nonumber \\
\end{eqnarray*}

\noindent
Further noting that $\y^T\left(I-\X Q^{-1}\X^T\right)\y\leq \y^T\y$ we have:
\begin{eqnarray*}
k\left(\bl_0,\bl\right) &\geq& 2^{-p}\left(\omega_*\right)^{-p/2}\left(\omega^*\right)^{-3p/2}d^{-p/\delta_0}\left(1+d^{\frac{2}{\delta_1}}\right)^{-p}\left(\frac{b}{\y^T\y/2+b}\right)^{a+\frac{n}{2}}\frac{b^{a+\frac{n}{2}}}{\Gamma\left(a+\frac{n}{2}\right)} \nonumber \\
&\times& \int_{[T,\infty)}\int_{\mathbb{R}_+}
\left(\sigma^2\right)^{-\left(a+\frac{n}{2}\right)-1}\exp{\left[-\frac{1}{\sigma^2}(\y^T\y+b)\right]} \nonumber \\
&\times& \prod_{j=1}^p\left\{\left(1+\frac{1}{\lambda_j^2}\right)^{-5/2}\left(\lambda_j^2\right)^{-2}\right\}\left(1+\tau^2\right)^{-2p} d\sigma^2d\tau^2\nonumber \\
\end{eqnarray*}

\noindent
Integrating wrt $\sigma^2$ we have:
\begin{eqnarray*}
k\left(\bl_0,\bl\right) &\geq& 2^{-p}\left(\omega_*\right)^{-p/2}\left(\omega^*\right)^{-3p/2}d^{-p/\delta_0}\left(1+d^{\frac{2}{\delta_1}}\right)^{-p}\left(\frac{b}{\y^T\y+b}\right)^{2a+n} \nonumber \\
&\times&\prod_{j=1}^p\left\{\frac{\sqrt{\lambda_j^2}}{\left(1+\lambda_j^2\right)^{5/2}}\right\}\int_T^{\infty}\left(1+\tau^2\right)^{-2p}\pi_\tau\left(\tau^2\right)d\tau^2\nonumber \\
&=& \epsilon^* h\left(\bl\right) \nonumber \\
\end{eqnarray*}

\noindent

where 
$$\epsilon^*= 3^{-p}\left(\omega_*\right)^{-p/2}\left(\omega^*\right)^{-3p/2}d^{-p/\delta_0}\left(1+d^{\frac{2}{\delta_1}}\right)^{-p}\left(\frac{b}{\y^T\y+b}\right)^{2a+n}\int_T^{\infty}\left(1+\tau^2\right)^{-2p}\pi_\tau\left(\tau^2\right) d\tau^2$$
and $h$ is a probability density on $\mathbb{R}^p_+$ given by
$$h\left(\bl\right)=\prod_{j=1}^p\left\{\frac{3}{2}\cdot\frac{\sqrt{\lambda_j^2}}{\left(1+\lambda_j^2\right)^{5/2}}\cdot I_{(0,\infty)}\left(\lambda_j^2\right)\right\}$$
Hence, the minorization condition for the MTD \eqref{MTD} is established. \qed

\section{Samplers for conditional posterior distributions of $\bl$ and $\tau^2$ for $K_{aug,reg}$} \label{Appendix D}

\subsection{Rejection sampler for $\bl$}

\noindent
Recall that the target distribution $g\left(\left.\cdot\right|\nu,\beta,\sigma^2,\tau^2,\y\right)$ has density proportion to the function $\phi(\cdot)$ where $$\phi(x)=\left(\frac{1}{c^2}+\frac{1}{\tau^2 x}\right)^{\frac{1}{2}}x^{-\frac{3}{2}}\exp{\left[-\frac{1}{x}\left(\frac{1}{\nu}+\frac{\beta^2}{2\sigma^2\tau^2}\right)\right]}$$
Consider a probability density function $\psi$ on $\mathbb{R}_+$ as follows:
\begin{small}
$$\psi(x)=\frac{\mid c\mid^{-1}\sqrt{\frac{1}{\nu}+\frac{\beta^2}{2\sigma^2\tau^2}}}{\sqrt{\pi}\left(\mid c\mid^{-1}+\left(\tau^2\right)^{-1/2}\right)} x^{-3/2}\exp{\left[-\frac{1}{x}\left(\frac{1}{\nu}+\frac{\beta^2}{2\sigma^2\tau^2}\right)\right]}+\frac{\left(\tau^2\right)^{-1/2}\left(\frac{1}{\nu}+\frac{\beta^2}{2\sigma^2\tau^2}\right)}{\left(\mid c\mid ^{-1}+\left(\tau^2\right)^{-1/2}\right)}x^{-2}\exp{\left[-\frac{1}{x}\left(\frac{1}{\nu}+\frac{\beta^2}{2\sigma^2\tau^2}\right)\right]}$$
\end{small}
Note that the above is a convex combination of two Inverse-Gamma densities and is easy to sample from. 
After simple algebraic manipulation, one can show that
$$\sup_{x\in (0,\infty)}\frac{\phi(x)}{\psi(x)}\leq M,$$
where $$M=\sqrt{\pi}\frac{\mid c\mid^{-1}+\left(\tau^2\right)^{-1/2}}{\sqrt{\frac{1}{\nu}+\frac{\beta^2}{2\sigma^2\tau^2}}}+\frac{\mid c\mid^{-1}+\left(\tau^2\right)^{-1/2}}{\frac{1}{\nu}+\frac{\beta^2}{2\sigma^2\tau^2}}$$
We apply the following algorithm: \\
For $i=1,2,\cdots$
\begin{enumerate}
    \item sample $X_i$ from $\psi(\cdot)$
    \item sample $U_i$ from the uniform distribution over $(0,1)$
    \item Accept $X_i$ if $$U_i\leq \frac{\phi\left(X_i\right)}{M\psi\left(X_i\right)}$$ for all $i$; otherwise, we reject $X_i$.
    \item Repeat the above three steps until we reach a sample of a predetermined size, say, $p$.
\end{enumerate}
\subsection{Metropolis sampler for $\tau^2$}
Recall that the target distribution $\pi\left(\left.\cdot\right|\bl,\y\right)$ has density proportion to the function $\phi(\cdot)$ where 
$$\phi(x)=\left|A_c\right|^{-\frac{1}{2}}\prod\limits_{j=1}^p\left\{\left(\frac{1}{c^2}+\frac{1}{x\lambda_j^2}\right)^{\frac{1}{2}}\right\}\left(\frac{\y^T\left(I_n-\X A_c^{-1}\X^T\right)\y}{2}+b\right)^{-\left(a+\frac{n}{2}\right)}\pi_\tau \left(x\right),$$ 

$A_c=\X^T\X+\left(x\Bl\right)^{-1}+c^{-2}I_p$ and $\pi_\tau(\cdot)$ is a probability density function supported on $\mathbb{R}_+$.
We will also need to pick what is called a “proposal distribution” that changes location at each iteration in the algorithm. We will call this $q\left(u\mid x\right)$. 
Then the algorithm is:
\begin{enumerate}
    \item Choose some initial value $x_0$.
    \item For $i=1,\cdots,p$
    \begin{enumerate}
        \item sample $x_i^*$ from $q\left(u\mid x_{i-1}\right)$.
        \item Set $x_i=x_i^*$ with probalility
        $$\alpha=\text{min}\left(\frac{\phi(x_i^*) q\left(x_{i-1}\mid x_i^*\right)}{\phi(x_{i-1}) q\left(x_i^*\mid x_{i-1}\right)},1\right)$$
        otherwise set $x_i=x_{i-1}$.
    \end{enumerate}
\end{enumerate}
Often times we choose $q$ to be a $\mathcal{N}(x,1)$ distribution. This has the convenient property of symmetry. Which means that $q\left(u\mid x\right)=q\left(x\mid u\right)$, so the quantity $\alpha$ can be simplified to 
$$\alpha=\text{min}\left(\frac{\phi(x_i^*)}{\phi(x_{i-1})},1\right)$$
which is much easier to calculate. This variant is a called a \emph{Metropolis sampler}.

\bibliographystyle{rusnat}
\bibliography{ref_HS_GE.bib}
\end{document}